\documentclass[11pt]{amsart}
\input{mystyle.sty}

\title[Model Free Mean Field RL]{Model-Free Mean-Field Reinforcement Learning: \\ Mean-Field MDP and Mean-Field Q-Learning}
\author{Ren\'e Carmona, Mathieu Lauri\`ere \& Zongjun Tan}
\thanks{This work has been supported by NSF grant DMS-1716673 and ARO grant W911NF-17-1-0578.}
\address{Program in Applied and Computational Mathematics \& ORFE}

\date{}

\begin{document}

\maketitle

\begin{abstract}
We study infinite horizon discounted Mean Field Control (MFC) problems with common noise through the lens of Mean Field Markov Decision Processes (MFMDP). We allow the agents to use actions that are randomized not only at the individual level but also at the level of the population. This common randomization allows us to establish connections between both closed-loop and open-loop policies for MFC and Markov policies for the MFMDP. In particular, we show that there exists an optimal closed-loop policy for the original MFC. 
Building on this framework and the notion of state-action value function, we then propose reinforcement learning (RL) methods for such problems, by adapting existing tabular and deep RL methods to the mean-field setting. The main difficulty is the treatment of the population state, which is an input of the policy and the value function. We provide convergence guarantees for tabular algorithms based on discretizations of the simplex. Neural network based algorithms are more suitable for continuous spaces and allow us to avoid discretizing the mean field state space. Numerical examples are provided.
\end{abstract}



\tableofcontents

\section{\textbf{Introduction}}

\vskip 6pt

In today's highly connected world, important applications often involve very large number of interacting rational agents. Understanding how individual decisions aggregate to create global outcome and how, in turn, the agents react to those outcomes to adjust their behavior is a major challenge for numerous applications. Theoretical analyses distinguish between competitive and cooperative scenarios using game-theoretic notions. From a computational viewpoint, solving games with multiple players becomes infeasible when their number grows, in no small part because the number of pairwise interactions increases exponentially. To cope with this issue, mean field games (MFG) and mean field control (MFC) problems, also called McKean-Vlasov (MKV) control have been introduced (see e.g.~\cite{MR2295621,MR2346927,MR3134900,CarmonaDelarue_book_I,CarmonaDelarue_book_II}). Assuming that the population is homogeneous (the agents have the same transition and cost functions), and that the interactions are symmetric (these functions depend only on the empirical distribution of the other agents), the main idea is to use a mean-field approximation of the population's state in order to simplify the model. It is then sufficient to study the interactions between one representative player and the population distribution rather than the interactions between every pair of players. In the past decade, theoretical results and potential applications have received a growing level of interest. Numerical methods, which are a crucial tool for applications, have also been developed, mostly based on deterministic methods for partial differential equations (see e.g.~\cite{MR2679575,MR2888257,MR3575615}). Computational methods based on deep learning have recently been introduced and seem particularly suitable to tackle high-dimensional MFG and MFC problems (see e.g.~\cite{carmonalauriere2021convergenceI,carmonalauriere2019convergenceII,alaradi2019applications,germain2019numericalMKV,fouquezhang2020deep,ruthotto2020machineMFG,agram2020deepMFC}).   

\vskip 6pt

In the past few years, the question of learning solutions of mean field problems in a model-free way has gained momentum (see e.g.,~\cite{SubramanianMahajan-2018-RLstatioMFG,guo2019learningMFG,fu2019actorcriticMFG,elie2020convergence,gu2019dynamicmfc,gu2020meanQ,anahtarci2020qregu,perrin2020fictitious,motte2019mean,carmona2020policyCDC}). Roughly speaking, the main point is to develop computational methods that can compute MFG or MFC solutions only by sampling realizations of trajectories and without having access to the model (i.e., the transition and cost functions). Model-free methods for a single agent control problem have been developed in the framework of reinforcement learning (RL), building on the formalism of Markov Decision Processes (MDP). One of the most well-known methods is the so-called Q-learning which exploits the dynamic programming principle satisfied by the state-action value function (also called Q-function). Multi-agent reinforcement learning (MARL) extends RL methods to situations in which several agents are simultaneously learning. Breakthrough results have been obtained in games with a small number of players (such as chess or go).

\vskip 6pt

In this work, we focus on a setting with an infinite cooperative population modeled by an MFC problem in discrete time with an infinite horizon discounted cost. Not only is the dynamics is subject to idiosyncratic and common noise, but the actions too are subject to both idiosyncratic and common randomness. The problem can be interpreted as one posed to a central planner who helps a very large population of agents to minimize its social cost and, to this end, can first sample a random policy for the whole population before letting each agent sample their action based on this common policy. This is a distinctive feature compared with the existing literature, in particular~\cite{gu2019dynamicmfc,gu2020meanQ} which has studied RL for MFC without common noise and~\cite{motte2019mean} which has studied MFMDP in the presence of common noise but without common randomization. We use this extra source of randomness to connect open-loop and closed-loop policies for the MFC problem to a mean field MDP (MFMDP) whose state is the population distribution. Defining properly this MFMDP and dealing rigorously with common noise and common policy randomization leads us to carry out a careful probabilistic analysis of this type of problems. Besides the development of the theoretical framework, we also investigate how RL methods can be adapted to the MFC setting, using tabular methods or neural network based methods. Since policies are common to the whole population in the MFMDP, randomization used in RL methods translates into common randomization from the point of view of the MFC problem. We illustrate the performance of two methods on several numerical examples.

\vskip 6pt

The main contributions are threefold.  
{\bf (1)} We introduce a MFMDP and study its connection with the original MFC problem: we prove a DPP for the MFMDP value function (Theorem~\ref{thm:MFMDP_DPP_BOREL}), on which we build to enables to prove equality of the open-loop and closed-loop value functions for the MFC problem (Theorem~\ref{thm:equality_value_function_open_closed_Markov}). Furthermore, we show existence of a stationary closed-loop policy (Proposition~\ref{proposition:existence_opt_policy_mean_field}). 
{\bf (2)} We study the state-action value function (or Q-function) of the MFMDP, for which we prove a DPP (Theorem~\ref{prop:opt_Bellman_eq_Q}). {\bf (3)} We propose several RL methods: a tabular Q-learning relying on a discretization of the mean-field state simplex (Theorem~\ref{thm:main-cv-tabular}), and a deep RL method to deal with continuous state or action spaces, which allows us to avoid simplex discretization or to deal with randomized actions. 

The rest of the paper is organized as follows. Section~\ref{sec:model_description} introduces the main concepts for the MFC problem, including the probabilistic framework and the notions of open-loop and closed-loop policies. Then, the corresponding MFMDP and its DPP are given in Section~\ref{sec:MFMDP}. The connections between the MFC and the MFMDP are developed in Section~\ref{sec:relations-models}. We then turn our attention to the numerical aspects. In Section~\ref{se:Q_learning}, we introduce the state-action value function, prove it satisfies a DPP, and propose computational methods based on RL. Several numerical examples are provided in Section~\ref{sec:numres} to illustrate original features of the MFMDP with common noise and randomized actions.

\section{\textbf{Model Description and Notations \label{sec:model_description}} }

Throughout the paper we work with Borel spaces, namely spaces homeomorphic to a non-empty Borel subset of some Polish space. If $C$ is such a space, we denote by $\cB_C$ its Borel $\sigma$-field and $\cP(C)$ the space of probability measures on $(C,\cB_C)$ implicitly assumed to be equipped with the topology of the weak convergence and its corresponding Borel $\sigma$-field $\cB_{\cP(C)}$.  
For all the measurability issues we refer the reader to any of the textbooks \cite{BertsekasShreve} or \cite{Kallenberg_RM}.

\subsection{Probabilistic set-up of the MFC model}\,

In this section, we specify what we mean by mean-field models with common noise. We first introduce the major building blocks, leaving the description of the dynamics for later on.

\begin{definition}[MFC model]
	\label{def:MFC_problem}
	An infinite horizon discounted \defi{mean-field control (MFC) model} with common noise is based on the following elements $(S, A, E, E^0, F, f, \gamma)$:
	\begin{itemize}
		\item A Borel space $(S, \cB_S)$ for the state space.
		\item A Borel space $(A, \cB_A)$ for the action space. 
		\item Two Borel spaces $(E, \cB_E)$ and $(E^0, \cB_{E^0})$ for the values of the idiosyncratic and common noise.
		\item A Borel measurable function $F: S \times A \times \cP(S \times A) \times E \times E^0 \to S$ called the system function.
		\item A bounded Borel measurable function $f: S \times A \times \cP(S \times A) \to \RR$ called the one-stage cost function.
		\item A discount factor $\gamma \in (0,1)$.
	\end{itemize}
\end{definition}
The system function $F$ is used to describe the evolution of the state process based on a state, an action, a mean-field interaction term\footnote{The interactions are through the joint state-action distribution, which is sometimes referred to as ``extended'' MFC or MFC of controls.} and two noise terms. 

Even if we are willing to postpone the regularity conditions on $F$ and $f$, the above definition is still incomplete. It introduces the building blocks of the model, but does not explain how the actions are taken and how the system evolves over time.  In order to motivate the nature of the assumptions we are about to introduce, we take an informal excursion in the world of finitely many actors.
Our goal is to motivate the following important features of our model:
\begin{enumerate}[label=(\roman*)]
	\item The mean-field interactions are conditioned on all shared information.
	\item The actions are randomized  by additional sources of randomness.
\end{enumerate} 

\subsubsection{Motivation from Finitely Many Player Models}

Let us imagine that $N$ robots with states $X^1_n,\cdots, X^N_n$, take actions  $\alpha^1_n,\cdots, \alpha^N_n$ at time $n$ and that
their next state is given by:
$$
    X^i_{n+1}=F\Bigl(X^i_n,\alpha^i_n,\frac1N\sum_{j=1}^N\delta_{(X^i_n,\alpha^i_n)},\epsilon^i_{n+1}\Bigr),\quad n\ge 0,\;\;i=1,\cdots,N
$$
$\epsilon^1_n,\cdots, \epsilon^N_n,\cdots$ being independent and identically distributed (i.i.d. from now on)  random shocks with distribution $\nu\in\cP(E)$.
Let us also assume that at each time $n$, a central unit collects the information about the states, the actions and the costs (or rewards) incurred by the individual robots, and minimizes the overall average cost to the system as given by:
$$
    \frac1N\sum_{j=1}^N \EE\Bigl[\sum_{n=0}^\infty \gamma^n f\bigl(X^j_n,\alpha^j_n,\frac1N\sum_{i=1}^n\delta_{(X^i_n,\alpha^i_n)}\bigr)\Bigr].
$$
In the limit $N\to\infty$, using standard propagation of chaos arguments, we expect that the coupled evolutions of the states of the individual robots will become independent, and that each state evolves according to the dynamics:
$$
    X_{n+1}=F(X_n,\alpha_n,\PP_{(X_n,\alpha_n)},\epsilon_{n+1}),\quad n\ge 0,
$$
and one can then imagine that the optimization of the central unit reduces to the minimization:
$$
    \inf_{\balpha=(\alpha_n)_{n\ge 0}} \EE\Bigl[\sum_{n=0}^\infty \gamma^n f\bigl(X_n,\alpha_n,\PP_{(X_n,\alpha_n)}\bigr)\Bigr],
$$
where $\PP_{(X_n,\alpha_n)}$ denotes the joint law of the state-action couple at time $n$. This is the formulation of (discrete time) MFC problem, which we will not call and MDP because of the lack of Markov property due to the presence of $\PP_{(X_n,\alpha_n)}$ in the equation, responsible for the McKean-Vlasov nature of the dynamics.

Still, such a model does not account for the fact that the robots evolve in an environment which is most likely random. Since all the robots face the same random shocks due to the randomness of the common environment, we introduce a sequence
$(\epsilon^0_n)_{n\ge 1}$ of i.i.d. random elements in $E^0$, with common law $\nu^0\in\cP(E^0)$,  independent of the idiosyncratic shocks $(\epsilon^i_n)_{n\ge 1,i=1,\cdots,N}$ of the individual robots, and with this addition to the model, the individual robot state dynamics
become:
$$
X^i_{n+1}=F\Bigl(X^i_n,\alpha^i_n,\frac1N\sum_{j=1}^N\delta_{(X^i_n,\alpha^i_n)},\epsilon^i_{n+1},\epsilon^0_{n+1}\Bigr),\quad n\ge 0,\;\;i=1,\cdots,N.
$$
Even though we added significant sources of coupling between the robots experiences, exchangeability among the robots persists, and in the limit $N\to\infty$, 
a conditional form of the propagation of chaos theory should lead to conditional independence of  individual dynamics which should take the generic form:
$$
X_{n+1}=F(X_n,\alpha_n,\PP^0_{(X_n,\alpha_n)},\epsilon_{n+1},\epsilon^0_{n+1}),\quad n\ge 0,
$$
$\PP^0_{(X_n,\alpha_n)}$ being the conditional law of the state-action couple $(X_n,\alpha_n)$ given the common noise  $(\epsilon^0_n)_{n\ge 1}$. Accordingly, the optimization of the central unit should be:
$$
\inf_{\balpha=(\alpha_n)_{n\ge 0}} \EE\Bigl[\sum_{n=0}^\infty \gamma^n f\bigl(X_n,\alpha_n,\PP^0_{(X_n,\alpha_n)}\bigr)\Bigr].
$$
This form of Mean Field discrete time Control problem with common noise is now very close to the model we investigate in this paper.
The last question we would like to address before turning to the theoretical analysis of the model is:  \emph{Could the conditional distributions $\PP^0_{(X_n,\alpha_n)}$ depend upon some extra sources of randomness?}
Indeed, if individual robots and the central unit are allowed to use mixed strategies, they need independent sources of randomness to randomize their actions and decisions at each time. So since the cost optimization and the choice of an action strategy are performed by the central unit,
\begin{itemize}
\item the central unit needs a source of randomness to \emph{randomize} the choice of a mixed policy to be dispatched to the individual robots, and
\item the individual robots need to sample their actions from the policy sampled for them by the central unit.
\end{itemize}
So in the limit $N\to\infty$, the \emph{idiosyncratic}  randomizations of the individual robots average out, while the \emph{central unit randomization} remains.

As a result, the conditioning in $\PP^0_{(X_n,\alpha_n)}$ should be with respect to the common noise $\epsilon^0_1,\cdots,\epsilon^0_n$ \emph{as well as} the sources of randomness used by the central unit to randomize the policy handed out to the individual robots for implementation.

\subsubsection{Back to the MFC problem formulation}
\label{sec:back-to-MFC}
Motivated by the previous $N$-agent model discussion when the size of the population $N$ tends to infinity, we propose the following set-up to disentangle clearly the various sources of randomness.

\vskip 2pt
We assume that all the sources of randomness are from a probability space $(\Omega,\cF,\PP)$ supporting
\begin{enumerate}[label=(\roman*)]
	\item  an i.i.d. sequence $(\epsilon_n)_{n\ge 0}$ with distribution $\nu \in \cP(E)$  modeling the idiosyncratic random shocks;
	\item an i.i.d. sequence $(\epsilon^0_n)_{n\ge 0}$ with distribution $\nu^0 \in \cP(E^0)$ modeling the common noise;
	\item a random variable $\mathscr{U}$ with distribution $\PP_{\mathscr{U}} $ in a Borel space $(\Upsilon, \cB_{\Upsilon})$ providing the randomization for the initial state;
	\item  an i.i.d. sequence $(\vartheta_n)_{n\ge 0}$ of random variables in a Borel space  $(\Theta, \cB_{\Theta})$ with distribution $\PP_{\vartheta}$ providing the generic robot with a source of randomization for their action choices;
	\item  an i.i.d. sequence $(\vartheta^0_n)_{n\ge 0}$ of  random variables in a Borel space  $(\Theta^0, \cB_{\Theta^0})$ with distribution $\PP_{\vartheta^0}$ providing the central unit with a source of randomization for the choices of policies.
\end{enumerate} 
We assume that all these random sequences are independent of each other. We also assume that  $\PP_{\vartheta}$ and $\PP_{\vartheta^0}$ are both atomless. This guarantees the existence of Borel measurable functions $h:\Theta\mapsto[0,1]$ and $h^0:\Theta^0\mapsto[0,1]$ which are uniformly distributed when viewed as random variables on the probability spaces $(\Theta,\cB_\Theta,\PP_\vartheta)$ and $(\Theta^0,\cB_{\Theta^0},\PP_{\vartheta^0})$  respectively.
The uniform random variables constructed with the functions $h$ and $h^0$ will be used repeatedly with the following classical result from measure theory which we state for the sake of later reference. 
\begin{lemma}[Blackwell-Dubins Lemma]
\label{le:BlackwellDubins}
For any Polish space $B$, there exists a measurable function $\rho_B:\cP(B)\times [0,1]\mapsto B$, which we shall call the Blackwell-Dubins function of the space $B$, satisfying
\vskip 0pt
i) for each $\nu\in\cP(B)$ and each uniform random variable $U\sim U(0,1)$, the $B$ - valued random variable $\rho_B(\nu,U)$ has distribution $\nu$;
\vskip 0pt
ii) for almost every $u\in[0,1]$, the function $\nu\mapsto \rho_B(\nu,u)$ is continuous for the weak topology of $\cP(B)$.
\end{lemma}

\vskip 6pt
For the sake of illustration, some of the computations will be done in the canonical probability space $(\Omega^c, \cF^c, \PP^c)$ defined by: 
	\begin{equation*}
		\Omega^c =  \Upsilon \times \Theta \times \Theta^0 \times (E \times E^0 \times \Theta \times \Theta^0)^\infty, 
\qquad
		 \cF^c= \cB_{\Upsilon} \times \cB_{\Theta} \times \cB_{\Theta^0} \times (\cB_E \times \cB_{E^0} \times \cB_{\Theta} \times \cB_{\Theta^0})^\infty, 
	\end{equation*} 
	and the product probability $\PP^c$ given by:
	\begin{equation*}
		\PP^c = \PP_{\sU} \otimes \PP_{\vartheta} \otimes \PP_{\vartheta^0} \otimes \  ( \nu \otimes \nu^0 \otimes \PP_{\vartheta} \otimes \PP_{\vartheta^0} )^\infty.
	\end{equation*}
A generic element $\omega \in \Omega$ reads
	$$
		\omega = ( u, \theta_0, \theta_0^0, e_1, e_1^0, \theta_1, \theta_1^0, e_2, e_2^0 \ldots,  e_{n}, e_{n}^0, \theta_n, \theta_n^0,  \ldots),
	$$
and we realize the random variables as coordinate mappings, $\sU(\omega) = u$, and for $n \geq 0$, 
	\begin{equation*}
		\quad \vartheta_n(\omega) = \theta_n, 
		\quad \vartheta_n^0(\omega) = \theta_n^0, 
		\quad \varepsilon_{n+1} (\omega) = e_{n+1}, 
		\quad \varepsilon_{n+1}^0(\omega) = e_{n+1}^0 .	
	\end{equation*}
\vskip 6pt
 We will use the short hand notations $\underline \vartheta_n = (\vartheta_0, \ldots, \vartheta_n)$, $\underline \vartheta_n^0 = (\vartheta_0^0, \ldots, \vartheta_n^0)$ for every $n \geq 0$, and $\underline \varepsilon_n = (\varepsilon_1, \ldots, \varepsilon_n)$, $\underline \varepsilon^0_n = (\varepsilon^0_1, \ldots, \varepsilon_n^0)$ for every $n \geq 1$.

 \subsection{Probabilistic framework and classes of policies}
\subsubsection{Filtrations, action processes, and control processes}\,
\label{subsec:filtractions-action-controls}
We now introduce the framework that will be used to rigorously study the MFC problem. As explained intuitively in the introduction, we will distinguish several types of randomness. We will also distinguish between actions (elements of $A$), and controls (probability measures on $A$). These are the building blocks to define later the notion of policy (see \S~\ref{subsec:open-closed-policies}).

The $\sigma-$field for the initial state is denoted by $\cF_{x_0} = \sigma\{\sU\}$. 
We introduce four filtrations to define the action and control processes.
The filtrations of the idiosyncratic and common noises are: 
		$$
			\cF_0^{\idioNoise} = \cF_0^{\commNoise} = \{ \emptyset, \Omega \}, \quad \cF_n^{\idioNoise} = \sigma\{ \underline \varepsilon_n\}, \quad \cF_n^{\commNoise} = \sigma\{ \underline \varepsilon^0_n \}, \quad \  n \geq 1.
		$$
The filtration of the (idiosyncratic) action randomization is:
		$$
			\cF_n^\Theta = \sigma\{ \Theta_0, \ldots, \Theta_n \}, \qquad \  n \geq 0.
		$$
The filtration of the (common) policy  randomization is:
		$$
			\cF_n^{\Theta^0} = \sigma\{\Theta_0^0, \ldots, \Theta_n^0 \}, \qquad \  n \geq 0.
		$$

We also introduce three new  filtrations $\FF^0 = (\cF_n^0)_{n \geq 0}$, $\GG^c = (\cG_n^c)_{n \geq 0}$ and $\GG^a = (\cG_n^a)_{n \geq 0}$ defined by:
\begin{equation*}
\arraycolsep=1pt\def\arraystretch{1.5}
		\begin{array}{ll}
			\cF_0^0 = \sigma\{ \vartheta_0^0 \} & \quad \cF_n^0 = \cF_n^\commNoise \vee \cF_n^{\Theta^0},
			= \sigma\{ \underline{\vartheta}_n^0, \underline{\varepsilon}_n^0 \}, 
			\quad n \geq 1,
			\\
			\cG_0^{c} = \sigma\{ \sU, \vartheta_0^0 \}, & \quad  \cG_n^c = \cF_{x_0} \vee \cF_n^{\idioNoise} \vee \cF_n^{\commNoise} \vee \cF_n^{\Theta^0} \vee \cF_{n-1}^{\Theta},
			= \sigma\{ \sU, \underline{\vartheta}_{n-1}, \underline{\vartheta}_n^0, \underline{\varepsilon}_{n}, \underline{\varepsilon}_{n}^0 \},   
			\quad n \geq 1,
			\\
			\cG_0^a = \sigma \{ \sU, \vartheta_0^0,  \vartheta_0 \},  & \quad  \cG_n^a = \cF_{x_0} \vee \cF_n^{\idioNoise} \vee \cF_n^{\commNoise} \vee \cF_n^{\Theta^0} \vee \cF_{n}^{\Theta}
			= \sigma\{ \sU, \underline{\vartheta}_{n-1}, \underline{\vartheta}_n^0, \underline{\varepsilon}_{n}, \underline{\varepsilon}_{n}^0, \vartheta_n \},
			\qquad n \geq 1.
		\end{array}
\end{equation*}

\vskip 6pt
Next, we introduce the terminology which is going to help us characterize the information available to a generic robot and the central unit to make their choices of actions and policies. Incidentally, we start referring to the robots as agents and everything related to agents (e.g. actions, controls, policies, costs, etc) will be referred as level-0. This is in contrast with the lifted stochastic optimization model which we introduce in Section \ref{sec:MFMDP} whose elements will be referred to as level-1. For reasons which will become clear later, the central unit controlling the robots will be called the \emph{level-1 controller}. This new terminology frees our presentation of the theoretical results from the gory details of the robotic application we used to motivate the set-up.

\begin{definition}
A \defi{level-0 action} is an element of $A$. A \defi{level-0 (mixed) control} is a random probability measure on $(A, \cB_A)$, that is, any random variable with values in the Borel space $(\cP(A), \cB_{\cP(A)})$. 
\end{definition}
Typically, a level-0 action is denoted by $a$ and a level-0 control is denoted by $\fa$. Unless specified otherwise, all the controls we consider are mixed. So for the sake of brevity, we will omit this term.

We now consider the notion of action and control processes for a representative agent. Intuitively, an action process is the realization of a control process, where the sampling is done using $(\vartheta_n)_{n \ge 0}$. 

\begin{definition}%
 \label{de:level-0}
 A \defi{level-0 action process} is a sequence of random variables $\balpha = (\alpha_n)_{n \geq 0}$ with values in $A$ which is adapted to the filtration $\GG^a$. The set of such action processes is denoted by $\AA$.
A \defi{level-0 control process} is a sequence $\bfa = (\fa_n)_{n \geq 0}$ of level-0 controls which is adapted to the filtration $\GG^c$.
Finally, an  action process $\balpha = (\alpha_n)_{n \geq 0}$ is said to be a \defi{realization} of a level-0 control process $\bfa = (\fa_n)_{n \geq 0}$ if 
			$\cL \big( \alpha_n \, | \, \sigma\{ \sU, \underline{\vartheta}_{n-1}, \underline{\vartheta}_n^0, \underline{\varepsilon}_{n}, \underline{\varepsilon}_{n}^0\} \big)  = \fa_n$, $\PP-a.s.$,  for every $n \geq 0$.	
\end{definition}
Here and in the following, we use the notations $\PP_\xi$ and $\cL(\xi)$ interchangeably for the distribution of a random element $\xi$, and we use natural extensions to denote conditional distributions.

It can be shown (see Lemma~\ref{lemma:identitly_level_0_action_realization} in the appendix) that, for any level-0 control process $\fa$ and any two realizations $\balpha,\balpha'$ of $\fa$, every bounded Borel measurable function $h$, 
$\EE \left[ h( \alpha_n' ) \, | \, \cG_n^c \right] = \int_A h(\alpha) \fa_n(d \alpha) = 	\EE \left[ h( \alpha_n ) \, | \, \cG_n^c \right]$, $\PP$-a.s., $n \ge 0$.

\subsubsection{Conditional distribution and state process}
We are now in a position to describe precisely the mean-field interactions in the system function, and provide a clear definition of the state process driven by a mixed control process in the mean-field model with common noise.

\begin{definition}%
	\label{def:state_process_from_control_process}
For any initial distribution $\mu_0 \in \cP(S)$ and level-0 action process $\balpha = (\alpha_n)_{n \geq 0}$, we say that  a process $\bX^{\balpha,\mu_0} = (X_n^{\balpha, \mu_0})_{n \geq 0}$ is a \defi{state process associated to} $(\balpha, \mu_0)$ for the MFC model if: 
			$X_0^{\balpha, \mu_0} $ is an $S$-valued $\cF_{x_0}$ - random variable with distribution $ \mu_0$,
and for every $n \geq 0$, 
			\begin{equation}
			    \label{eq:system_dynamics_level_0}
				X_{n+1}^{\balpha, \mu_0} = F\big( X_n^{\balpha, \mu_0}, \alpha_n, \PP^0_{(X_n^{\balpha, \mu_0}, \alpha_n)}, \varepsilon_{n+1}, \varepsilon_{n+1}^0 \big),
			\end{equation}
			where $\PP^0_{(X_n^{\balpha, \mu_0}, \alpha_n)}$ is a regular version of $\cL\bigl( (X_n^{\balpha, \mu_0}, \alpha_n) \, | \, \cF_n^0 \bigr)$, the conditional joint distribution of state-action at time $n$ with respect to common noise and common randomization up to current time.
\end{definition}
Such a state process $\bX^{\balpha, \mu_0}$ is adapted to the filtration $\GG^x = (\cG_n^x)_{n \geq 0}$, defined by
\begin{equation}
	\label{eq:filtration_G^x}
	\cG_0^x = \sigma\{ \sU \}, \qquad \cG_n^x = 
	\sigma\{ \sU, (\vartheta_k, \vartheta_k^0, \varepsilon_{k+1}, \varepsilon_{k+1}^0 )_{k = 0, \ldots, n-1} \}, \quad \  n \geq 1.
\end{equation}

For each level-$0$ action process $\balpha$ and each  $n \geq 0$, we denote by $\PP^0_{X_n^{\balpha, \mu_0}}$ a regular version of the conditional distribution $\cL(X_n^{\balpha, \mu_0}\,|\,\cF^0_n)$. It holds:
\begin{equation}
\label{fo:P0X_n}
    \PP^0_{X_n^{\balpha, \mu_0}}
    =\cL(X_n^{\balpha, \mu_0}\,|\,\sigma\{\underline{\epsilon}^0_n,\underline{\vartheta}^0_{n-1}\}),  \qquad \mathbb{P}-a.s..
\end{equation}
this is due to the fact that $X_n^{\balpha, \mu_0}$ is $\cG_n^x$-measurable, hence $X_n^{\balpha, \mu_0} \perp_{ \cF_n^\commNoise \vee \cF_{n-1}^{\commRandom} } \vartheta_n^0$.

\subsubsection{Open-loop and closed-loop policies}
\label{subsec:open-closed-policies}
We now introduce two concepts of policies: open-loop and closed-loop ones.

We first consider open-loop policies. For each $n\ge 0$, let $\Xi_n = \Upsilon \times (\Theta \times \Theta^0 \times E \times E^0)^n \times \Theta^0$, with the convention that $\Xi_0 = \Upsilon \times \Theta^0$, and let us define $\bxi=(\xi_n)_{n\ge 0}$ by: 
$$
    \xi_0 = (\sU, \vartheta_n^0), 
    \qquad\text{and}\qquad \xi_n = (\sU, (\vartheta_k, \vartheta_k^0, \varepsilon_{k+1}, \varepsilon_{k+1}^0)_{k = 0, \ldots, n-1}, \vartheta_n^0), \qquad n \ge 1.
$$

\begin{definition}%
\label{def:open_loop_policy} 
An \defi{level-0 open-loop policy} is a sequence $\bpi= (\pi_n)_{n \geq 0}$ of deterministic measurable functions $\pi_n : \Xi_n \to \cP(A)$, called the \defi{open-loop strategy functions} at time $n$ of policy $\bpi$. 
The set of all open-loop policies is denoted by $\bPi^\tinyol$. 
	A level-$0$  control process $\bm{\fa} = (\fa_n)_{n \ge 0}$ is said to be \defi{generated by} the open-loop policy $\bpi$ if: %
	\begin{equation}
		\label{eq:def_admissible_open_loop_control}
		\fa_n = \pi_n( \xi_n ), \quad \PP-a.s., \qquad n \geq 0.
	\end{equation}
\end{definition}

Notice that because of measurability restrictions ($\fa_n$ must be adapted to the filtration $\GG^c$), there is a one-to-one correspondence between level-0 control processes $\bfa$ and open-loop policies $\bpi$, and both objects can be identified  through equation~\eqref{eq:def_admissible_open_loop_control}. We shall use one object or the other depending on whether we want to emphasize the stochastic process or the sequence of deterministic functions defining the process. 
To be consistent with Definition \ref{de:level-0}, we shall often short-circuit the control process $\bfa$ and say that an action process $\balpha = (\alpha_n)_{n \geq 0}$ is (a realization of the control process) \defi{generated by} $\bpi$ if $\cL( \alpha_n\, | \, \cG_n^c ) = \pi_n(\xi_n),$ $\PP-a.s.$, for all $n \geq 0.$

We now consider closed-loop policies. 

\begin{definition}%
	\label{def:Markovian_policy}
	A \defi{closed-loop Markov strategy function} is a measurable function from $S \times \cP(S) \times \Theta^0$ into $\cP(A)$. 
	A \defi{closed-loop Markov policy} $\bpi=(\pi_n)_{n\ge 0}$  is a sequence of such functions. %
The set of all closed-loop Markov policies is denoted by $\bPi^\tinycl$. 
\end{definition} 

The choice of this form of closed-loop Markov strategy function is suggested by the mean-field nature of the dynamics~\eqref{eq:system_dynamics_level_0} and the form of the costs~\eqref{eq:J_alpha}. Taking values in $\cP(A)$ instead of $A$ indicates that the strategy functions are mixed. Typically, they suggest that at each time $n \geq 0$, the action $\alpha_n \in  A$ taken according to such a policy should be sampled from a probability measure depending directly on the values of $X_n^{\balpha, \mu_0}$ and $\PP^0_{X_n^{\balpha, \mu_0}}$, and the random variable $\vartheta_n^0$ used by the level-$1$ controller to randomize their choice. We formalize this procedure in Definition \ref{def:admissible_state_action_processes} below.

\begin{definition}%
	\label{def:admissible_state_action_processes}
	
	For a closed-loop Markov policy $\bpi\in \bPi^\tinycl$ and an initial distribution $\mu_0 \in \cP(S)$, a pair of state and action processes $(\bX, \balpha) = ( X_n, \alpha_n)_{n \geq 0}$ is said to be \defi{generated by} $(\bpi, \mu_0)$ if 
	\begin{enumerate}[label=\roman*)]
		\item $\bX$ is a state process associated to $(\balpha, \mu_0)$ in the sense of Definition~\ref{def:state_process_from_control_process}. 
		\item The action process $\balpha$ is adapted to $\GG^a$ and satisfies
		\begin{equation}
			\label{eq:definition_action_Markov_closed_loop}
			\cL \big( \alpha_n \, | \,  \cG_n^c \big) = \pi_n \big(  X_n, \, \PP^0_{X_n},\, \vartheta_n^0 \big), \qquad \PP-a.s., \qquad n \geq 0.
		\end{equation}
	\end{enumerate}
\end{definition}

In Definition~\ref{def:admissible_state_action_processes}, we can view the state and action processes as constructed simultaneously by alternatively invoking the system dynamics~\eqref{eq:system_dynamics_level_0} and the sampling procedure consistent with~\eqref{eq:definition_action_Markov_closed_loop}. 
A convenient way to construct an action process $\balpha$ satisfying \eqref{eq:definition_action_Markov_closed_loop} is to use the Blackwell-Dubin's Lemma~ \ref{le:BlackwellDubins}. Indeed, if $\rho_A$ is the Blackwell-Dubin's function of $ A$ and the uniformly distributed random variables $U_n$ is given by $U_n =h(\vartheta_n)$, we can choose $\alpha_n= \rho_A \big( \pi_n  \big(  X_n, \,  \PP^0_{X_n} ,\, \vartheta_n^0 \big), U_n \big)$, $\PP$-a.s.,  $n \geq 0$.

\begin{remark}
Even though we call a policy $\bpi\in \bPi^{\tinycl}$ a ``Markov" policy, 
	it does not imply any Markov property for the state process $\bX$ associated to such a policy. This abuse of terminology can be explained by our intention to work with level-1 Markov policies which will imply the Markov property for a lifted measure-valued state process constructed in the next section.  Also, since we only use the term ``Markov policy'' in the closed-loop setting, we shall most often drop the term closed-loop hereafter and only call them simply Markov policies. 
\end{remark}

\subsection{Optimization and value functions}\,
\label{subsec:optim-valuefct}

We now move to the definition of the optimization problems for an MFC model with open-loop and closed-loop Markov policies.

We now introduce the value function associated to an action process. 
\begin{definition}%
	For every level-0 action process $\balpha$, the \defi{value function}  $J^{\balpha}: \cP(S) \to \RR$ is  defined for every $\mu_0 \in \cP(S)$ by:
	\begin{equation}
		\label{eq:J_alpha}
			J^\balpha(\mu_0) := \EE \left[\sum_{n \geq 0} \gamma^n f\left(X_n^{\balpha, \mu_0}, \alpha_n,\PP^0_{(X_n^{\balpha, \mu_0},\alpha_n)} \right)  \right],
	\end{equation} 
	where the state process $\bX^{\balpha, \mu_0}$ is associated to $(\balpha, \mu_0$) according to the dynamics~\eqref{eq:system_dynamics_level_0}.
\end{definition}

The value of $J^\balpha(\mu)$ is well-defined for every $\mu \in \cP(S)$ because $f$ is measurable and bounded. Furthermore, this value depends only upon the sequence of one-dimensional marginal distributions of the $\cP(S\times A)$-valued process $\Bigl(\PP^0_{(X_n^{\balpha, \mu_0},\alpha_n)} \Bigr)_{n\ge 0}$. %

For any open-loop or closed-loop policy $\bpi$, we can show that $\mathbb{P}^0_{(X_n,\alpha_n)}$ depends on the action process only through the policy $\bpi$ provided $(\bX,\balpha)$ is generated by $\bpi$. See Lemma~\ref{le:policy_value} in the appendix. As a consequence, we can define, for any level-$0$ action process $\balpha$  generated by $\bpi$, 
$$
    J^{\bpi}(\mu) = J^{\balpha}(\mu), \qquad \mu \in \cP(S).
$$

Accordingly, we define the optimal open-loop value function and the optimal closed-loop value function as:
	\begin{equation*}
		J^{\tinyol,*}(\mu) = \inf_{ \bpi \in \bPi^\tinyol }  J^{\bpi}( \mu ), \qquad J^{\tinycl,*}(\mu) = \inf_{ \bpi \in \bPi^\tinycl }  J^{\bpi}( \mu ), \qquad  \mu \in \cP(S),
	\end{equation*}
	which are finite because we assume that the one-stage cost function $f$ is bounded and $\gamma \in (0,1)$.

\section{\textbf{Mean-Field MDP }}
\label{sec:MFMDP} 

In this section, we introduce a Markov Decision Process (MDP) which we use to identify an optimal closed-loop Markov policy for our original MFC model.

\subsection{Mean-field MDP framework}\, 

The key observation is that, for a mixed Markov closed-loop policy $\bpi \in \bPi^\tinycl$, the associated value function $ J^{\bpi}$  can be viewed as the value function of an MDP with state space $\cP(S)$, state process $(\PP^0_{X_n})_{n \geq 0}$ and  action process $\bigl(\PP^0_{(X_n,\alpha_n)}\bigr)_{n \geq 0}$ with values in the new action space $\cP(S \times A)$. Actions need to be consistent with the state in the sense that the first marginal of any action which can be taken while in a given state, has to be equal to the state itself. We provide a rigorous definition of this MDP, and we show its connection to the original MFC model. 

\vskip 6pt
A  \defi{mean-field MDP (MFMDP)} consists of a six-tuple $(\bar S, \bar A, \bar \Gamma, P, \bar f, \gamma)$ as described below:
	\begin{itemize}
		\item The state space is the Borel space $\bar S := \cP(S)$; a generic element in $\bar S$ is denoted by $\mu$.
		\item The action space is the Borel space $\bar A := \cP(S \times A)$; a generic element in $\bar A$ is denoted by $\bar a$.		
		\item The control constraint is a set-valued function ${\bar U}$ from ${\bar S}$ into the set of non-empty subsets of ${\bar A}$ defined by:
		\begin{equation}
			\label{eq:Ubar}
			{\bar U}(\mu) :=\{{\bar a}\in{\bar A};\; \text{pr}_1({\bar a})=\mu\}, \qquad \forall \, \mu\in{\bar S};
		\end{equation}
		where $\pr_1: \bar A \to \bar S$ is the projection function that maps $\bar a \in \bar A$ onto its first marginal distribution on $S$. 
		We shall also use the notation
		\begin{equation}
			\label{eq:bar_Gamma}
			\bar\Gamma :=\{(\mu,{\bar a})\in{\bar S}\times{\bar A};\; {\bar a}\in{\bar U}(\mu)\}.
		\end{equation} 		
		\item The transition probability kernel $P: \bar \Gamma \to \cP( \bar S)$, which is a Borel measurable function. %
		\item The one-stage cost function $\bar f: \bar \Gamma \to \RR$, which is a bounded measurable function. 
		\item The discount coefficient $\gamma \in (0,1)$.
	\end{itemize}

\begin{remark}\,
	The \emph{projection map}  $\pr_1$ is continuous, so the constraint set $\bar U(\mu)$ is closed in $\bar A$ for every $\mu \in \bar S$. The graph $Gr(pr_1) := \{ (\bar a, \mu): pr_1(\bar a) = \mu\} \subset \bar A \times \bar S$ is closed, so $\bar \Gamma$ is also closed in $\bar S \times \bar A$. Hence $\bar\Gamma$ is an analytic subset of $\bar S \times \bar A$, and a Polish space on its own. We assume that $\bar \Gamma$ is endowed with the induced topology as well as the trace $\sigma$-field inherited from $\bar S \times \bar A$.
\end{remark}

\begin{definition}%
	\label{def:MFMDP_problem}
The six-tuple $(\bar S, \bar A, \bar \Gamma, P, \bar f, \gamma)$ 
is said to be the  \defi{MFMDP  lifted from the MFC model} $(S, A, E, E^0, F, f, \gamma)$ of Definition~\ref{def:MFC_problem} if it satisfies:
	\begin{itemize}
		\item The transition kernel $P$ is given by
		\begin{equation}
		\label{def:MFMDP_transition_kernel_from_system_func}
	        P( \mu, \bar a)(d \mu') = \big( \nu^0 \circ \bar F( \mu, \bar a, \cdot)^{-1})(d \mu'), \qquad (\mu, \bar a) \in \bar \Gamma,
	    \end{equation}
		where $\bar F: \bar \Gamma \times E^0 \to \bar S$ is the  system function defined in terms of $F$ by: 
		\begin{equation}
			\label{eq:Fbar_pushforward}
			{\bar F}(\mu, \bar a ,e^0)  =  (\bar a \otimes \nu)\circ F(\cdot, \cdot,\bar a, \cdot,e^0)^{-1}, \qquad (\mu, \bar a, e^0 ) \in \bar \Gamma \times E^0.
		\end{equation}
		\item The one-stage cost function $\bar f: \bar\Gamma \to \RR$ of the MFMDP satisfies:
		\begin{equation}
			\label{def:bar_f}
			\bar f(\mu, \bar a )= \int_{S\times A} f(x, \alpha, \bar a) \bar a(dx,d\alpha), \qquad (\mu, \bar a) \in \bar \Gamma.
		\end{equation}
	\end{itemize}
\end{definition}

Here and in the following we denote by $\nu\circ g^{-1}$ the push-forward of a measure $\nu$ by a measurable function $g$. %
We defined the dynamics using a transition kernel $P$ for two reasons: 1) to conform with the standard literature on MDPs which seems to prefer transition kernels to system functions; 2) we shall restrict ourselves to Markovian policies and control processes given by feedback functions of the state for the analysis of the lifted MDP.  

We can check that $\bar F$ is Borel measurable; see e.g. \cite[Proposition~7.29]{BertsekasShreve} and Remark~\ref{rem:barF-Borel} in the appendix.

\begin{remark}
In anticipation for what is going to come next, we want to emphasize that the
above MFMDP satisfies the assumptions of \cite[Chapter 8-9]{BertsekasShreve}. We will use the results therein to derive a form of Dynamic Programming Principle (DPP) for the MFMDP. %
\end{remark}

To highlight the tight connections between the lifted MFMDP and the original MFC, we define the notion of mixed strategy and mixed Markov policy for MFMDP.

\begin{definition}%
	\label{def:MFMDP_non_mixed_strategy_function}
We call \defi{level-1 pure strategy function} any Borel measurable function from $\bar S$ into $\bar A$ whose graph is contained in $\bar\Gamma$.	
We call \defi{level-1 mixed strategy function} any Borel measurable function $\bar\pi$ from $\bar S$ into $\cP(\bar A)$ satisfying
	$
		\bar \pi(\mu) (\bar U(\mu) )  = 1,$ $\mu \in \bar S.
	$ 
We denote by $\overline{\Pi}^{p}$ (resp.  $\overline{\Pi}$) the set of pure (resp. mixed) strategy functions. 
\end{definition}
We shall identify every $\beta\in\overline{\Pi}^{p}$ with the corresponding level-1 mixed strategy $\bar\pi$ defined by $\bar \pi(\mu) = \delta_{ \beta(\mu) }$ for $\mu \in \bar S$. %
For the sake of brevity, we sometimes omit the term ``mixed" or ``randomized'' when $\bar \pi \in \overline \Pi$, but we always keep the term ``pure" or ``non-randomized'' when $\bar \pi \in \overline \Pi^{p}$. We now define policies as sequences of strategy functions.

\begin{definition}%
	\label{def:MFMDP_mixed_Markov_policy}	
    A \defi{mixed Markov policy} for MFMDP is an element of $\bar\bPi := (\bar\Pi)^{\NN}$.
	Similarly, a \defi{pure policy} is an element of $\bar\bPi^{p} := (\bar\Pi^{p})^{\NN}$.
    We say that a policy $\bar\bpi=(\bar\pi_n)_{n\ge 0}$ is \defi{stationary} if the strategy functions $\bar \pi_n$ are equal for all $n$. 
\end{definition}
Keeping with the spirit of the previous section, these policies should be called ``Markov'' policies. We restrict ourselves to these policies and refrain from using history dependent policies because we are mostly interested in optimizing value functions, and we know that for MDPs like our lifted MFMDP, for each history dependent policy, there exists a Markov policy with the same value function (as defined below in Definition \ref{def:MFMDP_value_function_lsc}). See for example \cite[Proposition~9.1]{BertsekasShreve}.

\begin{definition}%
	\label{def:MFMDP_admissible_state_action_processes}	
	A pair of state and action processes $(\bmu, \bm{\bar a} ) = ( \mu_n, \bar a_n)_{n \geq 0}$ is said to be \defi{generated by} $(\bar \bpi, \mu) \in \bar\bPi \times \bar S$ if the following conditions are satisfied: $\mu_0 = \mu$, 
		\begin{equation}
			\label{eq:MFMDP_admissible_state_process}
			\mu_{n+1} = \bar F( \mu_n, \bar a_n, \varepsilon_{n+1}^0), \qquad \PP-a.s. \quad n\ge 0,
		\end{equation}
	and $ \bm{\bar a}$ is an $\bar A$-valued process adapted to $\FF^0$ satisfying:
		\begin{equation}
			\label{eq:MFMDP_admissible_action_process}
			\cL( \bar a_n | \mu_n ) = \bar \pi_n ( \mu_n ), \qquad \PP-a.s. \quad n\ge 0.
		\end{equation}
\end{definition}

\subsection{Assumptions and optimization problem for MFMDP}\,

We will sometimes rely on the following assumptions in the analysis of the model.

\begin{assumption}%
	\label{assumption:basic_assumption_MFC}
	\begin{itemize}
		\item \textbf{System function $F$:}   For every $(e, e^0) \in E \times E^0$, the function $F(\cdot, \cdot, \cdot, e, e^0)$ is  
		        continuous in its remaining variables.
		
		\item \textbf{One-stage cost function $f$:} $f: S \times A \times \cP(S \times A) \to \RR$ is continuous. 
	\end{itemize}
\end{assumption}

To show existence of an optimal policy, we will make use of the following extra assumption: %
\begin{assumption}
	\label{assumption:compactness_MFC} \textbf{Compactness}:	 The state space $S$ and the action space $A$ are compact metric spaces. 
\end{assumption}

In order to obtain a dynamic programming principle (DPP) with Borel measurable mixed Markov policies, measurability issues lead us to work under  Assumption~\ref{assumption:basic_assumption_MFC} for the MFC model. It can be shown that, under this assumption, $\bar F$ is Borel measurable and for every $e^0 \in E^0$, $\bar F(\cdot, \cdot, e^0)$ is continuous in its remaining variables, and $\bar f$ is bounded and lower semi-continuous. See Lemma~\ref{le:under_A_regu-barF-barf} in the appendix.

\begin{definition}%
	\label{def:MFMDP_value_function_lsc}
	For every $\bar\bpi \in \bar\bPi$, we define the \defi{value function} $J^{\bar \bpi}$ by:
	\begin{equation}
		\label{eq:MFMDP_J_BOREL}
		\bar J^{\bar \bpi}(\mu) := \EE \Big[ \sum_{n \geq 0} \gamma^n \bar f (\mu_n, \bar a_n) \Big], \qquad \mu \in \bar S,
	\end{equation}
	where $(\bmu, \bm{\bar a}) = (\mu_n, \bar a_n)_{n \geq 0}$ is any pair of state and action processes generated by $(\bar \bpi, \mu)$. If $\bar \bpi \in \bar\bPi$ is stationary with a strategy function $\bar \pi \in \bar\Pi$, then we let $\bar J^{\pi} = \bar J^{\bar \pi}$.%
\end{definition}

It can be shown that the value function  $\bar J^{\bar \bpi}$ given in~\eqref{eq:MFMDP_J_BOREL} is  well defined because the expectation in \eqref{eq:MFMDP_J_BOREL} does not depend upon the particular choice of the pair of state action processes $(\bmu,\bm{\bar a})$ generated by $(\bar \bpi, \mu)$.  See Lemma~\ref{lemma:MFMDP_well_definedness_value_function} in the appendix.

\begin{remark}
Using measurability arguments found for example in \cite[Chapter 7]{BertsekasShreve}, one can check that $\bar J^{\bar \bpi}$ is %
Borel  measurable when $\bar\bpi\in\bar\bPi$.

\end{remark}

With the value function for MFMDP at hand, we define the \defi{optimal value function of the MFMDP} as:  
	$$
	\bar J^{*}(\mu) = \inf_{ \bar \bpi \in \bar\bPi}{\bar  J}^{\bar \bpi}(\mu), \qquad \mu \in \bar S.
	$$

\subsection{Dynamic Programming principle for MFMDP}

We state and prove  the Dynamic Programming principle for the optimal value function with Borel measurable mixed Markov policies.

\begin{theorem}%
	\label{thm:MFMDP_DPP_BOREL}
	
	Assume that  \ref{assumption:basic_assumption_MFC} and \ref{assumption:compactness_MFC} hold. 
	Then, the function $\bar J^{*}$ is bounded and lower semi-continuous, and moreover it is the unique bounded and lower semi-continuous function satisfying the following dynamic programming equation with unknown $\bar J$: %
	\begin{equation}
		\label{eq:MFMDP_DPP_fixed_pt_BOREL}
		\bar J(\mu) = \inf_{ \bar a \in \bar U(\mu) } \left\{  \bar f(\mu, \bar a) + \gamma \EE \Big[ \bar J \big( \bar F( \mu, \bar a, \varepsilon^0) \big) \Big] \right\}, \qquad \mu \in \bar S.
	\end{equation}
	Furthermore, there exists a pure stationary $\bar \bpi^{*} = (\bar \pi^{*}, \bar \pi^{*}, \dots) \in \bar\bPi^{p}$ that is optimal, i.e., 
		$
		\bar J^{\bar \pi^*} = \bar J^{*}.
	$
\end{theorem}

The Dynamic Programming Principle \eqref{eq:MFMDP_DPP_fixed_pt_BOREL} is known to hold for universally measurable policies. See for example \cite[Proposition~9.8]{BertsekasShreve}.\footnote{Beware that the MFMDP setting is denoted with overlines in the notation of the present paper, but due to the common noise it corresponds to a stochastic model in \cite{BertsekasShreve}, which is denoted without overlines.} The gist of the above theorem is to show that it also holds for Borel measurable policies.

\begin{proof}

\textbf{Step 1: Bellman operator and fixed point with universal measurability. } 
We define the Bellman operator $\bar T$ by:
\begin{equation}
\label{fo:bar_T} 
    [\bar T \bar J](\mu) = \inf_{ \bar a \in \bar U(\mu) } \left\{  \bar f(\mu, \bar a) + \gamma \EE \Big[ \bar J \big( \bar F( \mu, \bar a, \varepsilon^0) \big) \Big] \right\}, \qquad \mu \in \bar S.
\end{equation}
Then, by \cite[Proposition~9.8]{BertsekasShreve}, $\bar T$ is a strict contraction on the space of bounded universally measurable functions on $\bar S$ and the fixed point coincides with the optimal value function for universally measurable policies, namely, $\bar J^{*,Univ}$ defined as:
$$
    \bar J^{*,Univ}(\mu) = \inf_{\bar\bpi \in \bar\Pi^{Univ}} \bar J^{\bar\bpi}, \qquad \mu \in \bar S
$$
where $\bar\Pi^{Univ}$ is the set of universally measurable mixed strategy functions, and for every $\bar\bpi \in \bar\Pi^{Univ}$, $\bar J^{\bar\bpi}$ is defined as in~\eqref{eq:MFMDP_J_BOREL} for the Boreal measurable case.

\textbf{Step 2: Fixed point with Borel measurability.} We will apply the Banach fixed point theorem for $\bar T$ on $\Lsc(\bar S)$, which denotes the set of real valued,  bounded and lower semi-continuous functions on $\bar S$. This set is a closed subset of the Banach space of real valued bounded functions on $\bar S$ endowed with the sup norm.

The key point is to show that $\bar T$ leaves $\Lsc(\bar S)$ invariant.  This follows by the measurable selection theorem for lower semi-continuous functions given in~\cite[Proposition~7.33]{BertsekasShreve}, since we can show that the content of the curly bracket in \eqref{fo:bar_T} is lower semi-continuous if $\bar J \in \Lsc(\bar S)$. Indeed, $\bar f$ is lower semi-continuous (see Lemma~\ref{le:under_A_regu-barF-barf} in the Appendix).  Moreover, since $\bar F$ is continuous for $\epsilon^0$ fixed (again by Lemma~\ref{le:under_A_regu-barF-barf}), the expectation is a continuous function of $(\mu,\bar a)$ whenever $\bar J$ is continuous by the dominated convergence theorem. Now since a function is lower semi-continuous if and only if it is the pointwise limit of an increasing sequence of continuous functions, we can use the monotone convergence theorem to show that the expectation is the limit of an increasing sequence of continuous functions, hence that it is lower semi-continuous. 

Last, $\bar T$ is a strict contraction for the sup norm. We thus conclude by the Banach fixed point theorem that $\bar T$ has a unique fixed point in $\Lsc(\bar S)$, which we denote by $\bar J^{*,bd,lsc}$. In other words, $\bar J^{*,bd,lsc}$ is the unique bounded lower semi-continuous function satisfying the dynamic programming equation~\eqref{eq:MFMDP_DPP_fixed_pt_BOREL}.

\textbf{Step 3: $\bar J^{*,bd,lsc} = \bar J^{*}$.} Being lower semi-continuous, $\bar J^{*,bd,lsc}$ is universally measurable, so it is also a fixed point in the space of universally measurable functions. Hence, by uniqueness, it coincides with $\bar J^{*,Univ}$. Furthermore $\bar\bPi \subseteq \bar\bPi^{Univ}$. So we deduce:
\begin{equation}
\label{fo:*}
    \bar J^{*,bd,lsc}(\mu)
    =\bar J^{*,Univ}(\mu)
    \le \bar J^{*}(\mu), \qquad \mu \in \bar S.
\end{equation} 
Assumption~\ref{assumption:compactness_MFC} implies that the lifted MFMDP satisfies the assumptions of \cite[Corollary 9.17.2]{BertsekasShreve} so there exists a stationary pure (at the level-1) Borel measurable policy $\bar\bpi^*=(\bar\pi^*,\bar\pi^*,\cdots)$ which is optimal in the sense that: 
$$
    \bar J^{\bar\bpi^*}(\mu)
    =\bar J^{*, Univ}(\mu),\qquad \mu\in\bar S.
$$
Consequently:
\begin{equation*}
    \bar J^{*}(\mu)
    \le \bar J^{\bar\bpi^*}(\mu)
    =\bar J^{*,Univ}(\mu)
    =\bar J^{*,bd,lsc} (\mu),\qquad \mu\in\bar S,
\end{equation*}
which together with \eqref{fo:*} gives the equality  $\bar J^{*} = \bar J^{*,Univ}$. 

Hence $\bar J^{*}$ satisfies the dynamic programming principle~\eqref{eq:MFMDP_DPP_fixed_pt_BOREL} and by Step 4, it is the unique bounded lower semi-continuous function satisfying it. 

\end{proof}

\begin{remark}
\label{remark:explain_universally_strategy_function}
	When the  mixed strategy function $\bar \pi_n \in \overline \Pi^{Univ}$ is only universally measurable, for each time $n$, as in the above proof of Theorem~\ref{thm:MFMDP_DPP_BOREL}, the understanding of condition~\eqref{eq:MFMDP_admissible_action_process} requires a modicum of care.
 Let $q_n = \cL( \mu_n ) \in \cP(\bar S)$ be the distribution of the random state $\mu_n$ with values in $\bar S$.
 We consider a Borel measurable kernel, $\bar \pi_{q_n} : (\bar S, \cB_{\bar S}) \to (\cP(\bar A), \cB_{\cP(\bar A)} )$, such that $\bar \pi_{q_n}(\mu) = \bar \pi_n(\mu)$ for $q_n$-almost every $\mu \in \bar S$ (see~\cite[Lemma~7.28 (c)]{BertsekasShreve} for existence). Then condition~\eqref{eq:MFMDP_admissible_action_process} says that $\bar \pi_{q_n}$ is a regular version of the conditional probability of $\bar a_n$ given $\mu_n$. 
 Furthermore, the integration of a function $\phi(\mu_n, \cdot)$ with respect to $\bar \pi_n(\mu_n)$ should be understood in the following sense:
 	$$
	     \int_{\bar A} \phi(\mu, \bar a) \bar \pi_n( \mu) (d \bar a)  = \int_{\bar A} \phi( \mu, \bar a) \bar \pi_{q_n}( \mu ) (d \bar a), 
	$$
	for $q_n-$almost every $\mu \in \bar S$, where $q_n = \cL(\mu_n)$.
 \end{remark}

\section{\textbf{Relations between the models}}\,
\label{sec:relations-models}

\subsection{Relations between MFC closed-loop policies and MFMDP policies}\,
\label{sec:relation-models-closedloop-barpi}

In this section, we discuss some of the connections between the original level-0 MFC model and the lifted MFMDP model. 

We start by highlighting that, intuitively, a closed-loop policy for the MFC can be viewed as sampled from a policy for the MFMDP by picking the common randomness. The following definition formalizes this idea. 

\begin{definition}%
	\label{def:correspondence_between_policies}
Let $\bpi \in \bPi^\tinycl$ and $\bar \bpi \in \bar\bPi$. We say that they \defi{correspond to each other} if  for each $\mu\in\bar S$ and $n\ge 0$, $\bar\pi_n(\mu)\in\cP(\bar A)$ is equal to the push forward of $\PP_{\vartheta^0}$ by the map:
$$
    \Theta^0\ni\theta^0\mapsto \mu(dx)\pi_n(x,\mu,\theta^0)(d\alpha)\in \bar A.
$$
\end{definition}
Note that if $\bpi$ and $\bar\bpi$ correspond to each other, then one is stationary if and only if the other one is. The main result of this section is the following.

\begin{theorem}%
	\label{thm:identical_value_function_MFC_MFMDP}
Assume~\ref{assumption:basic_assumption_MFC} holds. Then for every $\mu \in \cP(S)$, 
$
    \{J^{\bpi}(\mu) \,:\, \bpi \in \bPi^{\tinycl}\} = \{\bar J^{\bar \bpi}(\mu) \,:\, \bar \bpi \in \bar\bPi\}.
$
Similarly, for stationary policies, we have: for every $\mu \in \cP(S)$, 
$
    \{J^{\bpi}(\mu) \,:\, \bpi \in \bPi^{\tinycl} \hbox{ stationary}\} = \{\bar J^{\bar \bpi}(\mu) \,:\, \bar \bpi \in \bar\bPi \hbox{ stationary}\}.
$
\end{theorem}
This result means that for every $\bpi \in \bPi^{\tinycl}$, there exists  $\bar\bpi \in \bar\bPi$ such that $J^{\bpi} = \bar J^{\bar \bpi}$ and conversely, for every $\bar \bpi \in \bar\bPi$, there exists $\bpi \in \bPi^{\tinycl}$ such that this equality holds, with the same result in the stationary case.

In preparation for the proof of this result, we give two useful technical lemmas whose proofs are deferred to the appendix. These lemmas describe the properties of conditional distributions of the level-0 state and action processes. 

\begin{lemma}
\label{lemma:dynamics_of_conditional_distribution}
Assume~\ref{assumption:basic_assumption_MFC} holds. Let $\balpha \in \AA$, $\mu_0 \in \cP(S)$, and let $\bX$ be the associated state process. Then:
\begin{equation}
     \label{eq:formula_dynamics_in_lemma_MFC_MFMDP}
        \PP^0_{X_{n+1}} = \bar F( \PP^0_{X_n}, \PP^0_{(X_n, \alpha_n)}, \varepsilon_{n+1}^0 ), \qquad \PP-a.s. \qquad n \geq 0.
    \end{equation}
So $\cL(\PP^0_{X_{n+1}}) = P(\PP^0_{X_n}, \PP^0_{(X_n, \alpha_n)})$, $n \geq 0$, where the transition kernel $P$ was defined in~\eqref{def:MFMDP_transition_kernel_from_system_func}.
\end{lemma}

\begin{lemma}\,
\label{lemma:idenfiy_conditional_joint_dist_with_kernels}
 Assume~\ref{assumption:basic_assumption_MFC} holds. Let $\balpha \in \AA$, $\mu_0 \in \cP(S)$, and let $\bX$ be the associated state process. For every $n \ge 0$, let $\kappa_n: \bar S \to \cP(\bar A)$ be the Borel measurable disintegration kernel of $\cL( \PP^0_{X_n}, \PP^0_{(X_n, \alpha_n)} )$ along its first marginal. 
Then, if $(\bzeta, \bm{\bar \eta} )$ is an $(\bar S \times \bar A)$ - valued pair of stochastic processes  which are $\FF^0$ - adapted, and satisfy: 
    	    $\zeta_0 = \mu_0,$ $\PP-a.s.,$ 
    	    $\zeta_{n+1} = \bar F( \zeta_n, \bar \eta_n, \varepsilon_{n+1}^0)$, $\PP-a.s.$,  $n \ge 0$,
    	and if
    	    $\cL( \bar \eta_n  | \zeta_n ) = \kappa_n( \zeta_n ),$ $\PP-a.s.$ $n \ge 0$,
we have: 
    \begin{equation}
        \label{eq:identify_conditional_joint_distribution}
\cL (\zeta_n, \bar \eta_n) =\cL \big(\PP^0_{X_n}, \PP^0_{(X_n, \alpha_n)} \big), \qquad n \geq 0.
    \end{equation}
\end{lemma}

Intuitively, $(\kappa_n)_{n \ge 0}$ plays the role of the conditional law of the action process. We will come back to this interpretation in the proof of Lemma~\ref{lemma:construct_closed_loop}.  We are now ready to prove Theorem~\ref{thm:identical_value_function_MFC_MFMDP}.

\begin{proof}[Proof of Theorem~\ref{thm:identical_value_function_MFC_MFMDP}]

\textbf{Step 1:} Let $\bpi \in \bPi^\tinycl$. Let $\bar\bpi$ corresponding to $\bpi$ in the sense of Definition~\ref{def:correspondence_between_policies}. 
We now check the equality of the value functions $J^{\bpi}$ and $\bar J^{\bar\bpi}$. Note that if $\bar\bpi$ is stationary, then so is $\bpi$. Let $(\bX,\balpha)$ be a pair of state and action processes generated by $(\bpi,\mu_0)$. Then
	\begin{equation*}
		\begin{split}
			J^{\bpi}(\mu_0) 
			& = \EE \Bigl[ \sum_{n \geq 0} \gamma^n f \bigl( X_n, \alpha_n,\PP^0_{(X_n,\alpha_n)} \bigr) \Bigr]
			= \EE \Bigl[ \sum_{n \geq 0} \gamma^n \bar f \bigl( \PP^0_{X_n}, \PP^0_{(X_n, \alpha_n)} \bigr) \Bigr].
		\end{split}
	\end{equation*}
	We now check that $\mu_n=\PP^0_{X_n}$ and $\bar a_n=\PP^0_{(X_n,\alpha_n)}$ form a pair of state and action processes generated by $(\bar\bpi,\mu_0)$. This is indeed the case because \eqref{eq:MFMDP_admissible_state_process} is implied by Lemma~\ref{lemma:dynamics_of_conditional_distribution} and \eqref{eq:MFMDP_admissible_action_process} is implied by the definition of $\bar\pi_n$. Consequently,
	\begin{equation*}
    \bar J^{\bar \bpi}(\mu_0)  
    = \EE \Bigl[ \sum_{n \geq 0} \gamma^n \bar f \bigl( \PP^0_{X_n}, \PP^0_{(X_n, \alpha_n)}\bigr) \Bigr]
    = J^{\bpi}(\mu_0).
	\end{equation*}

\vskip 6pt
\textbf{Step 2:} Conversely, let $\bar \bpi = (\bar \pi_n)_{n \geq 0} $ in $\bar\bPi$. For every $n \geq 0$, 
	 $\bar \pi_n: \bar S \to \cP( \bar A)$ is a Borel measurable map such that for every $\mu\in \bar S$ we have  $\text{pr}_1(\bar\pi_n(\mu))=\mu$. According to the
 universal disintegration theorem~\cite[Corollary 1.26]{Kallenberg_RM}, there exists a Borel measurable probability kernel  $K: S \times \cP(S \times A) \times \cP(S) \to \cP(A)$ such that for every $\rho\in\cP(S\times A) $ and $\mu \in \cP(S)$ such that $\text{pr}_1(\rho)=\mu$, we have $\rho=\mu\measprod K(\cdot,\rho,\mu)$, where $\measprod$ denotes the product of a measure and a kernel.
So for every integer $n\ge 0$, $x\in S$, $\mu\in \bar S$ and $\theta^0\in\Theta^0$, we define:
	\begin{equation}
		\label{eq:construction_closed-loop_Markov_policy_MFC}
		\pi_n( x, \mu, \theta^0) := K\Bigl( x, \rho_{\bar A}\bigl(\bar \pi_n (\mu), h^0(\theta^0)\bigr),\mu\Bigr), 	
	\end{equation}
where $\rho_{\bar A}$ is the Blackwell-Dubins function of $\bar A$. Note that if $\bar\bpi$ is stationary, then so is $\bpi$. 
Because the functions $K$, $h^0$ and $\rho_{\bar A}$ are Borel measurable, so is the strategy function $\pi_n$ for every $n \geq 0$. Hence $\bpi = (\pi_n)_{n\ge 0} \in \bPi^{\tinycl}$. Recall that the function $ h^0$ was introduced in Section~\ref{sec:back-to-MFC}, and that $ h^0(\vartheta^0)$ is uniformly distributed on $[0,1]$ by construction. Notice that for every $\mu \in \bar S$ and for almost every $\theta^0 \in [0,1]$, the definition of the universal disintegration kernel $K$ implies that:
	\begin{equation}
	\label{eq:pi-K-barphi}
	    \Big(\rho_{\bar A}\bigl(\bar \pi_n (\mu), \theta^0\bigr) \Big) (d x, d\alpha)
	    =
		\mu(dx) K \Big( x,\rho_{\bar A}\bigl(\bar \pi_n (\mu), \theta^0\bigr) , \mu \Big)(d \alpha),
	\end{equation}
and as a result, we have:
	\begin{equation}
	\label{eq:mu-pin-rhoAbar}
		\Big(\rho_{\bar A}\bigl(\bar \pi_n (\mu),h^0(\theta^0\bigr) \Big) (d x, d\alpha)
		=\mu(dx)\pi_n( x, \mu, \theta^0)(d\alpha).
	\end{equation}
When $\theta^0$ is replaced by $\vartheta_n^0$, by Blackwell-Dubins lemma (see first point in Lemma~\ref{le:BlackwellDubins}), the left hand side of \eqref{eq:mu-pin-rhoAbar} is a random variable with values in $\bar A=\cP(S\times A)$ with distribution $\bar \pi_n (\mu)$. 

\vskip 2pt
Next, we show that $J^{\bpi} = \bar J^{\bar\bpi}$. Let $\mu_0 \in \bar S$. Let $(\bzeta, \bm{\bar{\eta}} )$ be state and action processes generated by $(\bar \bpi, \mu_0)$ (see Definition~\ref{def:MFMDP_admissible_state_action_processes}). Let $(\bX, \balpha)$ be a pair of  state and action processes generated by $\bpi$ and $\mu_0$.
Using the fact that $\PP^0_{(X_n,\alpha_n)}=\PP^0_{X_n}\measprod\pi_n(\cdot,\PP^0_{X_n},\vartheta^0_n)$, and the fact that $\vartheta^0_n$ is independent of $\PP^0_{X_n}$ by equation~\eqref{fo:P0X_n}, we have:
	\begin{equation*}
		\begin{split}
            J^{\bpi}(\mu_0)
            &= \sum_{n\ge 0}\gamma^n\EE\Bigl[\int_{S\times A}f(x,\alpha,\PP^0_{(X_n,\alpha_n)})\PP^0_{(X_n,\alpha_n)}(dx,d\alpha) \Bigr] %
            \\
            &= \sum_{n\ge 0}\gamma^n\EE\Bigl[\int_{S\times A}f\bigl(x,\alpha,\PP^0_{X_n}\measprod\pi_n(\cdot,\PP^0_{X_n},\vartheta^0_n)\bigr) \PP^0_{X_n}(dx)\pi_n(x,\PP^0_{X_n},\vartheta^0_n)(d\alpha)\Bigr]
            \\
            &= \sum_{n\ge 0}\gamma^n\EE\Bigl[\int_{S\times A}f(x,\alpha,\bar a)\bar a(dx,d\alpha)
\bar\pi_n\bigl(\PP^0_{X_n}\bigr)(d\bar a)\Bigr].
		\end{split}
	\end{equation*}
The last equality holds by the fact that both sides of~\eqref{eq:mu-pin-rhoAbar} with $\theta^0=\vartheta_n^0$ are random variables with values in $\bar A=\cP(S\times A)$ with distribution $\bar \pi_n (\mu)$. On the other hand:
	\begin{equation*}
		\begin{split}
			\bar J^{\bar \bpi}(\mu_0) 
			& = \EE \left[ \sum_{n \geq 0} \gamma^n \bar f ( \zeta_n, \bar \eta_n ) \right]
			\\
			& =  \sum_{n \geq 0} \gamma^n  \EE \left[ \int_{\bar A} \bar f ( \zeta_n, \bar a ) \cL( \bar \eta_n \, | \, \zeta_n)(d \bar a) \right]
			\\
			& = \sum_{n \geq 0} \gamma^n  \EE \left[ \int_{\bar A} \bar f ( \PP^0_{X_n}, \bar a ) \bar \pi_n ( \PP^0_{X_n} ) (d \bar a) \right],
		\end{split}
	\end{equation*}
where the last equality holds by~\eqref{eq:MFMDP_admissible_action_process} because $(\bzeta, \bm{\bar{\eta}} )$ are generated by $(\bar \bpi , \mu_0)$. This completes the proof.
\end{proof}

\vskip 2pt
At this stage, we need to emphasize the crucial role played by the common randomization provided by the sequence $(\vartheta^0_n)_{n\ge 0}$. Its presence is what allowed us to prove that the value of a policy for the lifted MDP can always be achieved by a closed loop policy of the original MFC.

Comparing to~\cite{motte2019mean}, while they prove equality of the optimal value functions without the central randomization, the latter allows us to prove the identity of the value functions, policy by policy, even before taking the optimum values.
\subsection{Relations between MFC closed-loop and open-loop policies}\,
\label{sec:back-to-MFC-2}

We first prove existence of optimal closed-loop Markov policies and then we prove equality of the open loop and closed loop value functions.

\begin{proposition}%
\label{proposition:existence_opt_policy_mean_field}

Assume~\ref{assumption:basic_assumption_MFC} and \ref{assumption:compactness_MFC} hold.
	There exists a stationary closed loop Markov policy for the original MFC that is optimal, i.e., $\bpi^{*} = (\pi^*,\pi^*,\dots) \in \bPi^{\tinycl}$ such that:
	$
		J^{\pi^{*}} = J^{\tinycl, *}.
	$
\end{proposition}

\begin{proof}
Let $\bar \bpi^{*}$ be an optimal pure stationary Markov policy for MFMDP whose existence is given in Theorem~\ref{thm:MFMDP_DPP_BOREL}, and let $\bpi^{*}\in\bPi^{CL}$ be a closed-loop Markov policy whose value function is the same and whose existence is given in Theorem~\ref{thm:identical_value_function_MFC_MFMDP} (for the case of stationary policies). We have:
$$
    J^{\bpi^{*}}(\mu) = \bar J^{\bar \bpi^{*}}(\mu) = \inf_{\bar \bpi \in \bar\bPi} \bar J^{\bar \bpi}(\mu), \qquad \mu \in \cP(S),
$$
Using Theorem~\ref{thm:identical_value_function_MFC_MFMDP} again, for every $\bpi\in \bPi^\tinycl$ for MFC, there exists $\bar \bpi \in \bar\bPi$ for MFMDP such that $    \bar J^{\bar \bpi}= J^{\bpi}$.
So, for every $\bpi \in \bPi^\tinycl$, 
$$
     J^{\bpi}(\mu) \geq \inf_{\bar \bpi \in \bar\bPi} \bar J^{\bar \bpi}(\mu) = J^{\bpi^{*}}(\mu), \qquad \mu \in \cP(S),
$$
which concludes the proof.
\end{proof}

\begin{remark}
Notice that because $\bar\bpi^*$ is pure, since the Blackwell-Dubins function $\rho_{\bar A}$ does not depend upon its second argument when the first is a point mass, we can conclude that $\pi_n$ does not depend upon the common randomization as given by $\vartheta^0_0$. In other words, the above result would still hold even if we did not have the common randomization.
\end{remark}

We now show the equality of the open-loop and closed-loop optimal value functions of the MFC.

\begin{theorem}%
		\label{thm:equality_value_function_open_closed_Markov}
Assume~\ref{assumption:basic_assumption_MFC} holds.
		Then %
		$
		J^{\tinyol, *} = J^{\tinycl, *}.
		$
\end{theorem}

This result is a direct consequence of the following Lemma~\ref{lemma:construct_open_loop} and Lemma~\ref{lemma:construct_closed_loop}, together with Theorem~\ref{thm:identical_value_function_MFC_MFMDP}. Indeed, by Lemma~\ref{lemma:construct_open_loop} and Lemma~\ref{lemma:construct_closed_loop},
$$
    J^{\tinycl,*} \ge J^{\tinyol,*} \ge \bar J^{*}. 
$$ 
Moreover, by Theorem~\ref{thm:identical_value_function_MFC_MFMDP}, 
$$
    \bar J^{*} = J^{\tinycl,*}.
$$
Hence the above inequalities are equalities, which proves the first part of Theorem~\ref{thm:equality_value_function_open_closed_Markov}. The existence of an optimal closed-loop policy stems from Theorem~\ref{thm:identical_value_function_MFC_MFMDP}, which entails the existence of an optimal open-loop policy by Lemma~\ref{lemma:construct_open_loop}.
Again, while this type of equality between the optimal value functions could be expected to hold under different assumptions and without the central randomization, we prove it here by leveraging the equalities proven in Lemma~\ref{lemma:construct_open_loop} and Lemma~\ref{lemma:construct_closed_loop} policy by policy, before computing optima over sets of policies.

\vskip 6pt
First, it is expected that every closed-loop policy for the MFC can be viewed as an open-loop policy for the MFC, which leads to  the following result in terms of value functions. 

\begin{lemma}%
\label{lemma:construct_open_loop}
Assume~\ref{assumption:basic_assumption_MFC} holds.
For every $\tilde{\bpi} \in \bPi^{\tinycl}$, there exists $\bpi \in \bPi^{\tinyol}$ such that: 
    $
        J^{ \bpi} = J^{\tilde{\bpi}}.
    $
\end{lemma}
The proof is deferred to the Appendix~\ref{app:proofs-back-MFC-2}.
Next, every open-loop policy for the MFC corresponds to a policy for the MFMDP, as we show in the following result. 

\begin{lemma}%
	\label{lemma:construct_closed_loop}
Assume~\ref{assumption:basic_assumption_MFC} holds. For every $\bpi \in \bPi^{\tinyol}$, there exists  $\bar \bpi \in \bar\bPi$ such that: 
	$
	    \bar J^{\bar \bpi} = J^{\bpi}.
	$
\end{lemma}

\begin{proof} 
Let $\bpi \in \bPi^{\tinyol}$. 
Let us fix an initial distribution $\mu_0 \in \cP(S)$, let 
$ \balpha$ be an action process generated by $\bpi$, and $\bX$ be the state process associated with $(\balpha, \mu_0)$ (recall Definition~\ref{def:state_process_from_control_process}).
For each $n\ge 0$, we consider the probability kernel $\kappa_n: \bar S \to \cP(\bar A)$ defined in the statement of Lemma~\ref{lemma:idenfiy_conditional_joint_dist_with_kernels}.  
We construct by induction an $(\bar S \times \bar A)$-valued pair of processes $(\bzeta, \bm{\bar \eta})$ in the following way.
For $n=0$ we set $\zeta_0 = \mu_0$ and $\bar\eta_0=\rho_{\bar A}\bigl(\kappa_0(\zeta_0),h^0(\vartheta^0_0)\bigr)$ where $\rho_{\bar A}$ is the Blackwell-Dubins function of the space $\bar A$ introduced in Lemma~\ref{le:BlackwellDubins}. Then for any $n\ge 0$ we define: Obviously, each time we involve $\vartheta^0_0$, we rely on the central randomization. Still, the conclusion of this lemma should not be considered as obvious because it proves that we can pack all the dependence on the past carried by the open loop controls at level $0$ into $\PP^0_{X_n},  \PP^0_{\tinymath[ \big( X_n, \alpha_n \big)]}$ and  a probability measure on the space of actions at level $1$.
$$
\zeta_{n+1} = \bar F( \zeta_n, \bar \eta_n, \varepsilon_{n+1}^0 ),
\quad\text{and}\quad
\bar\eta_{n+1}=\rho_{\bar A}\bigl(\kappa_{n+1}(\zeta_{n+1}),h^0(\vartheta^0_{n+1})\bigr).
$$
The process $\bm{\bar \eta}$ is adapted to $\FF^0$ and by Lemma~\ref{le:BlackwellDubins} it satisfies
$$
    \cL( \bar \eta_n \, | \, \zeta_n ) = \kappa_n( \zeta_n), \qquad \PP-a.s.\,, \quad n \geq 0, 
$$
because $\vartheta^0_{n+1}$ is independent of $\zeta_n$. Thus, by Lemma~\ref{lemma:idenfiy_conditional_joint_dist_with_kernels}:
\begin{equation}
    \label{eq:prop-zeta-bareta-law}
    \cL(\zeta_n, \bar \eta_n) =\cL \Big(\PP^0_{X_n},  \PP^0_{\tinymath[ \big( X_n,  \alpha_n \big)] }  \Big), \qquad n \geq 0.
\end{equation}
Let $\bar \bpi = (\kappa_n)_{n \geq 0}$. Since $\kappa_n(\bar U(\mu) ) = 1$ for every $n \geq 0$, we see that $\bar \bpi \in \bar\bPi$. We conclude by noting that:
$$
    \bar J^{\bar \bpi}(\mu_0) = \sum_{n=0}^\infty \gamma^n \EE\left[ \bar f( \zeta_n, \bar \eta_n) \right] = \sum_{n=0}^\infty \gamma^n \EE\left[ \bar f \Big( \PP^0_{X_n},  \PP^0_{\tinymath[ \big( X_n, \alpha_n \big)] } \Big) \right]  = J^{\bpi}(\mu_0),
$$
where the second equality holds by~\eqref{eq:prop-zeta-bareta-law}. 
\end{proof}

\section{\textbf{Mean-Field Q-Learning}}
\label{se:Q_learning}

\subsection{State-action value function}\,

We now turn our attention to the question of \emph{learning} the solution of the MFC problem in a model-free setting, i.e., assuming the model is unknown while still having access to sample realizations of state trajectories and associated rewards. Before considering algorithms, we first study the so-called state-action value function.

\vskip 2pt
In order to take advantage of the strongest results proven so far, we now assume both \ref{assumption:basic_assumption_MFC} and \ref{assumption:compactness_MFC} hold.
Under these assumptions, recall that the DPP given by Theorem~\ref{thm:MFMDP_DPP_BOREL} holds. In this section, we restrict ourselves to non-randomized stationary policies. When $\bar \bpi = (\bar \pi, \bar \pi, \ldots)$, we use freely the notation $\bar J^{\bar \pi} := \bar J^{\bar \bpi}$. Note however that we will use common randomization in the numerical section for the purpose of exploration.

In this section, without any loss of generality, we restrict the search for optimal policies to the set:
	$$
	\bar U_B(\bar A | \bar S) = \{  \bar \pi \in \bar\Pi^p \   | \  \forall \mu \in \bar S, \ \bar \pi (\mu) \in \bar U(\mu) \  \}.
	$$ 
For each $\bar \pi \in \bar U_B(\bar A | \bar S)$, the mapping
$\bar S \ni \mu \mapsto \delta_{\bar \pi(\mu)} \in \cP(\bar A)$ which assigns to each $\mu \in \bar S$ the Dirac point mass at the point $\bar \pi(\mu) \in \bar A$ is a Borel measurable function by definition of the Borel $\sigma$-field of $\cP(\bar A)$.

Now, for each $\bar\pi \in \bar\Pi$, we introduce the state-action value function $\bar Q^{\bar \pi} : \bar\Gamma \to \RR$ defined by: 
\begin{equation}
\label{eq:Q_function_stationary_control}
    \bar Q^{\bar \pi} (\mu, \bar a) := \bar f(\mu, \bar a) + \sum_{n \geq 1} \gamma^{n} \EE [\bar f(\mu_n, \bar\pi(\mu_n))], \qquad (\mu, \bar a) \in \bar\Gamma,
\end{equation}
where the process $(\mu_n)_{n\geq 0}$ starting at $\mu_0 = \mu$ satisfies $\mu_1 = \bar F(\mu, \bar a, \varepsilon_1^0)$, and for every $n \geq 1$:
$
	\mu_{n+1} = \bar F \bigl( \mu_n, \bar \pi(\mu_{n}), \varepsilon_{n+1}^0 \bigr).
$
Next we define the optimal state-action value function by:
\begin{equation}
\label{fo:optimal_Q}
    \bar Q^* (\mu, \bar a) := \inf_{\bar\pi\in\bar U_B(\bar A|\bar S)}\bar Q^{\bar \pi} (\mu, \bar a), 
\qquad
(\mu,\bar a)\in \bar\Gamma.
\end{equation}

The main goal of this section is to prove the following dynamic programming principle for $\bar Q^*$.

\begin{theorem} 
 	\label{prop:opt_Bellman_eq_Q}
Assume \ref{assumption:basic_assumption_MFC} and \ref{assumption:compactness_MFC} hold. The optimal state-action value function $\bar Q^*$ satisfies the so-called \defi{Bellman equation for state-action value function}:
\begin{equation}
\label{eq:Bellman_equation_Q}
	\bar Q^*(\mu, \bar a) = \bar f(\mu, \bar a) + \gamma \EE \Bigl[ \inf_{ \bar a' \in \bar U \big(\bar F(\mu, \bar a, \varepsilon^0) \big) }  \bar Q^*\big( \bar F(\mu, \bar a, \varepsilon^0) , \bar a'\big) \Bigr], \qquad (\mu, \bar a) \in \bar\Gamma.
\end{equation}
\end{theorem}

We will prove this result by showing that $\bar Q^*$ is the unique fixed point of state-action Bellman operator $T$ defined on the set $\Lsc(\bar\Gamma)$ of bounded lower semi-continuous functions on $\bar\Gamma$, by:
\begin{equation}
\label{fo:state-action_Belman}
    [T \bar Q](\mu,\bar a) := \bar f(\mu,\bar a) 
+\gamma\EE \Bigg[ \inf_{\bar a’\in \bar U \big( \bar F(\mu, \bar a, \varepsilon^0) \big) } \bar Q\bigl(\bar F(\mu,\bar a,\varepsilon^0),\bar a’ \bigr) \Bigg],\qquad
(\mu,\bar a)\in \bar\Gamma.
\end{equation}
We first justify in Lemma~\ref{le:strict_contraction} the fact that the operator $T$ is well-defined.

Since $\bar\Gamma$ is a closed subset of the Polish space $\bar S \times \bar A$, it is a Borel space, and the space $\cBD_u(\bar\Gamma)$ of bounded real-valued universally measurable functions on $\bar\Gamma$ endowed with the sup norm $\| f \|_{\infty} = \sup_{(\mu, \bar a) \in \bar\Gamma} | f(\mu, \bar a) |$ is a Banach space. While the set $\Lsc(\bar\Gamma)$ is not a vector space, it is a closed subset of $\cBD_u(\bar\Gamma)$, 
hence a complete metric space for the metric $d_\infty(f, f') = \| f - f' \|_{\infty}$.

\begin{lemma}
\label{le:strict_contraction}
Assume \ref{assumption:basic_assumption_MFC} and \ref{assumption:compactness_MFC} hold. The set $\cL sc(\bar\Gamma)$ is invariant under the state-action Bellman operator $T$, which is a strict contraction on this metric space.
\end{lemma}

\begin{proof}
We first claim that the set  $\cL sc(\bar\Gamma)$ is invariant under $T$. We need to show that $T\bar Q$ is lower semi-continuous whenever $\bar Q$ is. To wit, by the projection property for infima of lower semi-continuous functions (see for example \cite[Proposition~7.33]{BertsekasShreve}), the function $\bar S  \ni \mu' \mapsto \inf_{\bar a' \in \bar U(\mu') }\bar Q(\mu',\bar a')$ is lower semi-continuous.
Since $\bar\Gamma \ni (\mu,\bar a)\mapsto \mu'=\bar F(\mu,\bar a,e^0)\in \bar S$ is continuous for $e^0\in E^0$ fixed, the infimum in formula \eqref{fo:state-action_Belman} is then a lower semi-continuous function of $(\mu,\bar a)$ for fixed $e^0\in E^0$. Finally, Fatou's theorem implies that the expectation in \eqref{fo:state-action_Belman} is a lower semi-continuous function of $(\mu,\bar a) \in \bar\Gamma$. Furthermore $\bar f$ is continuous. Consequently, $T\bar Q$ is also lower semi-continuous.
Now if $\bar Q_1$ and $\bar Q_2$ are elements of $\cL sc(\bar\Gamma)$ we have:
\begin{align*}
	\| T \bar Q_1 - T \bar Q_2 \|_{\infty} 
	& \le   \gamma  \EE \Bigg[ \sup_{(\mu, \bar a) \in \bar \Gamma} \ \Bigl|\inf_{\bar a'\in \bar U(\bar F(\mu,\bar a,\varepsilon^0) )} \bar Q_1\bigl(\bar F(\mu,\bar a,\varepsilon^0),\bar a'\bigr) - \inf_{\bar a'\in\bar U(\bar F(\mu,\bar a,\varepsilon^0) )} \bar Q_2\bigl(\bar F(\mu,\bar a,\varepsilon^0),\bar a'\bigr) \Bigr| \ \Bigg]
	\\
	& \le   \gamma  \EE \Bigg[ \sup_{(\mu, \bar a) \in \bar \Gamma}
	\sup_{\bar a'\in \bar U(\bar F(\mu,\bar a,\varepsilon^0) )}\ \Bigl| \bar Q_1\bigl(\bar F(\mu,\bar a,\varepsilon^0),\bar a'\bigr) - \bar Q_2\bigl(\bar F(\mu,\bar a,\varepsilon^0),\bar a'\bigr) \Bigr|\ \Bigg]\\
	& \le \gamma \| Q_1 - Q_2 \|_{\infty}.
\end{align*}
Since $\gamma <1$, this proves that $T$ is a strict contraction on $\cL sc(\bar\Gamma)$. We conclude the proof of the result using the Banach fixed point theorem.
\end{proof}

Using the Markov property, we can rewrite the state-action value function $\bar Q^{\bar \pi}$
in terms of the state value function $\bar J^{\bar \pi}$: 

\begin{equation}
    \label{eq:connection_J_Q}
    \bar J^{\bar \pi}(\mu) = \bar Q^{\bar \pi}\bigl( \mu, \bar \pi(\mu)\bigr), \qquad  \bar \pi \in \bar\Pi^p, \mu \in \bar S.
\end{equation}

Now, for the optimal value functions, we have the following.
\begin{lemma}
	\label{le:relationship_Q*_J*}
Assume \ref{assumption:basic_assumption_MFC} and \ref{assumption:compactness_MFC} hold. For all $(\mu, \bar a) \in \bar\Gamma$: $\bar Q^*(\mu, \bar a) = \bar f(\mu, \bar a) + \gamma \EE[  \bar J^{*}(\bar F(\mu, \bar a, \varepsilon^0)) ]$. 
\end{lemma}

\begin{proof}
We show the inequalities in both directions. First, for  $(\mu, \bar a) \in \bar\Gamma$ given, let $q = \cL( \mu_1 ) \in \cP(\bar S)$ be the distribution of the random measure $\mu_1 = \bar F(\mu, \bar a, \varepsilon^0)$.  By~\cite[Corollary 9.5.2]{BertsekasShreve} and the fact that $\{ (\varphi, \varphi, \ldots) \text{ with } \varphi \in \bar U_B(\bar A | \bar S) \} \subseteq \bar\bPi$, we have:
$$
    \int_{\bar S} \bar J^{*}(\mu) q(d\mu) 
    = \inf_{ \bar \bpi \in \bar\bPi} \int_{\bar S} \bar J^{\bar \bpi}(\mu) q(d \mu)
    \le \inf_{ \substack{\bar \bpi = (\varphi, \varphi, \ldots); \\ \varphi \in \bar U_B(\bar A | \bar S) } } \int_{\bar S} \bar J^{\bar \bpi}(\mu) q(d \mu).
$$
Hence 
$
    f(\mu, \bar a) + \gamma \EE\left[ \bar J^{*}(\mu_1) \right] 
     \leq 
     \bar Q^*(\mu, \bar a).
$

Conversely, under the standing assumptions, by Theorem~\ref{thm:MFMDP_DPP_BOREL} there exists $\bar \pi^* \in \bar U_B( \bar A | \bar S)$ that is optimal. So we have:
 \begin{equation*}
   \bar Q^*(\mu, \bar a)  =  \inf_{\bar \pi \in \bar U_B(\bar A | \bar S) } \left\{ f(\mu, \bar a) + \gamma \EE\left[ \bar
    J^{\bar \pi}(\mu_1) \right] \right\}
    \leq f(\mu, \bar a) + \gamma \EE\left[ \bar J^{\bar \pi^* }(\mu_1) \right]
     = f(\mu, \bar a) + \gamma \EE\left[ \bar J^{*}(\mu_1) \right],
\end{equation*}
for every $(\mu, \bar a) \in \bar\Gamma$, which concludes the proof.
\end{proof}

\begin{lemma}
\label{le:Qstar_is_lsc}
 Assume \ref{assumption:basic_assumption_MFC} and \ref{assumption:compactness_MFC} hold. $\bar Q^*$ is lower semi-continuous and, as a result, there exists $\tilde \pi\in\bar U_B(\bar A|\bar S)$ such that for every $\mu \in \bar S$,
$
\tilde\pi(\mu) \in {\arg\inf}_{ \bar a \in \bar U(\mu) } \bar Q^*(\mu,\bar a).
$
\end{lemma}

\begin{proof}
For each fixed $e^0 \in E^0$, $(\mu,\bar a)\mapsto \bar J^*\bigl(\bar F(\mu,\bar a,e^0)\bigr)$ is lower semi-continuous by lower semi-continuity of $\bar J^{*}$ (see Theorem~\ref{thm:MFMDP_DPP_BOREL}) and continuity of $\bar F(\cdot,\cdot,e^0)$.
As in the proof of Lemma~\ref{le:strict_contraction}, Fatou's theorem implies that the function $(\mu,\bar a)\mapsto \EE\bigl[\bar J^*\bigl(\bar F(\mu,\bar a,e^0)\bigr)\bigr]$ is also lower semi-continuous. 
Furthermore, $\bar f$ is continuous. So, by the expression in Lemma~\ref{le:relationship_Q*_J*}, $\bar Q^*$ is a lower semi-continuous function on $\bar\Gamma$.

Since $\bar\Gamma$ is a closed subset of $\bar S \times \bar A$ and $\bar A$ is compact, by applying a selection theorem for lower semi-continuous function~\cite[Proposition~7.33]{BertsekasShreve} on $\bar Q^*: \bar\Gamma \to \RR$, we obtain that there exists a Borel measurable function $\tilde \pi \in \bar U_B(\bar A | \bar S)$ whose graph is contained in $\bar\Gamma$ and 
    $\bar Q^* \big(\mu, \tilde \pi(\mu) \big) = \inf_{\bar a \in \bar U(\mu) } \bar Q^*(\mu, \bar a),$
    for all $\mu \in \bar S.$
\end{proof}

\begin{lemma}
\label{lemma:infimum_action_opt_barQ_and_opt_barJ}
Assume \ref{assumption:basic_assumption_MFC} and \ref{assumption:compactness_MFC} hold. For all $\mu \in \bar S$, 
$\inf_{\bar a\in \bar U(\mu)} \bar Q^*(\mu, \bar a) = \bar J^{*}(\mu)$.
\end{lemma}

\begin{proof}
    We first show the inequality $ \inf_{\bar a\in \bar U(\mu)} \bar Q^*(\mu, \bar a) \geq \bar J^*(\mu)$. 
    Let us denote by ${\tilde\pi} \in \bar U_B(\bar A | \bar S)$ the strategy function in Lemma~\ref{le:Qstar_is_lsc}.
    By definition, 
    $$
    	\inf_{\bar a \in \bar U(\mu) } \bar Q^*(\mu, \bar a) = \bar Q^*(\mu, {\tilde\pi}(\mu) ) = \inf_{\bar \pi \in \bar U_B(\bar A | \bar S) } \bar Q^{\bar \pi} (\mu, {\tilde\pi}(\mu) ), \qquad \mu \in \bar S.
	$$ 
	Then for each $\bar \pi \in \bar U_B(\bar A | \bar S)$, we denote 
$
         \bar \bpi^{{\tilde\pi}} = ( {\tilde\pi}, \bar \pi, \bar \pi, \bar \pi, \ldots ) \in \bar\bPi.
$
    So, %
    $$
        \bar Q^{\bar \pi}(\mu, {\tilde\pi}(\mu) ) 
        = f(\mu, {\tilde\pi}(\mu) ) + \EE\left[ \sum_{n=1}^\infty \gamma^n \bar f( \mu_n, \bar \pi(\mu_n)  ) \right] 
        = \bar J^{ \bar \bpi^{{\tilde\pi}} } (\mu) \geq \bar J^*(\mu), \qquad \  \mu \in \bar S,
    $$
    which provides the first inequality. 
	To prove the converse inequality, let $\bar \pi^* \in \bar U_B(\bar A | \bar S)$ be an optimal non-randomized stationary Markov policy whose existence is given by Theorem~\ref{thm:MFMDP_DPP_BOREL}, and  notice that for every $\mu \in \bar S$, 
    $$
        \bar J^{*}(\mu) = \bar J^{\bar \pi^*}(\mu) = \bar Q^{\bar \pi^*}(\mu, \bar \pi^*(\mu) ) \geq \bar Q^*(\mu, \bar \pi^*(\mu) ) \geq \inf_{\bar a \in \bar U(\mu)} \bar Q^*(\mu, \bar a),
    $$
    by equation~\ref{eq:connection_J_Q}, Lemma~\ref{le:relationship_Q*_J*}. %
This concludes the proof.
\end{proof}

We can now complete the proof of Theorem~\ref{prop:opt_Bellman_eq_Q}. 
\begin{proof}[Proof of Theorem~\ref{prop:opt_Bellman_eq_Q}]
The Bellman equation~\eqref{eq:Bellman_equation_Q} is a direct consequence of Lemma~\ref{le:relationship_Q*_J*} and 
 Lemma~\ref{lemma:infimum_action_opt_barQ_and_opt_barJ}. Since $T$ is a strict contraction mapping on $\Lsc(\bar\Gamma)$ by Lemma~\ref{le:strict_contraction} and since $\Lsc(\bar\Gamma)$ is closed in the Banach space $\cBD_u(\bar\Gamma)$, by the Banach fixed point theorem we conclude that $\bar Q^*$ is the unique fixed point of $T$ on $\Lsc(\bar\Gamma)$.
\end{proof}

	Next, we build upon the previous results to propose reinforcement learning algorithms for the original MFC problem. From now on, we assume that the state and action spaces are finite, unless otherwise specified. 

\subsection{Controls for finite state and action spaces}\,

In this rest of this section, we assume that $S$ and $A$ are finite, we denote their numbers of elements by $|S|$ and $|A|$ respectively, and we denote by $x^{(1)},\dots,x^{(|S|)}$ and  $\alpha^{(1)},\dots,\alpha^{(|A|)}$  their elements. We first revisit the description of the action space and then propose two reinforcement learning methods in this setting. We shall explain later how to adapt reinforcement learning techniques to the case of continuous spaces.  

	Before introducing the mean-field Q-learning algorithm, we first provide a representation of the set $\bar \Gamma \subseteq \bar S \times \bar A = \cP(S) \times \cP(S \times A)$ on which the $\bar Q^*$ function is defined. 
	
	Since we assume that $S$ finite, its lifted space $\cP(S)$ can be identified with a simplex $\frak{S}$ in $\RR^{| S |}$. In other words, we treat a distribution $\mu \in \cP(S)$ as an $|S|$-dimensional vector $(\mu^{(i)})_{i=1,\dots,|S|}$ whose non-negative coordinates sum up to one. Similarly, since $A$ is finite, we identify $\cP(A)$ to a simplex $\frak{A}$ in $\RR^{|A|}$. 
	However, representing admissible actions $\bar a \in \bar U(\mu) \subseteq \cP(S \times A) $  of the lifted MDP requires a modicum of care due to the constraint. A first approach is to identify $\cP(S \times A)$ with a simplex in $\RR^{|S| \times |A|}$ and to view a lifted action $\bar a$ as a $|S|\times |A|$ matrix $ \big( \bar a( x^{(i)}, \alpha^{(j)}) \big)_{1 \leq i \leq |S|, 1 \leq j \leq |A|}$ of non-negative numbers summing up to $1$. Then a pair $(\mu, \bar a) \in \bar S \times \bar A$ is in $\bar \Gamma$ if and only if the following linear constraint is satisfied: $\sum_{j=1}^{|A|} \bar a( \mu^{(i)}, \alpha^{(j)} ) = \mu^{(i)}$ for all $i = 1, \ldots, |S|$.
The above transformation is straightforward but not sufficient for our purposes because it provides only a representation of the actions and controls of the central planner, and it does not address the strategy functions of non-randomized stationary mixed Markovian closed-loop policies for an individual agent in our original optimization problem.  

	For any pair $(\mu, \bar a) \in \bar \Gamma$, we can define the mapping $k_{\mu}: S \to \cP(A)$: for $i=1,\dots,|S|$, 
	\begin{equation*}
	k_\mu(x^{(i)}) = \bar a( \mu^{(i)}, \alpha^{(j)} )/\mu(x^{(i)}), \qquad \hbox{ if } \mu(x^{(i)})>0,
	\end{equation*}
	and any value otherwise. 
	Note that here, there is no common randomization. As proved above (see Theorem~\ref{thm:MFMDP_DPP_BOREL}), there exists a non-randomized stationary policy for the lifted MDP. So the central planner can look for strategy functions within the set: 
	\begin{equation}
	\label{eq:def-calA-S-PA}
		\cA := \{ \tilde a: S \to \cP(A) \  | \   \tilde a \text{ Borel measurable} \}.
	\end{equation}
	Consider the function $\tilde Q^*: \cP(S) \times \cA \to \RR$ defined by: 
	\begin{equation}
		\tilde Q^*(\mu, \tilde a) := \bar Q^*(\mu, \mu \measprod \tilde a).
	\end{equation}
	Then the Bellman equation~\eqref{eq:Bellman_equation_Q} becomes:
	\begin{equation}
		\label{eq:optimal_Q_kernel_version}
		\tilde Q^*(\mu, \tilde a) = \int_{S \times A} f(x, \alpha, \mu \measprod \tilde a )  \tilde a(x, d\alpha) \mu( dx)+ \gamma \EE \left[ \inf_{\tilde a'\in \cA} \tilde Q^*( \mu_1, \tilde a')\right], \qquad (\mu, \tilde a) \in \cP(S) \times \cA,
	\end{equation}
	where $\mu_1 = \bar F( \mu, \mu \measprod \tilde a, \varepsilon^0)$, keeping in mind that the integral over $S\times A$ is in fact a finite sum.
Even though $S$ and $A$ are finite,	equation \eqref{eq:optimal_Q_kernel_version} still needs to be understood as a fixed point in the space of bounded lower semi-continuous functions on a closed subset of a finite dimensional Euclidean space, as the measurability issues addressed in deriving equation~\eqref{eq:Bellman_equation_Q} still remain. 
We also introduce the function $\tilde f: \cP(S) \times \cA \to \RR$ such that:
	$$
		\tilde f(\mu, \tilde a) := \bar f(\mu, \mu \measprod \tilde a) = \int_{S \times A} f(x, \alpha, \mu \measprod \tilde a)  \tilde a(x, d\alpha) \mu( dx), \qquad (\mu, \tilde a) \in \cP(S) \times \cA.
	$$
In the rest of this section, we propose two model-free algorithms relying on the optimal state-action value function $\bar Q^* : \bar \Gamma \to \RR$ or equivalently $\tilde Q^*: \cP(S) \times \cA \to \RR$.

\subsection{Simplex discretization and tabular MFQ-learning}\,
\label{sec:tabularQ-discrete}

We consider two settings, depending on whether the controls at level-0 are mixed or pure. In both cases, we prove convergence of a tabular Q-learning algorithm, after suitable discretization of the simplexes. When using pure controls, we can prove not only convergence of the value function but also of the optimizer. 

\subsubsection{Q-learning with controls that are mixed at level-0}
\label{sec:tabularQ-mixed}
Since the simplexes $\frak{S}$ and $\frak{A}$ are not finite, it is not possible to directly apply a tabular version of Q-learning algorithm to approximate $\tilde Q^*$. A possible workaround is to first replace these simplexes by finite subsets $\check{\frak{S}} \subset \frak{S}$ and $\check{\frak{A}} \subset \frak{A}$. Let $\check \cA = \{ \check a : S \to \check{\frak{A}} \}$. In particular, $|\check \cA| = | \check{\frak{A}} |^{|S|}$ because we identify functions in $\check{\cA}$ with $|S|$-dimensional vectors  whose entries take values in the finite set $\check{\frak{A}}$. To ensure that the mean-field term takes values in the finite set $\check{\frak{S}}$, we use a projection: 
at time $n$, given $\mu_n \in \check{\frak{S}}$, we compute $\mu_{n+1} = \bar F( \mu_n, \mu_n \measprod \check a, \varepsilon_{n+1}^0)$, and then we project $\mu_{n+1}$ back on $\check{\frak{S}}$ using a projection operator $\proj_{\check{\frak{S}}}: \cP(S) \to \check{\frak{S}}$. Precise definitions of the discretization and the projection are provided below, after introducing a discrete version of the original MFC problem.  

More precisely, we consider the \defi{projected MFC problem}: 
$$
	\inf_{ \check \pi \in \check \Pi} \check J_{\check \pi}(\mu_0),\qquad \mu_0 \in \check{\frak{S}},
$$
where $\check \Pi = \{ \check \pi: S \times \check{\frak{S}} \to \check{\frak{A}} \}$, and for every strategy function $\check \pi : S \times \check{\frak{S}} \to \check{\frak{A}}$,
$\check J^{\check \pi}: \check{\frak{S}} \to \RR$ is defined by:
\begin{equation}
\label{eq:generic-MFC-fctmeasure-reward-proj}
	\check J^{\check \pi}(\mu_0) =  \EE \left[ \sum_{n \geq 0} \gamma^n  \tilde f \Big( \mu^{\mu_0, \check \pi}_n,  \check \pi(\cdot, \mu^{\mu_0, \check \pi}_n)  \Big) \right]
\end{equation}
where 
\begin{equation}
\label{eq:generic-MFC-fctmeasure-dyn-proj}
	\mu^{\mu_0, \check \pi}_{n+1} = \proj_{\check{\frak{S}}} \circ  \bar F \Big( \mu_n^{\mu_0, \check \pi}, \mu_n^{\mu_0, \check \pi} \measprod \check \pi( \cdot, \mu_n^{\mu_0, \check \pi}), \varepsilon_{n+1}^0 \Big) =: \check \Phi^{\check \pi, \varepsilon_{n+1}^0}(\mu_n^{\mu_0, \check \pi}).
\end{equation}
We will denote by $\check J^*$ and $\check Q^*$ respectively the optimal state and state-action value functions of this projected MFC problem. Here $ \check Q: \check{\frak{S}} \times \check \cA \to \RR$ can be represented by a matrix (also called a table) in $\RR^{ |\check{\frak{S}}| \times | \check \cA| }$ and is viewed as an approximation of $\tilde{Q}^* : \cP(S) \times \cA \to \RR$ of the original MFC problem. 

This problem can be viewed as an MDP with finite state and action spaces. In this case, a straightforward adaptation of the tabular Q-learning algorithm leads to Algorithm~\ref{algo:Qtable-projection}. Note that, even in the absence of common noise, this algorithm is possibly stochastic since at each episode, the order in which the state-action pairs are picked is potentially random. In practice, the order could be fixed in advance or stem from a sampled trajectory.  

\begin{algorithm}[h!]
\DontPrintSemicolon
\KwData{A number of episodes $N_{\mathrm{epi}}$; a sequence of learning rates $(\eta_n)_{n=0,\dots,N_{\mathrm{epi}}-1}$; a sequence of state-action pairs $(\check\mu_n,\check a_n)_{n \ge 0} \in \frak{S} \times \check \cA$.}
\KwResult{$\check Q_{N_{\mathrm{epi}}}$, an approximation of $\tilde{Q}^*$ on $\check{\frak{S}} \times \check \cA$.}
\Begin{
  Initialize table $\check Q_0  \in \RR^{|\check{\frak{S}}| \times |\check \cA|}$, $\mu_{0} \in \frak{S}$ and $a_{0} \in \cA$\;
  \For{ $n =0, 1, \dots N_{\mathrm{epi}}-1$}{
  \vskip 2pt
			Execute action $\check a_n$, observe $\check \mu'_{n+1} = \proj_{\check{\frak{S}}} \circ \bar F(\check \mu_n, \check \mu_n \measprod \check a_{n}, \varepsilon_{n+1}^0)$ and cost $\tilde f(\check \mu_{n}, \check a_{n})$ \;
			Initialize $\check Q_{n+1} = \check Q_n$ on $\check{\frak{S}} \times \check \cA$\;
			Set $\check Q_{n+1}(\check \mu_n , \check  a_n) = (1- \eta_n) \check Q_n( \check \mu_{n} , \check  a_{n}) + \eta_n  \left( \tilde f( \check \mu_{n}, \check a_{n}) + \gamma  \min_{\check a' \in \check \cA} \check  Q_n( \check \mu'_{n+1}, \check a' ) \right)$ \;
    }
  \KwRet{$\check Q_{N_{\mathrm{epi}}}$}
  }
\caption{Mean-Field Q-learning (MFQ) with simplex discretization}
\label{algo:Qtable-projection}
\end{algorithm}

Algorithm~\ref{algo:Qtable-projection} returns the table $\check Q_{N_{\mathrm{epi}}}$ after $N_{\mathrm{epi}}$ episodes. We prove below that this table converges to the optimal Q-function $\tilde Q^*$ in a suitable sense. 
To keep the paper at a reasonable length, we will make the following simplifying assumptions. 

We endow the simplexes $\frak{S}$ and $\frak{A}$ respectively with the Euclidean distances $d_{\frak{S}}$ and $d_{\frak{A}}$ of the spaces $\RR^{|S|}$ and $\RR^{|A|}$. Because $S$ is finite, we can identify $\cA$ defined in~\eqref{eq:def-calA-S-PA} with $\cP(A)^{|S|}$ and endow it with the distance $d_{\cA}(\tilde a, \tilde a') = \sup_{x \in S} d_{\frak{A}}( \tilde a(x), \tilde a'(x)) $ for $\tilde a, \tilde a' \in \cA$. Furthermore, we consider the following discretizations of the simplexes. Let $\varepsilon_{\frak{S}}>0$ satisfying: for all $\mu \in \frak{S},$ there exists $\check \mu \in \check{\frak{S}}$ s.t. $d_{\frak{S}}(\mu, \check \mu) \le \varepsilon_{\frak{S}}$.  Similarly, let $\varepsilon_{\frak{A}} > 0$ satisfying: for all $\nu \in \frak{A}$, there exists $\check \nu \in \check{\frak{A}}$ such that $d_{\frak{A}}( \nu, \check \nu ) \leq \varepsilon_{\frak{A}}$.  Because $S$ is finite and the definition of the distance $d_{\cA}$, we have for every $\tilde a \in \cA$, there exists $\check a \in \check \cA$, s.t. $d_{\cA}( \tilde a, \check a) \leq \varepsilon_{\frak{A}}$.

\begin{assumption}\label{hyp:bdd-smooth-data} 
\textbf{Regularity of the data:} $\tilde f$ is bounded and Lipschitz continuous with respect to $(\mu, \tilde a)$ with constant $L_{\tilde f}$, namely for every $(\mu, \tilde a), (\mu', \tilde a') \in \frak{S} \times \cA$, we have
	$$
		| \tilde f( \mu, \tilde a) - \tilde f( \mu' , \tilde a')  | \leq L_{\tilde f} \left( \| \mu - \mu' \|_{d_{\frak{S}}}  +  d_{\cA}( \tilde a , \tilde a')  \right) 
	\qquad
	and 
	\qquad
		\tilde f(\mu, \tilde a) \leq L_{\tilde f}.
	$$
	Also, $\bar F$ is Lipschitz continuous with respect to $\mu$ and $\tilde a$ with constant $L_{\bar F}$ in expectation over the randomness of the common noise, namely:  for every $(\mu , \tilde a), (\mu', \tilde a')  \in \frak{S} \times \cA$,
	\begin{align*}
		&\EE_{\varepsilon^0} \left[ \| \bar F(\mu, \mu \measprod \tilde a, \varepsilon^0) -  \bar F(\mu', \mu' \measprod \tilde a', \varepsilon^0)  \|_{d_{\frak{S}}} \right] 
		\le L_{\bar F}  \left( \| \mu - \mu' \|_{d_{\frak{S}}} + d_{\cA}( \tilde a , \tilde a') \right)
	\end{align*}
\end{assumption}
\begin{assumption}\label{hyp:smooth-V} 
\textbf{Regularity of the value function:} $\bar J^*$ is Lipschitz continuous w.r.t. $\mu$ with constant $L_{\bar J^*}$. 
\end{assumption}
	
\begin{assumption}\label{hyp:covering-time} 
\textbf{Covering time:} There exists a finite $T_{cov}$ such that with probability $1/2$ (over the randomness of the common noise and of Algorithm~\ref{algo:Qtable-projection})  the following holds: For every starting point in  $\check{\frak{S}} \times \check \cA$, every element of $\check{\frak{S}} \times \check \cA$ has been visited before time $T_{cov}$ during the execution of Algorithm~\ref{algo:Qtable-projection}. 
\end{assumption}

The regularity of $\bar J^*$ in~\ref{hyp:smooth-V} can typically be ensured through suitable conditions on the data of the problem, as e.g. in~\cite{ChassagneuxCrisanDelarue_Master,MR3967062,CarmonaDelarue_book_II}. Assumption~\ref{hyp:covering-time} is similar to the covering time assumption in~\cite{MR2247972}. In practice, exploration can be enhanced by adjusting the greediness level and by using exploring starts (if the learner can query an oracle which simulates transitions from any $(\mu, \tilde a)$). 
Note that the boundedness of the one-stage cost $\tilde f$ from Assumption~\ref{hyp:bdd-smooth-data} together with the fact that $\gamma \in (0,1)$ ensures the existence of a finite bound $\check J_{bound}$ for the state value function of the projected MFC problem. We denote by $\beta = (1-\gamma)/2$ the horizon of the MDP corresponding to the projected MFC problem, and for $\delta \in (0,1)$, we let $T_{cov}(\delta) = \lceil{T_{cov} \log_2(1/(2\delta))\rceil}$. We consider projection operators $\proj_{\check{\frak{S}}}: \frak{S} \to \check{\frak{S}}$ and $\proj_{\check \cA}: \cA \to \check \cA$ such that $(\proj_{\check{\frak{S}}}(\mu), \proj_{\check \cA}(\tilde a) ) := (\check \mu, \check a)$ for every $(\mu, \tilde a) \in \frak{S} \times \cA$ where $(\check \mu, \check a)$ is the closest point (or one of the closest points, in case of equality) in $\check{\frak{S}} \times \check \cA$ with respect to $d_{\frak{S}}$ and $d_{\cA}$. Based on simplexes discretizations, this point  satisfies $\| \mu - \check \mu\|_{d_{\frak{S}}} \leq \varepsilon_{\frak{S}}$ and $d_{\cA}( \tilde a , \check a ) \leq \varepsilon_{\frak{A}}$.

\begin{theorem}
\label{thm:main-cv-tabular}
Let $\delta \in (0,1)$ and $\varepsilon >0$. Assume Assumptions~\ref{hyp:bdd-smooth-data}--\ref{hyp:covering-time} hold. Consider learning rates $(\eta_n)_n$ satisfying: There exists $\kappa \in (1/2,1)$ such that for every $(\check \mu, \check a) \in \check{\frak{S}} \times \check \cA$, $\eta_n := \eta_n(\check \mu, \check a) = 1/ \big ( 1 + C(n, \check \mu, \check a) \big)^\kappa$ for each $n \geq 0$, where $C(n, \check \mu, \check a)$ is the number of times up to $n$ that the pair $(\check \mu, \check a)$ has been visited in Algorithm~\ref{algo:Qtable-projection}. If the number of episodes $N_{\mathrm{epi}}$ is of order
\begin{equation}
\label{eq:LB-Nepi-tabular}
\Omega\left( 
	\left(\frac{(T_{cov}(\delta))^{1+3 \kappa} \check J_{bound}^2 \, \ln \left(|\check{\frak{S}}| \, |\check{\frak{A}}|^{|S|} \check J_{bound} / (2\delta \beta \varepsilon) \right)}{\beta^2 \varepsilon^2}\right)^{\frac{1}{\kappa}} 
	+ 
	\left(\frac{(T_{cov}(\delta))}{\beta} \ln \left( \frac{\check J_{bound}}{\varepsilon} \right)\right)^{\frac{1}{1-\kappa}} \right),
\end{equation}
then with probability $1-\delta$,  
for all $(\mu, \tilde a) \in \frak{S} \times \cA$, 
$$
	\left| \check Q_{N_{\mathrm{epi}}} \Big(\proj_{\check{\frak{S}}} (\mu), \proj_{\check \cA}(\tilde{a}) \Big) - \bar Q^*(\mu,  \mu \measprod \tilde a) \right| \le \varepsilon',
$$
where 
$ \displaystyle
	\varepsilon'
	= 
	\varepsilon 
	+ \left( \frac{\gamma}{1 - \gamma} L_{\bar J^*} + L_{\tilde f} + \gamma L_{\bar J^*} L_{\bar F} \right) \varepsilon_{\frak{S}} + \frac{1}{1- \gamma} \left( L_{\tilde f} + \gamma L_{\bar J^*} L_{\bar F} \right) \varepsilon_{\frak{A}}.
$
\end{theorem}

Note that $\varepsilon$ can be chosen as small as desired provided $N_{\mathrm{epi}}$ is large enough. 
The second and the third terms in the error $\varepsilon'$ are proportional to $\varepsilon_{\frak{S}}$ and $\varepsilon_{\frak{A}}$, which is somehow unavoidable in general due to the projection on the finite sets $\check{\frak{S}}$ and $\check{\frak{A}}$. However, this error vanishes as $\varepsilon_{\frak{S}} \to 0$ and $\varepsilon_{\frak{A}} \to 0$, i.e., as $\check{\frak{S}}$ and $\check{\frak{A}}$ are better and better approximations of $\cP(S)$ and $\cP(A)$ respectively.

We prove this result below. The proof can be summarized in the following three steps: 
\textbf{(1)} For $N_{\mathrm{epi}}$ large enough, we have 
			$
				\check Q_{N_{\mathrm{epi}}} \approx \check Q^*
			$ 
			on $\check{\frak{S}} \times \check \cA$;
\textbf{(2)} 
	        $
	            \check Q^* \approx \tilde Q^*
			$ 
			on $\check{\frak{S}} \times \check \cA$;
\textbf{(3)} For every $(\mu, \tilde a) \in \frak{S} \times \cA$,
			$
				 \tilde Q^*( \proj_{\check{\frak{S}}}(\mu), \proj_{\check \cA }(\tilde a ) ) \approx \tilde Q^*(\mu, \tilde a).
			$
The first step relies on standard Q-learning convergence results~\citep{MR2247972}, while the two other steps stem from the regularity assumptions and the approximation of $(\frak{S}, \cA)$ by $(\check{\frak{S}}, \check \cA)$.

\begin{proof}[Proof of Theorem~\ref{thm:main-cv-tabular}]

Recall that we denote by $\check J^*$ and $\check Q^*$ respectively the state value function and the state-action value function of the projected MFC problem defined by~\eqref{eq:generic-MFC-fctmeasure-reward-proj}--\eqref{eq:generic-MFC-fctmeasure-dyn-proj}.

	 We first note that, for every $(\mu, \tilde a) \in \frak{S} \times \cA$,
	\begin{align*}
	    & \left| \check Q_{N_{\mathrm{epi}}} \big( \proj_{\check{\frak{S}}}(\mu), \proj_{\check \cA}(\tilde a) \big)  - \tilde Q^* \big( \mu, \tilde a \big) \right| 
	    \\
	    \leq  & \left| \check Q_{N_{\mathrm{epi}}} \big( \proj_{\check{\frak{S}}}(\mu), \proj_{\check \cA}(\tilde a) \big) - \check Q^* \big( \proj_{\check{\frak{S}}}(\mu), \proj_{\check \cA}(\tilde a) \big) \right| 
	    \\
	    &\quad + \left|  \check Q^* \big( \proj_{\check{\frak{S}}}(\mu), \proj_{\check \cA}(\tilde a) \big) - \tilde Q^* \big(  \proj_{\check{\frak{S}}}(\mu), \proj_{\check \cA}(\tilde a) \big) \right| 
	    \\
	   &  \quad + \left| \tilde Q^* \big(  \proj_{\check{\frak{S}}}(\mu), \proj_{\check \cA}(\tilde a) \big) - \tilde Q^* \big( \mu, \tilde a \big)  \right|.
	\end{align*}
	We then split the proof into three steps, which consist in bounding from above each term in the right hand side.

\vskip 6pt
	\textbf{Step 1.} We first analyze the difference between $\check Q_{N_{\mathrm{epi}}}$ and $\check Q^*$. 
			This comes from standard convergence results on Q-learning for finite state-action spaces.  
			More precisely, under Assumptions~\ref{hyp:bdd-smooth-data} and \ref{hyp:covering-time}, with our choice of learning rates, and given that $N_{\mathrm{epi}}$ is of order~\eqref{eq:LB-Nepi-tabular}, we can apply Theorem~4 and Corollary~34 in~\cite{MR2247972} for asynchronous Q-learning and polynomial learning rates, and we obtain that, with probability at least $1-\delta$, 
			$$
				\|\check Q_{N_{\mathrm{epi}}} - \check Q^* \|_{\infty}= \sup_{ (\check \mu, \check a) \in \check{\frak{S}} \times \check \cA} \Big| \check Q_{N_{\mathrm{epi}}}(\check \mu, \check a) - \check Q^*(\check \mu, \check a) \Big| \le \varepsilon.
			$$

\textbf{Step 2. } We then turn our attention to the difference between $\check Q^*$ and $\tilde Q^*$. 
			The analysis amounts to say that the projection on $\check{\frak{S}}$ realized at each step does not perturb too much the value function.  
			Recall that for some given common noise $\varepsilon^0$, the operator $\check \Phi^{\varepsilon^0}: \check{\frak{S}} \times \check \cA \to \check{\frak{S}}$ is given by $\check \Phi^{\varepsilon^0}( \check \mu, \check a) = \proj_{\check{\frak{S}}} \circ \bar F ( \check \mu, \check \mu \measprod \check a, \varepsilon^0)$.
			Likewise, we denote the transition dynamic with $\bar F$ by a function $\Phi^{\varepsilon^0}: \frak{S} \times \cA \to \frak{S}$ such that:
			$$
				\Phi^{\varepsilon^0} (\mu, \tilde a) = \bar F( \mu, \mu \measprod \tilde a, \varepsilon^0), \qquad \forall (\mu, \tilde a) \in \frak{S} \times \cA.
			$$ 
			Let us start by noting that, for every $(\check \mu, \check a) \in \check{\frak{S}} \times \check \cA$,
			\begin{align*}
				& \left| \check Q^*(\check \mu, \check a) -  \tilde Q^*(\check \mu, \check a) \right| 
				\\
				& \le  \gamma    \EE \Bigg[ \Big| \check J^*(\check \Phi^{\varepsilon^0}(\check \mu, \check a)) - \bar J^*( \Phi^{\varepsilon^0}(\check \mu, \check a) ) \Big| \Bigg]
				\\
				& \le   \gamma   \EE \Bigg[  \Big|  \check J^*( \check \Phi^{\varepsilon^0}(\check \mu, \check a) ) - \bar J^*( \check \Phi^{\varepsilon^0}(\check \mu, \check a)) \Big|  + \Big| \bar J^* (\check \Phi^{\varepsilon^0}(\check \mu, \check a)) - \bar J^*( \Phi^{\varepsilon^0}(\check \mu, \check a)) \Big|  \Bigg]
				\\
				& \le \gamma   \EE \Bigg[  \Big|  \inf_{\check a' \in \check \cA} \check Q^* \left(\check \Phi^{\varepsilon^0}(\check \mu, \check a), \check a' \right) - \inf_{\tilde a' \in \cA } \tilde Q^* \left(\check \Phi^{\varepsilon^0}(\check \mu, \check a), \tilde a' \right) \Big| \Bigg]  +  \gamma L_{\bar J ^*} \EE \left[  \left\Vert \check \Phi^{\varepsilon^0}(\check \mu, \check a) -  \Phi^{\varepsilon^0}(\check \mu, \check a) \right\Vert_{d_{\frak{S}}} \right],
			\end{align*}
			where the last inequality holds by Lipschitz continuity of $\bar J^*$ on $\frak{S}$, see Assumption~\ref{hyp:smooth-V}. 
			
			The second term in the last inequality can be bounded using the simplex discretization properties and Assumption~\ref{hyp:bdd-smooth-data}: 
			\begin{align*}
			   \EE \left[  \| \check \Phi^{\varepsilon^0}(\check \mu, \check a) - \Phi^{\varepsilon^0}( \check \mu, \check a) \|_{d_{\frak{S}}}  \right]
 				& = \EE_{\varepsilon^0_1}\left[ \| \proj_{\check{\frak{S}}} \circ \bar F(\check \mu, \check a, \varepsilon^0) - \bar F(\check \mu, \check a, \varepsilon^0)  \|_{d_{\frak{S}}} \right] \leq \varepsilon_{\frak{S}}.
			\end{align*}
			 
    		For the first term, let $\check \mu' = \check \Phi^{\varepsilon^0}(\check \mu, \check a) \in \check{\frak{S}}$ to alleviate the notation, and let us consider $\check a^*_1 \in \check \cA$ and $\tilde a^*_2 \in \cA$ satisfying:
    		$
    		    \check Q^*(\check \mu', \check a^*_1) = \inf_{\check a' \in \check \cA} \check Q^* \left(\check \mu', \check a' \right)$
    		and
    		$
    	        \tilde Q^*(\check \mu', \tilde a^*_2) = \inf_{\tilde a' \in \cA } \tilde Q^* \left(\check \mu', \tilde a' \right).
    	    $
    	   The existence of $\check a^*_1$ and $\tilde a^*_2$ is guaranteed respectively by finiteness of  $\frak{S} \times \check \cA$ and by Lemma~\ref{le:Qstar_is_lsc}. 
    	   
    	   We observe that
    	    \begin{align*} 
	        & \check Q^*(\check \mu', \check a_1^*) - \tilde Q^*( \check \mu', \tilde a_2^*) 
    	        \\
    	        & = \Big( \check Q^*(\check \mu', \check a_1^*) - \check Q^*( \check \mu', \proj_{\check \cA} ( \tilde a_2^*) )  \Big) + \Big( \check Q^*( \check \mu', \proj_{\check \cA} ( \tilde a_2^*) ) - \tilde Q^*( \check \mu', \proj_{\check \cA}(\tilde a_2^*) ) \Big) 
    	        \\
    	        & \qquad + \Big( \tilde Q^*( \check \mu', \proj_{\check \cA}( \tilde a_2^*) ) - \tilde Q^*( \check \mu', \tilde a_2^*)  \Big)
    	        \\
    	        & \leq 0 +  \sup_{ (\check \mu, \check a) \in \check{\frak{S}} \times \check \cA} \left| (\check Q^* - \tilde Q^*)(\check \mu, \check a) \right| + \left( \tilde f (\check \mu', \proj_{\check \cA}(\tilde a^*_2) ) + \gamma \EE_{(\varepsilon^0)'} \left[ \bar J^*( \bar F( \check \mu', \proj_{\check \cA}(\tilde a_2^*), (\varepsilon^0 )' ) ) \right]  \right)
    	        \\
    	        & \qquad \qquad - \left( \tilde f (\check \mu', \tilde a^*_2 ) + \gamma \EE_{(\varepsilon^0)'} \left[ \bar J^*( \bar F( \check \mu', \tilde a_2^*, (\varepsilon^0)' ) ) \right] \right) 
    	        \\
    	        & \leq  \| \check Q^* - \tilde Q^* \|_{\infty}  + ( L_{\tilde f} + \gamma L_{\bar J^*} L_{\bar F} ) \varepsilon_{\frak{A}}.
    	    \end{align*}
    	   On the other hand, 
		 \begin{align*}
	       \check Q^*( \check \mu', \check a_1^*) -  \tilde Q^*(\check \mu', \tilde a_2^*)
    	        = - \Big( \tilde Q^*( \check \mu' ,  \tilde a_2^* ) - \tilde Q^*( \check \mu', \check a_1^*) \Big) - \Big( \tilde Q^*( \check \mu' , \check a_1^*) - \check Q^*( \check \mu', \check a_1^* ) \Big)
    	        \geq - \| \check Q^* - \tilde Q^* \|_{\infty}.
    	    \end{align*}
	    
Combining the above bounds yields that for every $(\check \mu, \check a) \in \check{\frak{S}} \times \check \cA$,
			$$
			    \left| \check Q^*(\check \mu, \check a) - \tilde Q^*(\check \mu, \check a) \right| \leq \gamma \left( \| \check Q^* - \tilde Q^* \|_{\infty} + (L_{\tilde f} + \gamma L_{\bar J^*} L_{\bar F} ) \varepsilon_{\frak{A}} \right) + \gamma L_{\bar J^*} \varepsilon_{\frak{S}}.  
			$$
			Consequently,
			$$
			\|\check Q^* - \tilde Q^* \|_\infty
			\le
			\frac{\gamma }{1-\gamma} \left( ( L_{\tilde f} + \gamma L_{\bar J^*} L_{\bar F} )\varepsilon_{\frak{A}} + L_{\bar J^*} \varepsilon_{\frak{S}} \right).
			$$

\textbf{Step 3. } Last, we look at the difference between $\tilde Q^*( \proj_{\check{\frak{S}}}(\mu), \proj_{\check \cA}(\tilde a )  )$ and $\tilde Q^*(\mu, \tilde a)$. 
			For every $\mu \in \frak{S}$ and $\tilde a \in \cA$, letting $\check \mu = \proj_{\check{\frak{S}}} (\mu)$ and $\check a = \proj_{\check \cA}(\tilde a)$ to alleviate the notation, we have $\| \check \mu - \mu \|_{d_{\frak{S}}} \leq \varepsilon_{\frak{S}}$ and $\| \check a - \tilde a\|_{d_{\cA}} \leq \varepsilon_{\frak{A}}$. We obtain
			\begin{align*}
				\left| \tilde Q^* (\check \mu, \check a) -  \tilde Q^*(\mu, \tilde a) \right|
				& \le \Bigg|  \tilde f(\check \mu, \check a)  -  \tilde f(\mu, \tilde a) \Bigg| +  \gamma  \EE \Bigg[  \left| \inf_{\tilde a' \in \cA} \tilde Q^* ( \Phi(\check\mu, \check a ), \tilde a') - \inf_{\tilde a' \in \cA}  \tilde Q^*(\Phi(\mu, \tilde a), \tilde a') \right| \Bigg]
				\\
				& \le L_{\tilde f} \left( \| \check \mu - \mu \|_{d_{\frak{S}}} +  \| \check a - \tilde a\|_{d_{\cA}} \right)
				 +  \gamma  \EE \left[  \left| \bar J^* ( \bar F ( \check \mu, \check a, \varepsilon^0 ) - \bar J^*( \bar F (\check \mu, \check a, \varepsilon^0 ) ) \right| \right] 
				\\
				& \le L_{\tilde f} (\varepsilon_{\frak{S}} + \varepsilon_{\frak{A}} ) + \gamma L_{\bar J^*} \EE \left[ \| \bar F(\check \mu, \check a, \varepsilon^0 ) - \bar F( \mu, \tilde a, \varepsilon^0) \|_{d_{\frak{S}}} \right]
				\\
				& \le  ( L_{\tilde f} + \gamma L_{\bar J^*} L_{\bar F} )(\varepsilon_{\frak{S}} + \varepsilon_{\frak{A}} ),
			\end{align*}
			where we used the Lipschitz continuity of $\tilde f, \bar J^*, \bar F$ and the assumption on $\check{\frak{S}}$, see Assumptions~\ref{hyp:bdd-smooth-data}, \ref{hyp:smooth-V} and the simplex discretization properties.

\end{proof}

\subsubsection{Q-learning with controls that are pure at level-0}
The above method is designed for the case where one looks for optimal actions that are potentially randomized at the individual level. Searching in the space $\cP(A)$ comes with a computational cost that is reflected in the bounds through the cardinality of the discrete simplex $\check{\frak{A}}$. In some situations it can be interesting to directly search for actions that are pure at the individual level. 

In this case, instead of~\eqref{eq:def-calA-S-PA}, the set of strategy functions is (for simplicity we keep the notation $\cA$):
 \begin{equation*}
		\cA := \{ \tilde a: S \to A \} = A^S.
\end{equation*}
In Algorithm~\ref{algo:Qtable-projection}, we replace $\check a \in \check{\cA}$ by $\tilde a \in \cA$. 

\begin{theorem}
\label{thm:main-cv-tabular-pure}
Let $\delta \in (0,1)$ and $\varepsilon >0$. Assume Assumptions~\ref{hyp:bdd-smooth-data}--\ref{hyp:covering-time} hold. Consider learning rates $(\eta_n)_n$ satisfying: There exists $\kappa \in (1/2,1)$ such that for every $(\check \mu, \tilde a) \in \check{\frak{S}} \times \tilde \cA$, $\eta_n := \eta_n(\check \mu, \tilde a) = 1/ \big ( 1 + C(n, \check \mu, \tilde a) \big)^\kappa$ for each $n \geq 0$, where $C(n, \check \mu, \tilde a)$ is the number of times up to $n$ that the pair $(\check \mu, \tilde a)$ has been visited in Algorithm~\ref{algo:Qtable-projection}. If the number of episodes $N_{\mathrm{epi}}$ is of order
\begin{equation}
\label{eq:LB-Nepi-tabular-pure}
\Omega\left( 
	\left(\frac{(T_{cov}(\delta))^{1+3 \kappa} \check J_{bound}^2 \, \ln \left(|\check{\frak{S}}| \, |A|^{|S|} \check J_{bound} / (2\delta \beta \varepsilon) \right)}{\beta^2 \varepsilon^2}\right)^{\frac{1}{\kappa}} 
	+ 
	\left(\frac{(T_{cov}(\delta))}{\beta} \ln \left( \frac{\check J_{bound}}{\varepsilon} \right)\right)^{\frac{1}{1-\kappa}} \right),
\end{equation}
then with probability $1-\delta$,  
for all $(\mu, \tilde a) \in \frak{S} \times \cA$, 
$$
	\left| \check Q_{N_{\mathrm{epi}}} \Big(\proj_{\check{\frak{S}}} (\mu), \tilde{a} \Big) - \bar Q^*(\mu,  \mu \measprod \tilde a) \right| \le \varepsilon',
$$
where 
$ \displaystyle
	\varepsilon'
	= 
	\varepsilon 
	+ \left( \frac{\gamma}{1 - \gamma} L_{\bar J^*} + L_{\tilde f} + \gamma L_{\bar J^*} L_{\bar F} \right) \varepsilon_{\frak{S}}.
$
\end{theorem}

The above result provides convergence guarantee for the Q-function. Let us now derive a consequence in terms of the optimizer. To this end, we will use the following additional assumption on the gap between the values of the best and second-best actions, which is rather standard in approximation algorithms based on tabular Q-functions~\citep{farahmand2011action,bellemare2016increasing}. 

\begin{assumption}\label{hyp:action-gap} \textbf{Action gap:} There exists $K_A > 0$ such that: 
	$$
		\tilde Q^*( \check \mu, \tilde a) - \inf_{\tilde a' \in \cA} \tilde Q^*( \check \mu, \tilde a' )  \ge K_A, \quad \check \mu \in \check{\frak{S}}, \tilde a \in  \cA \backslash \arg\inf_{\tilde a' \in \cA} \tilde Q^*( \check \mu, \tilde a').
	$$
\end{assumption}
To recover minimizers or approximate minimizers, it will be convenient to work with the following operators. In general, they are defined on the vector space $\RR^{\texttt{m}}$. For $\tau>0$ and $x = (x_1,\dots,x_{\texttt{m}}) \in \RR^{\texttt{m}}$, we define $\softmin_\tau : \RR^{\texttt{m}} \to \RR^{\texttt{m}}$ by
$$
	\softmin_\tau(x) = (e^{- \tau x_1}, \dots, e^{- \tau x_{\texttt{m}}}) / \sum_{j} e^{- \tau x_j}.
$$
For $x \in \RR^{\texttt{m}}$, we define $ \argmine: \RR^{\texttt{m}} \to [0,1]^{\texttt{m}}$ by
$$
	\argmine(x) = \left( \mathbf{1}_{i \in \argmin(x)} \right)_{i=1}^{\texttt{m}} / |\argmin(x)|.
$$
where $\argmin(x) = \{ j \in \{1,\ldots, {\texttt{m}} \} \, : \,   x_j = \min\{x_1,\ldots, x_{\texttt{m}} \} \}$.
In the sequel, we use these operators with the dimension ${\texttt{m}} = | \cA | = | A |^{|S|}$. For any function $q:   \cA \to \RR$, we identify $q$ with the vector $(q(\tilde a))_{\tilde a \in \cA}$.

\begin{corollary}
\label{coro:action-cv-tabular-pure}
Assume the same assumptions as in Theorem~\ref{thm:main-cv-tabular-pure} hold and, in addition, that Assumption~\ref{hyp:action-gap} holds. Let $\check Q_{N_{\mathrm{epi}}}$ be the table returned by Algorithm~\ref{algo:Qtable-projection}, and let $\varepsilon'$ be as in Theorem~\ref{thm:main-cv-tabular-pure}. Then for every $\check\mu \in \check{\frak{S}}$,
\begin{equation*}
	\big\| \softmin_\tau \big( \check Q_{N_{\mathrm{epi}}}( \check\mu, \cdot) \big)	-
	\argmine \big(  \tilde Q^* (\check\mu, \cdot) \big) \big\|_2
	\leq 
	 \tau \varepsilon' \sqrt{ | \cA | } + 2  e^{- \tau K_A}  |\cA|.
\end{equation*} 
\end{corollary}

The proof is provided below. The $\argmine$ in the second term is here in case there are several optimal controls. The $\softmin$ regularizes the best action predicted by the estimation $\check Q_{N_{\mathrm{epi}}}$ of the function $\tilde Q^*$.

\begin{remark}
\label{rem:action-gap-issue}
Imagine we want the error bound in Corollary~\ref{coro:action-cv-tabular-pure}, 
to be smaller than some $\delta >0$. It is sufficient to have: for the second term: $\displaystyle \tau \ge \frac{1}{K_A} \log \big( \frac{|\cA|}{\delta / 4} \big)$; and for the first term:  $\displaystyle \varepsilon' \le  \delta / (2\tau \sqrt{ | \cA | }) = \delta K_A / \big(2 |\cA|^{1/2} \log \big( \frac{ |\cA| }{\delta /4} \big) \big) $.  Then both terms in the error bound will be smaller than $\delta / 2$.
Notice that, contrary to Theorem~\ref{thm:main-cv-tabular}, here we do not need to approximate the probability space of action $\cP(A)$ by $\check{\frak{A}}$ with an $\varepsilon_{\frak{A}}$-net, hence the error bound in Theorem~\ref{thm:main-cv-tabular-pure} is independent of any  $\varepsilon_{\frak{A}}$. 
So it is possible to choose $\tau$ and to make $\varepsilon'$ as small as we want. 
\end{remark}

\begin{proof}[Proof of Corollary~\ref{coro:action-cv-tabular-pure}]
We use~\cite[Proposition~4]{gao2017properties},  which states that $\softmin_\tau$ is $\tau$-Lipschitz and~\cite[Lemma~7]{guo2019learningMFG}, which states that for $(x_i)_{i =1, \dots, {\texttt{m}}}$,
$$
	\|\softmin_\tau\, (x) - \argmine(x)\|_2 \le 2 {\texttt{m}} e^{- \tau \delta},
$$
where $\delta = \inf_{x_j > \inf(x)} x_j - \inf(x)$, and $\delta = \infty$ if all $x_i$ are equal. We can apply this latter result to $\tilde Q^*(\check \mu, \cdot)$ thanks to assumption~\ref{hyp:action-gap}, with ${\texttt{m}} = |\check \cA|$ and $\delta = K_A$. Combining this with Theorem~\ref{thm:main-cv-tabular},  we have, for every $\check \mu$,
\begin{align*}
	&\bigg\| \softmin_\tau \Big( \check Q^*(\check\mu, \cdot) \Big) - \argmine \Big( \tilde Q^*( \check \mu, \cdot)  \Big) \bigg\|_2
	\\
	&\le
	\bigg\| \softmin_\tau \Big(  \check Q^*(\check\mu, \cdot) \Big) - \softmin_\tau \Big(  \tilde Q^*(\check \mu, \cdot) \Big) \bigg\|_2
	+ \bigg\| \softmin_\tau \Big( \tilde Q^*(\check\mu, \cdot) \Big)  - \argmine \Big( \tilde Q^*(\check \mu, \cdot) \Big) \bigg\|_2
	\\
	&\le
	\tau \bigg\| \check Q^*(\check\mu, \cdot) - \tilde Q^* (\check \mu, \cdot) \bigg\|_2 + 2 |\cA| e^{- \tau K_A}
	\\
	& \le \tau  \sqrt{ | \cA | } \sup_{ \tilde a' \in \cA} \left| \check Q^*(\check\mu, \tilde a') - \tilde Q^* (\check \mu, \tilde a') \right| + 2 |\cA| e^{- \tau K_A}
	\\
	&\le
	\tau \varepsilon'  \sqrt{ | \cA | } + 2  e^{- \tau K_A} | \cA |.
\end{align*}

\end{proof}

\subsection{Deep reinforcement learning for MFMDP }\,

The above method has the advantage to be simple enough to let us carry out a detailed analysis. However, it cannot be used in practice for large state or actions spaces because of the prohibitive computational cost due to the discretization of the simplexes. An alternative is to work directly with continuous spaces, in which case the policies and value functions cannot be represented in a tabular way. Instead, we can rely on function approximation. 
To this end, we now propose to use methods from deep reinforcement learning which are more suitable for continuous spaces. The motivations are twofold. 

First, if $S$ and $A$ are finite but we want to learn an optimal policy that is potentially randomized at level-0, the discretization approach proposed in~\S~\ref{sec:tabularQ-mixed} has a complexity that increases with the number of points in the discretization of $\cP(A)$, which itself increases exponentially quickly with the cardinality of $A$. For this reason, it can be interesting to tackle directly $\bar A = \cP(A)$ as a continuous action space and to use deep RL methods for continuous action space MDPs.

Second, some MFC problems are naturally posed with a continuous state space $S$. In this case, under mild conditions, the optimal policy is in fact non-randomized not only at the level-1 but even at the level-0. 
However, the state of the mean field MDP is an element of the infinite dimensional space $\cP(S)$. From a numerical viewpoint, we need two ingredients: (1) a finite-dimensional approximation of the mean field and (2) a parameterized approximation of the value function or the policy taking this finite-dimensional representation of the mean field state as an input. For the second point, we will again use deep neural networks. For the first point, for the sake of definiteness, we choose to simply replace $\cP(S)$ by $\cP(\check S)$ where $\check S$ is a discretization of $S$ with a finite number of points. We assume that, given $\check\mu \in \cP(\check S)$ and $\bar a: S \to A$, one can get from the environment a sample of the next state and the associated cost $\tilde f(\check\mu, \bar a)$. The problem thus boils down to an MDP with finite dimensional (but potentially continuous) state and action spaces. Such MDPs can be solved with a variety of deep RL algorithms. In the sequel, we provide numerical illustrations based on the Deep Deterministic Policy Gradient (DDPG) proposed in~\cite{lillicrap2015continuous}. It relies on two neural networks, one for the Q-function (the critic) and one for the policy (the actor). The heart of the algorithm consists in updating alternatively the critic by minimizing an empirical square error and the actor by making one step of gradient descent. 
To improve exploration, a Gaussian noise $\epsilon^a_{n+1}$ is added to the action prescribed by the actor. Furthermore, for more stability, target networks are also added. The algorithm is summarized in our setting in Algorithm~\ref{algo:DDPG-MFC} in the Appendix~\ref{app:details-Qlearning}.

\section{\textbf{Numerical Examples}}
\label{sec:numres}

\subsection{Example 1: Cyber security model}\,

For a first testbed, we start with a finite state problem. 
We revisit the cyber security example introduced in~\cite{MR3575619}, but here from the point of view of a central planner (such as a large company or a state) trying to protect its computers against the attacks of a hacker. The situation can thus be phrased as a MFC problem.

In this model, the population consists of a large group of computers which can be either defended (D) or undefended (U), and either infected (I) or susceptible (S) of infection. Hence the set $S$ has four elements corresponding to the four possible combinations: DI, DS, UI, US. The action set is $A = \{0,1\}$, where $0$ is interpreted as the fact that the central planner is satisfied with the current level of protection (D or U) of the computer under consideration, whereas $1$ means that she wants to change this level of protection. In the latter case, the update occurs at a (fixed) rate $\lambda >0$. If the controls are pure at level-0, at each of the four states, the central planner only chooses one action per state and applies it to all the computers at that state. If the controls are mixed at level-0, then for each state, she chooses a distribution over actions and then each computer in this state picks independently an action according to the chosen distribution. 
When infected, each computer may recover at rate $q_{rec}^D$ or $q_{rec}^U$ depending on whether it is defended or not. On the other hand, a computer may be infected either directly by a hacker, at rate $v_H q_{inf}^D$ (resp. $v_H q_{inf}^U$) if it is defended (resp. undefended), or by undefended infected computers, at rate $\beta_{UU}\mu(\{UI\})$ (resp. $\beta_{UD}\mu(\{UI\})$) if it is undefended (resp. defended), or by defended infected computers, at rate $\beta_{DU}\mu(\{DI\})$ (resp. $\beta_{DD}\mu(\{DI\})$) if it is undefended (resp. defended). Here $v_H$ can be interpreted as the attack intensity parameter.

In short, the transition matrix is given by: 
\begin{equation}
    \label{eq:cybersecurity-MFC-transition-matrix}
	P^{\mu,a} = 
	\begin{pmatrix}
	\dots 	& 		P^{\mu,a}_{DS \rightarrow DI}	&	 \lambda a 	&	0
	\\
	q_{rec}^D 	& 	\dots 		&	 0	&	\lambda a
	\\
	\lambda a 	& 	0 		&	 \dots	&	P^{\mu,a}_{US \rightarrow UI}
	\\
	0	&	\lambda a	&	q_{rec}^U	& \dots
	\end{pmatrix}
\end{equation}
where 
\begin{align*}
	&P^{\mu,a}_{DS \rightarrow DI} = v_H q_{inf}^D + \beta_{DD} \mu(\{DI\})  + \beta_{UD} \mu(\{UI\}) ,
	\\
	&P^{\mu,a}_{US \rightarrow UI} = v_H q_{inf}^U + \beta_{UU} \mu(\{UI\}) + \beta_{DU} \mu(\{DI\}),
\end{align*}
and all the instances of $\dots$ should be replaced by the negative of the sum of the entries of the row in which $\dots$ appears on the diagonal. At each time step, the central planner pays a protection cost $k_D>0$ for each defended computer, and a penalty $k_I>0$ for each infected computer. The instantaneous cost in the MFMDP is thus defined as: 
$$
	\bar f(\mu, \bar a) = k_D \mu(\{DI, DS\}) + k_I \mu(\{DI, UI\}), \qquad (\mu, \bar a) \in \bar S \times \bar A.
$$
The optimal control and optimal flow of distributions can be characterized by a forward-backward ODE system which can be obtained in way similar to what is done in the MFG setting e.g. in ~\cite[\S~7.2.3]{CarmonaDelarue_book_I}. We will use this solution as a benchmark. 

\textbf{Tabular Q-learning.} For the sake of illustration, we present results obtained by tabular Q-learning with simplex discretization as described in \S~\ref{sec:tabularQ-discrete}. The state space for the population distribution is $\bar S$, which is identified with the simplex $\frak{S} = \{(\mu^{(i)})_{i=1,\dots,4} \in [0,1]^4\,:\, \sum_{i} \mu^{(i)} = 1\}$. To follow the original setting considered in~\cite{MR3575619}, we consider pure controls, both at the common and idiosyncratic levels (level-0 and level-1). So we identify $\bar A$ with the set of functions $A^S$, which is finite and of cardinality $2^4 = 16$. 

We replace $\frak{S}$ by the finite set:
$$
	\check{\frak{S}} = \Big\{ (\mu^{(i)})_{i=1,\dots,4} \in [0,1/N_m,\dots,1-1/N_m, 1]^4\,:\, \sum_{i} \mu^{(i)} = 1 \Big\},
$$
where $[0,1/N_m,\dots,1-1/N_m, 1]$ is a uniform grid over $[0,1]$ with $N_m+1 \ge 2$ points. We then aim at computing the Q-function for the projected MDP with finite state space $\check{\frak{S}}$ and action space $A^S$, that we still denote by $\tilde{Q}^*$ although we do not consider mixed actions at the level-0. We note that, in the absence of common noise, the MFMDP is completely deterministic hence it would be enough to query once each state-action pair from the environment in order to learn the level-1 reward function and transition function, and hence to be able to compute perfectly the Q-function. However, for the sake of illustration, we stick to applying Algorithm~\ref{algo:Qtable-projection}, replacing both $\cA$ and $\check \cA$ by $A^S$. 

In order to be able to compare with the benchmark solution obtained by the ODE method, we considerthat the time steps are of size not $1$ but $\Delta t = 0.01$. Although the problem is set on an infinite horizon, we truncate the training episodes and the plots at the horizon $T=10$. 

\begin{figure}[h]
	\begin{subfigure}{.45\columnwidth}
		\centering
		\includegraphics[width=\columnwidth]{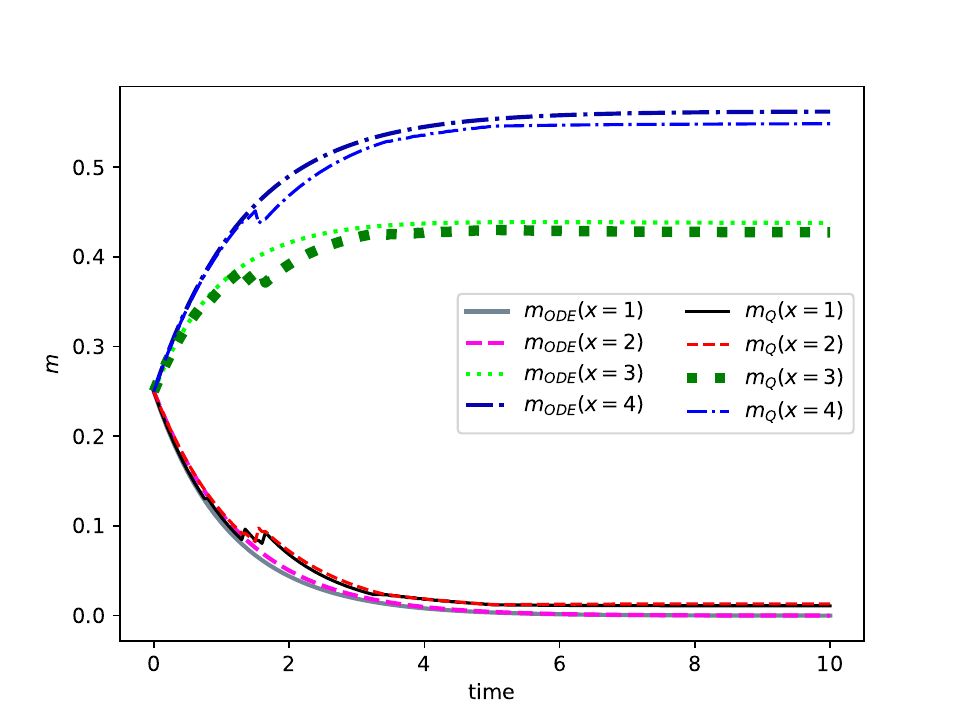}
	\end{subfigure}%
	\begin{subfigure}{.45\columnwidth}
		\centering 
		\includegraphics[width=\columnwidth]{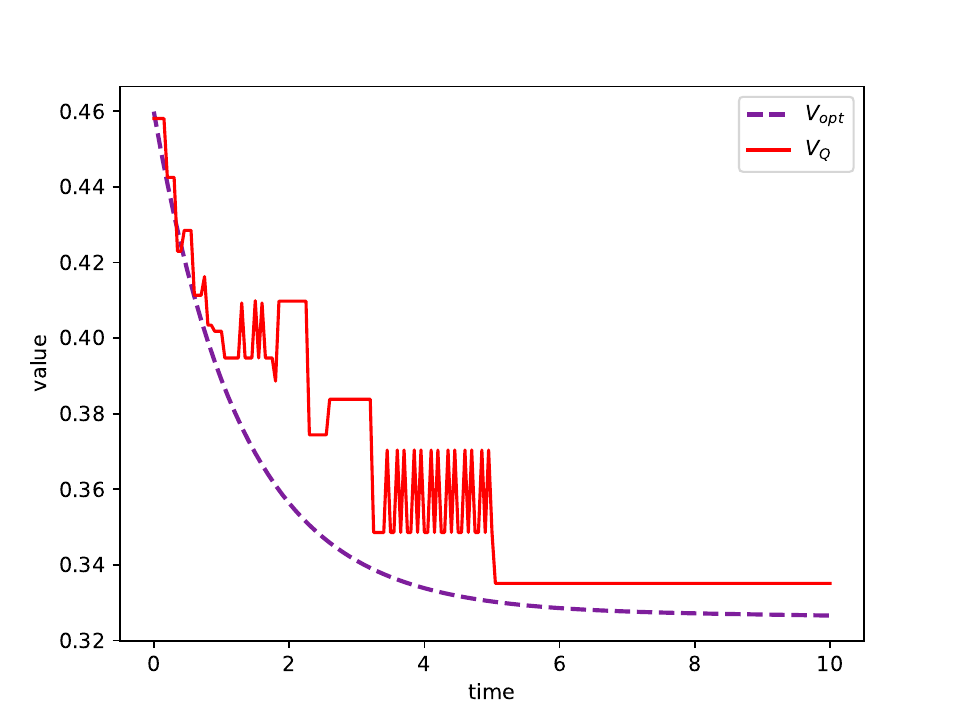}
	\end{subfigure}
	\caption{Example 1: Cyber security model. Test case 1: $m_0 = (1/4, 1/4, 1/4, 1/4)$. Left: Evolution of the distribution when using the benchmark optimal control ($m_{ODE}$) or the control recovered from the learnt Q-function ($m_Q$). Right: state value function using the benchmark solution ($V_{opt}$) or the learnt Q-function ($V_{Q}$) along the optimal mean field flow. The benchmark solution is obtained using the ODE method. }    
 	\label{AMS-num-fig:finiteMFCQ-cyber-Master-m0-1}
\end{figure}

\begin{figure}[h]
	\begin{subfigure}{.45\columnwidth}
		\centering
		\includegraphics[width=\columnwidth]{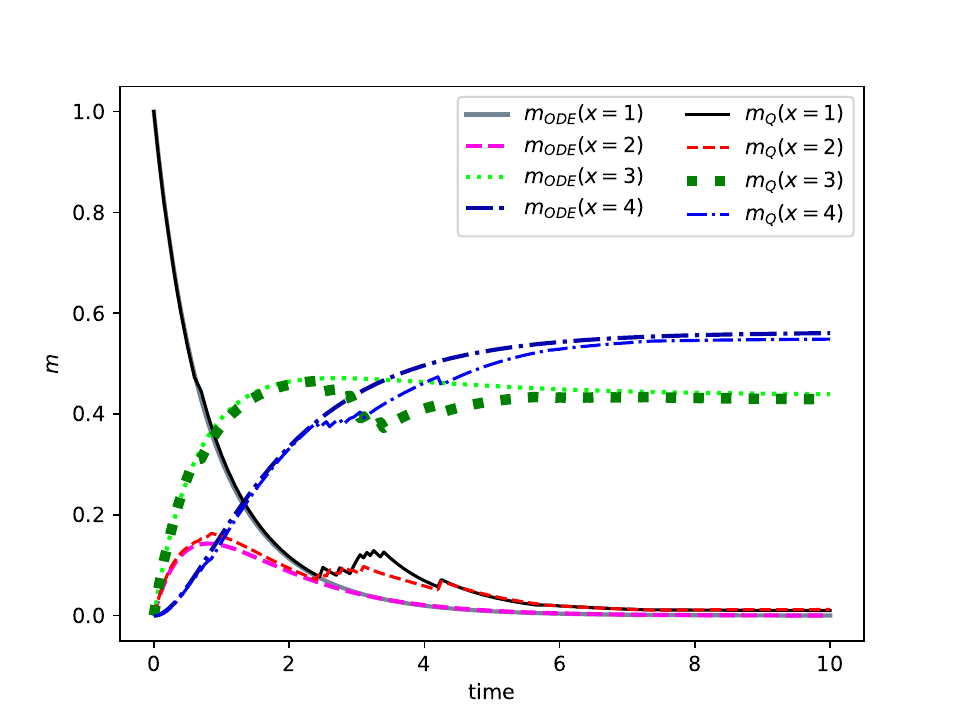}
	\end{subfigure}%
	\begin{subfigure}{.45\columnwidth}
		\centering 
		\includegraphics[width=\columnwidth]{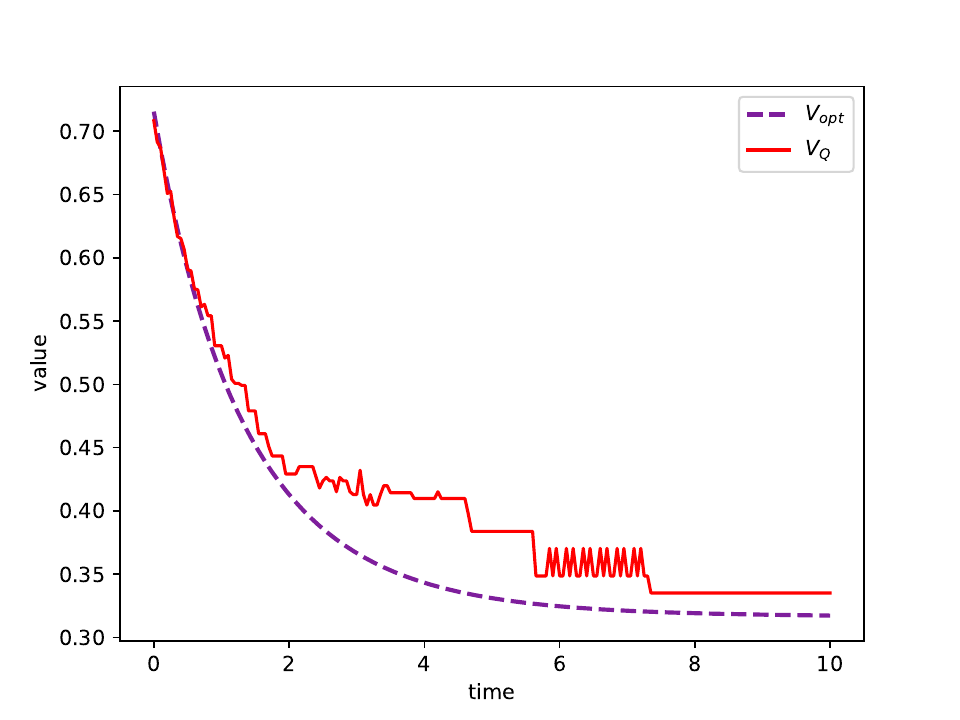}
	\end{subfigure}
	\caption{Example 1: Cyber security model. Test case 2: $m_0 = (1, 0, 0, 0)$. Left: Evolution of the distribution when using the benchmark optimal control ($m_{ODE}$) or the control recovered from the learnt Q-function ($m_Q$). Right: state value function using the benchmark solution ($V_{opt}$) or the learnt Q-function ($V_{Q}$) along the optimal mean field flow. The benchmark solution is obtained using the ODE method. }   
 	\label{AMS-num-fig:finiteMFCQ-cyber-Master-m0-2}
\end{figure}

\begin{figure}[h]
	\begin{subfigure}{.45\columnwidth}
		\centering
		\includegraphics[width=\columnwidth]{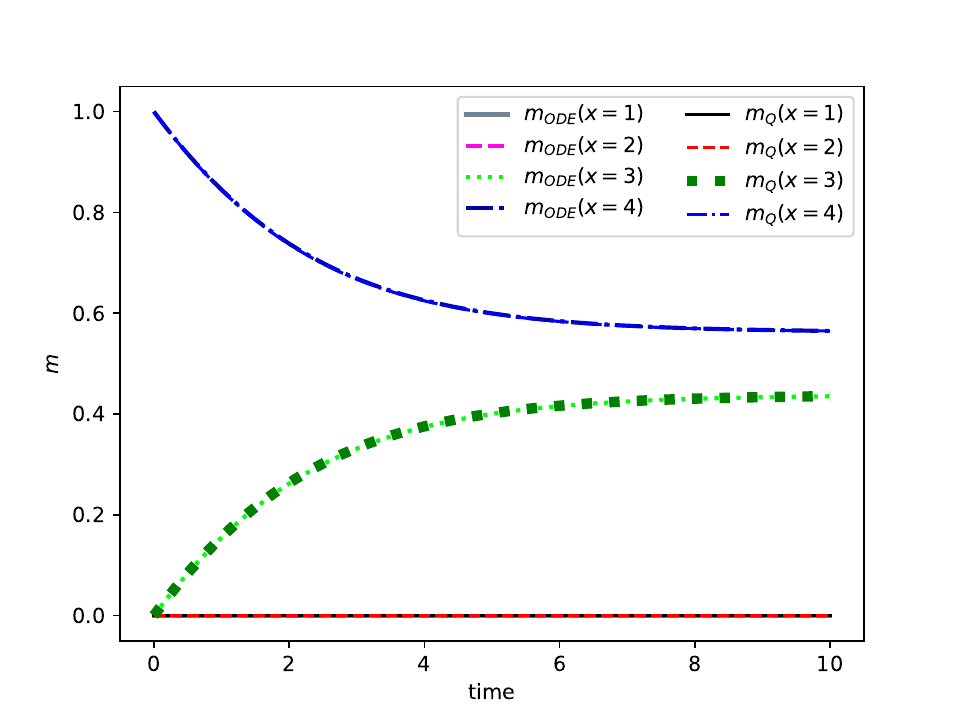}
	\end{subfigure}%
	\begin{subfigure}{.45\columnwidth}
		\centering 
		\includegraphics[width=\columnwidth]{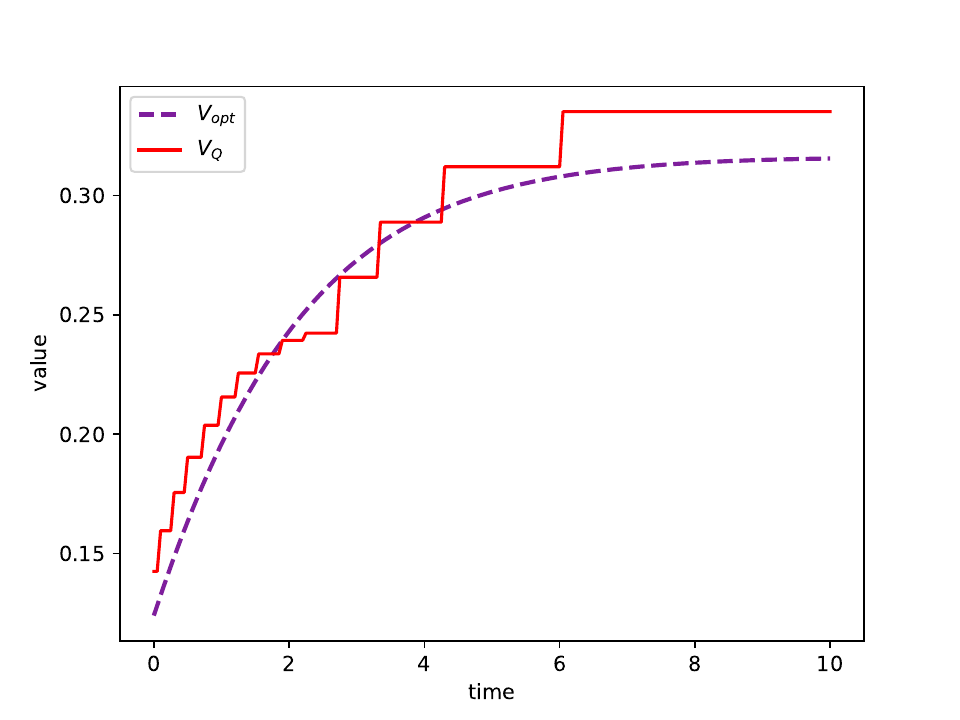}
	\end{subfigure}
	\caption{Example 1: Cyber security model. Test case 3: $m_0 = (0, 0, 0, 1)$. Left: Evolution of the distribution when using the benchmark optimal control ($m_{ODE}$) or the control recovered from the learnt Q-function ($m_Q$). Right: state value function using the benchmark solution ($V_{opt}$) or the learnt Q-function ($V_{Q}$) along the optimal mean field flow. The benchmark solution is obtained using the ODE method. }   
 	\label{AMS-num-fig:finiteMFCQ-cyber-Master-m0-3}
\end{figure}

After $N_{\mathrm{epi}}$ episodes of Q-learning, we obtain an approximation $\check Q_{N_{\mathrm{epi}}}$ of the Q-function, from which we can recover an approximation $\bar a_{N_{\mathrm{epi}}}$ of the optimal control by taking the argmax, namely: $\bar a_{N_{\mathrm{epi}}}(\mu, \cdot) = \argmax_{\check a \in A^S} \check Q_{N_{\mathrm{epi}}}( \proj_{\check{\frak{S}}}(\mu), \check a)$. 
 We compare the flow of distributions induced by this control $\bar a_{N_{\mathrm{epi}}}$ with the optimally controlled flow computed by the ODE method.  This method also allows us to compute for each $t$ the value $\bar J^{*}(\mu^*_t)$ along the optimal flow $(\mu^*_t)_{t\in[0,T]}$, in line with Lemma~\ref{lemma:infimum_action_opt_barQ_and_opt_barJ}, we compare it with $\max_{A^S} \check Q_{N_{\mathrm{epi}}}(\proj_{\check{\frak{S}}} (\mu^*_t), \cdot) = \check Q_{N_{\mathrm{epi}}}(\proj_{\check{\frak{S}}} (\mu^*_t), \bar a_{N_{\mathrm{epi}}}(\mu^*_t))$. Figures~\ref{AMS-num-fig:finiteMFCQ-cyber-Master-m0-1}--\ref{AMS-num-fig:finiteMFCQ-cyber-Master-m0-3} show the results for three initial conditions $\mu_0$. We see that the learnt value function approximately matches the $\bar J^{*}$ value function, and the induced flows of distributions approximately match the benchmark ones.    
For these simulations, we used $N_m = 30$, $\gamma = 0.5$, and the following parameters:
\begin{equation*}
\left\{
\begin{split}
&\beta_{UU} = 0.3,
\beta_{UD} = 0.4,
\beta_{DU} = 0.3,
\beta_{DD} = 0.4,
\\
&q_{rec}^D = 0.5, 
q_{rec}^U = 0.4, 
q_{inf}^D = 0.4, 
q_{inf}^U = 0.3, 
\\
&v_H = 0.6,
\lambda = 0.8,
k_D = 0.3, 
k_I = 0.5.
\end{split}
\right.
\end{equation*}

\textbf{Deep Deterministic Policy Gradient.}  The solution can also be learnt by using the deep RL algorithm given in Algorithm~\ref{algo:DDPG-MFC} instead of mean field Q-learning given in Algorithm~\ref{algo:Qtable-projection}. In the present case, this approach has the advantage of avoiding the  discretization of $\cP(\cS)$ since we instead directly deal with the distribution as a vector in dimension $4$. Since, with this method, it is possible to allow the control to take continuous values, we replace $A=\{0,1\}$ by $A =[0,1]$.  

Furthermore, to make things more interesting, we consider that the attack intensity parameter $v_H$ is stochastic. Since its value affects the evolution of the whole population, we model this using a common noise. We replace the state distribution by a conditional state distribution, conditioned on the realization of $\epsilon^0$ up to the current time.  To wit, let $(\epsilon^0_n)_{n \ge 1}$ be a sequence of i.i.d. random variables with Gaussian distribution. Let $v_{H,n+1} = v_{H,n} + \epsilon^0_{n+1}$, $n \ge 0$, $v_{H,0}$ given. The evolution from time $n$ to time $n+1$ of the state distribution is by the transition matrix defined in~\eqref{eq:cybersecurity-MFC-transition-matrix} but with the constant $P^{\mu,a}_{DS \rightarrow DI}$, $P^{\mu,a}_{US \rightarrow UI}$ replaced by the following stochastic coefficients that evolve in time due to the fact that $v_H$ is replaced by a stochastic process:
\begin{align*}
	&P^{\mu,a}_{DS \rightarrow DI, n} = v_{H,n} q_{inf}^D + \beta_{DD} \mu(\{DI\})  + \beta_{UD} \mu(\{UI\}) ,
	\\
	&P^{\mu,a}_{US \rightarrow UI, n} = v_{H,n} q_{inf}^U + \beta_{UU} \mu(\{UI\}) + \beta_{DU} \mu(\{DI\}),
\end{align*}

\begin{figure*}
	\centering
	\begin{subfigure}{.3\textwidth}
		\centering\captionsetup{width=.8\linewidth}
		\includegraphics[width=\textwidth]{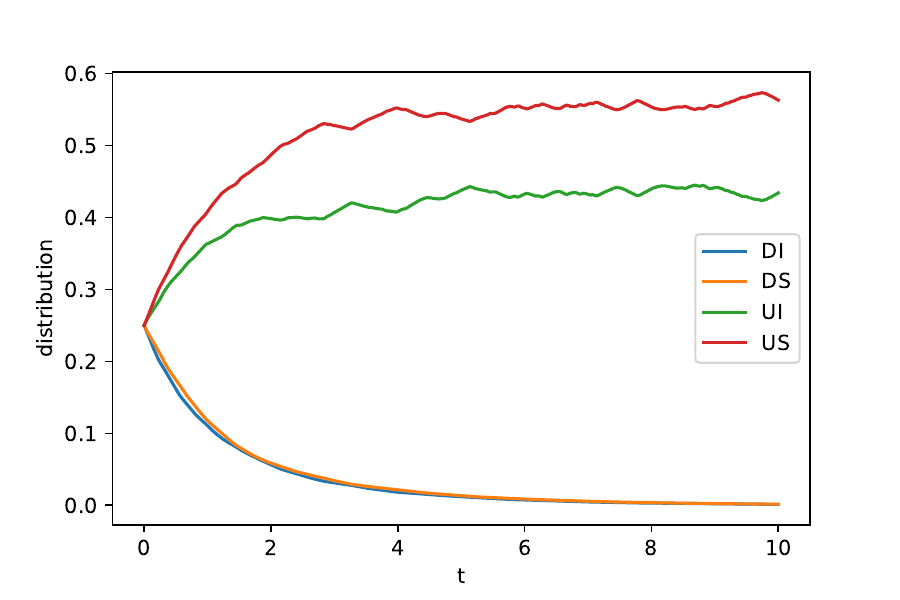}
		\caption{Test distribution 1}
	\end{subfigure}%
	\begin{subfigure}{.3\textwidth}
		\centering\captionsetup{width=.8\linewidth}
		\includegraphics[width=\textwidth]{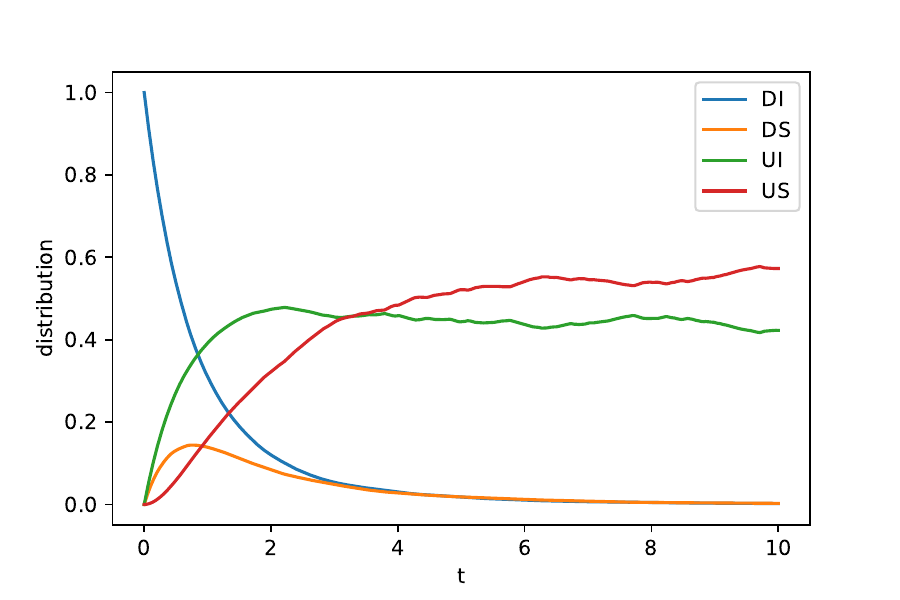}
		\caption{Test distribution 2}
	\end{subfigure}%
	\begin{subfigure}{.3\textwidth}
		\centering\captionsetup{width=.8\linewidth}
		\includegraphics[width=\textwidth]{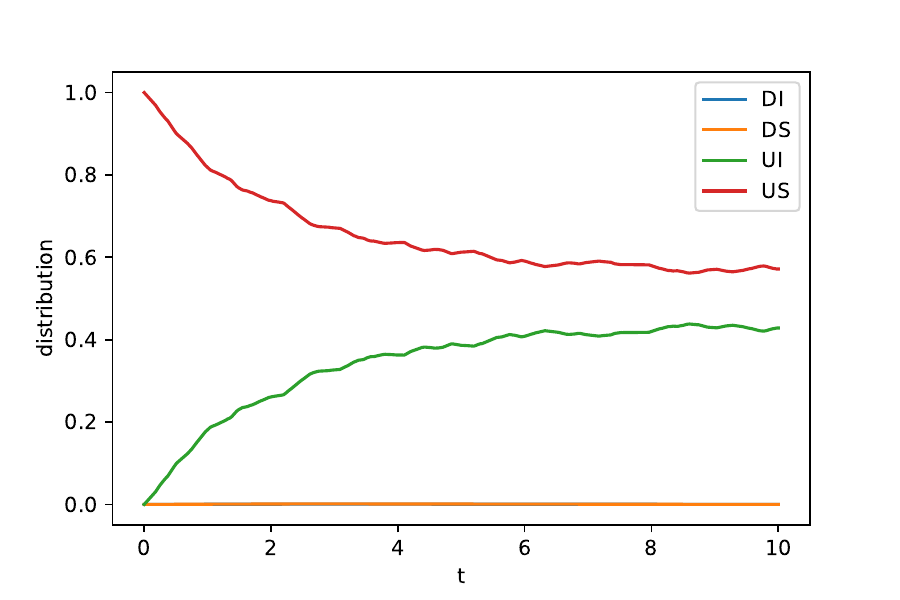}
		\caption{Test distribution 3}
	\end{subfigure}%
	\\
	\begin{subfigure}{.3\textwidth}
		\centering\captionsetup{width=.8\linewidth}
		\includegraphics[width=\textwidth]{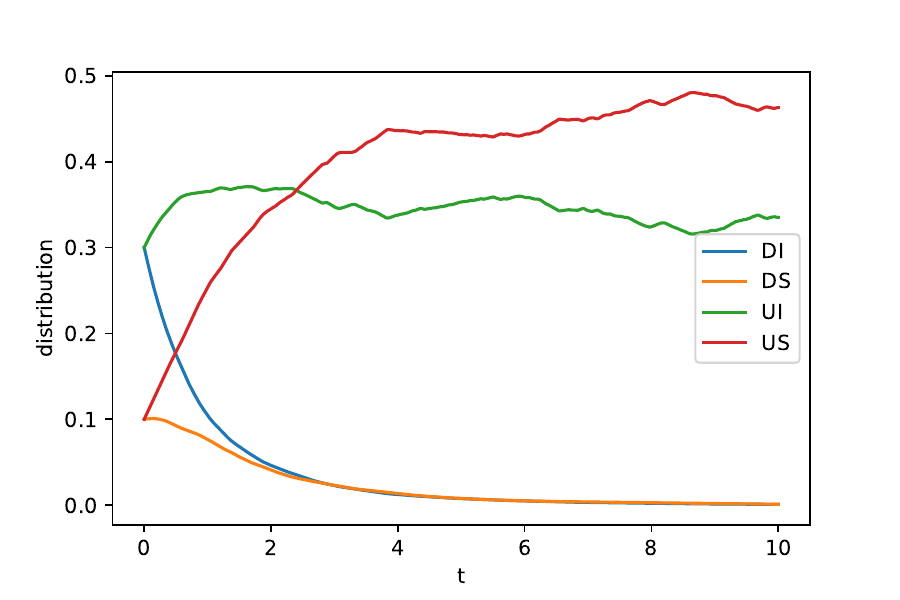}
		\caption{Test distribution 4}
	\end{subfigure}%
	\begin{subfigure}{.3\textwidth}
		\centering\captionsetup{width=.8\linewidth}
		\includegraphics[width=\textwidth]{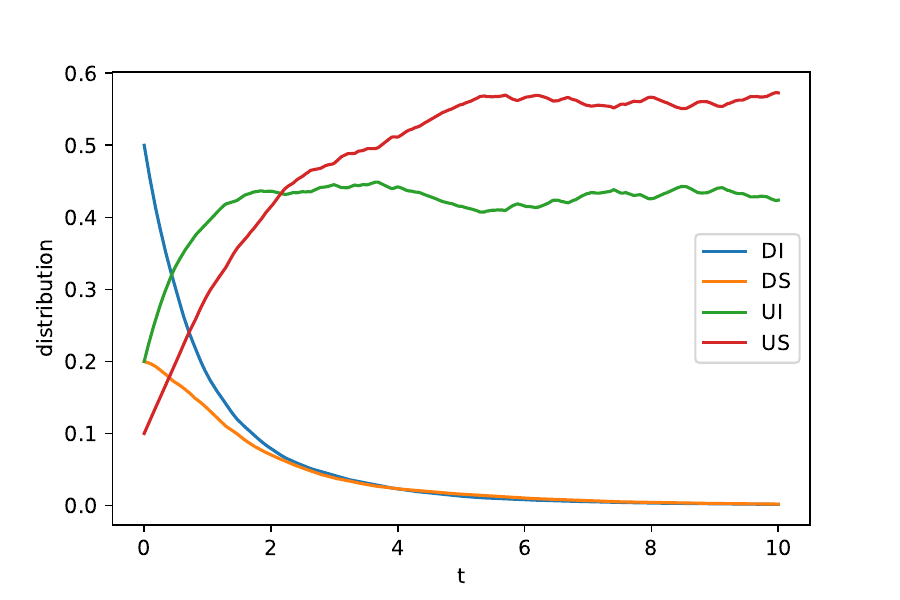}
		\caption{Test distribution 5}
	\end{subfigure}%
	\begin{subfigure}{.3\textwidth}
		\centering\captionsetup{width=.8\linewidth}
		\includegraphics[width=\textwidth]{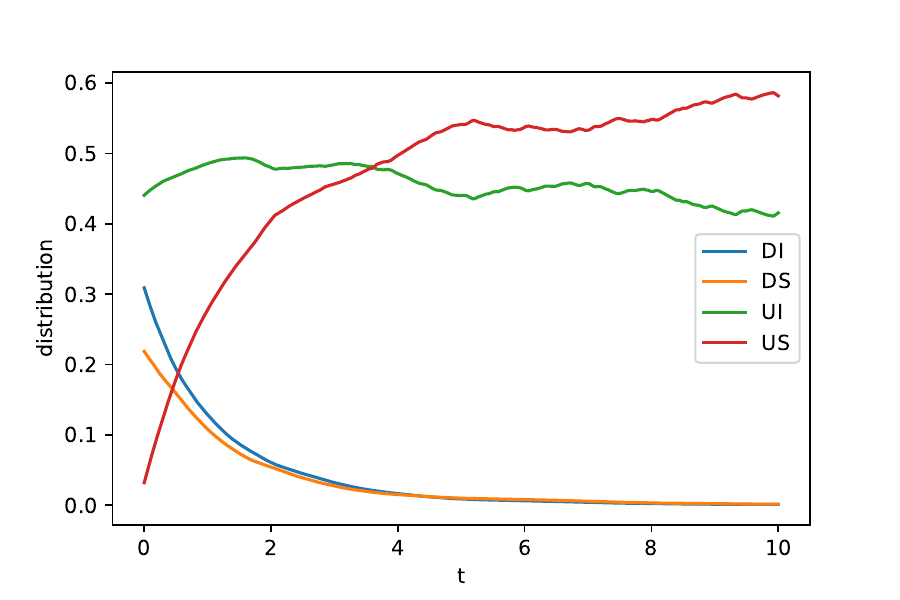}
		\caption{Random training distribution}
	\end{subfigure}
	
	\caption{Cyber security example: Evolution of the distribution in the presence of common noise when applying the control learnt by DDPG on a testing set of five initial distributions and one random initial distribution of the training set. }
	\label{fig:cyber-ddpg-cn-testing}
\end{figure*}

\begin{figure*}
	\centering
	\begin{subfigure}{.48\textwidth}
		\centering\captionsetup{width=.8\linewidth}
		\includegraphics[width=\textwidth]{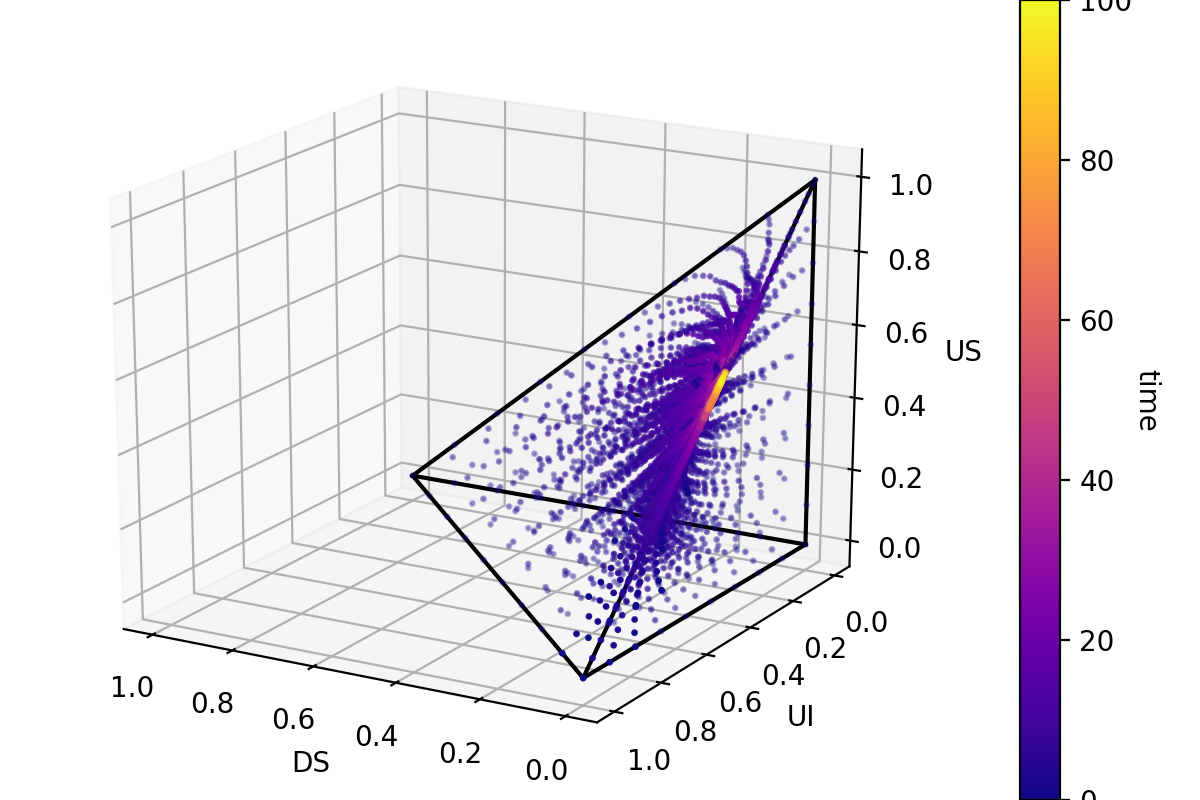}
		\caption{Scatter plot of the distribution trajectories. The color gradient correspond to time.}
	\end{subfigure}%
	\begin{subfigure}{.48\textwidth}
		\centering\captionsetup{width=.8\linewidth}
		\includegraphics[width=\textwidth]{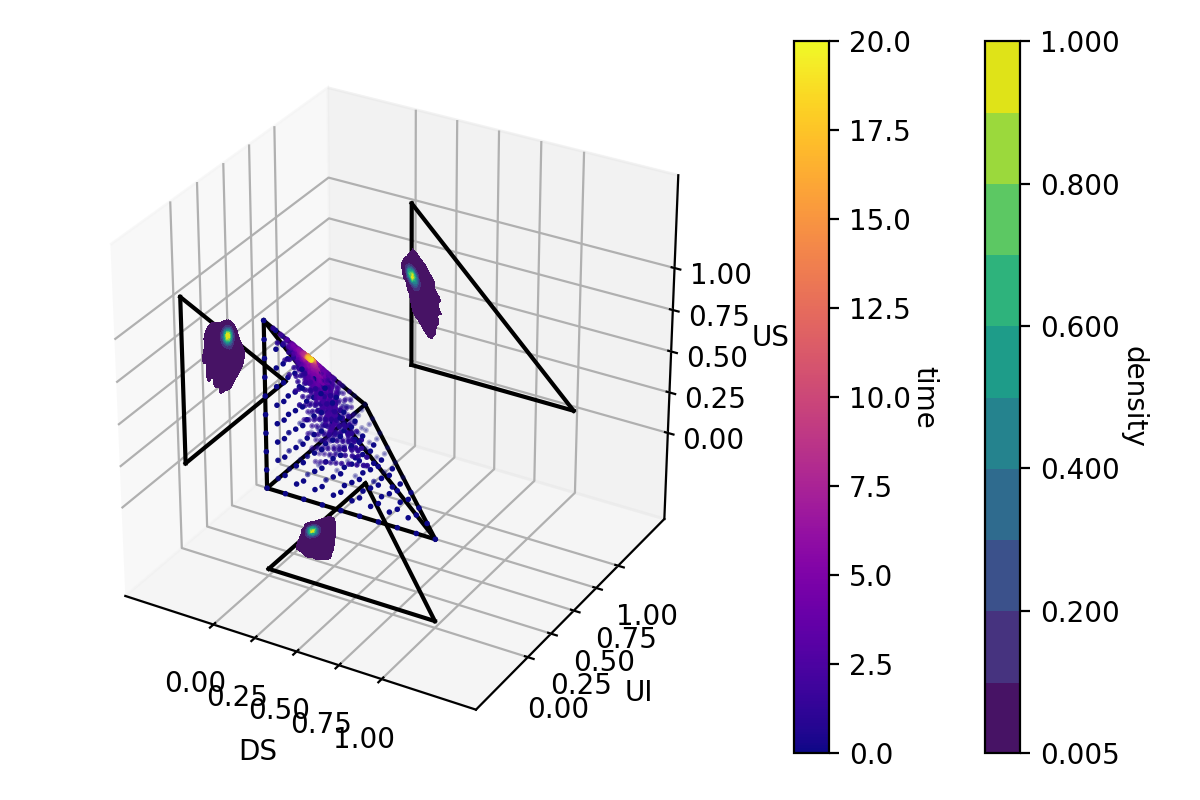}
		\caption{Scatter plot with projection on each dimension of a density estimate using Gaussian kernel density estimation.}
	\end{subfigure}%
	\caption{Cyber security example: Trajectories of the distribution's last three coordinates (DS, UI, US) in the simplex, when subject to common noise and controlled by the control learnt with DDPG. The initial distributions are from a uniform discretization of the simplex with mesh size $0.1$ in each dimension.  }
	\label{fig:cyber-ddpg-3d-simplex}
\end{figure*}

Using the DDPG method described above, we train the neural networks by picking at each episode a random initial distribution $\mu$ and a random sequence of common noises $\epsilon^0$. Fig.~\ref{fig:cyber-ddpg-cn-testing} displays the evolution of the population when using the learnt control starting from five initial distributions of the testing set and one initial distribution of the training set. The testing set of initial distributions is: $\{(0.25,0.25,0.25,0.25),$ $(1, 0, 0, 0),$ $(0, 0, 0, 1),$ $(0.3, 0.1, 0.3, 0.1),$ $(0.5, 0.2, 0.2, 0.1)\}$.   Consistently with the case without common noise, we see that the distribution always evolves towards a configuration in which there is no defended agents, and the proportion of undefended infected and undefended susceptible are roughly $0.43$ and $0.57$ respectively. Due to the common noise, the distribution is not perfectly stable; it oscillates around these values. Figure~\ref{fig:cyber-ddpg-3d-simplex} shows from two perspectives the evolution of the mean field state dynamics when applying the learnt optimal control. The initial distributions are on a uniform grid of the simplex and time steps are distinguished by colors. We see that for any initial distribution, as time increases, the mean field states concentrate around the aforementioned point, which lies on an edge of the simplex. 

\subsection{Example 2: Discrete distribution planning}\,

We now consider an MFC problem in which the goal is to match a target distribution. We take a model with $N_{states} = 10$ states and $3$ actions (left, stay, right): $X = \{1,\dots,10\}$ and $A = \{L,S,R\}$. When at the leftmost state (resp. rightmost), the agents cannot go left (resp. right). The incurred to a representative agent is $1$ if they move (i.e., they use action left or right) plus the $L^2$ distance between the population distribution and a target distribution. Here we chose: $(0, 0, 0.05, 0.1, 0.2, 0.3, 0.2, 0.1, 0.05, 0, 0)$ for the target distribution. There is no idiosyncratic noise. A key point is that, in this setting, except for some specific pairs of initial and target distributions, it is not possible for the population to match the target distribution unless the agents are allowed to randomize their actions at the individual level. So we use $\cP(A)^X$ for the level-1 action space. Hence the action space is naturally continuous and this justifies, here again, the use of the DDPG method.

We present results without and with common noise in the dynamics. In the second case, the common noise is an i.i.d. additive $\epsilon^0_n$ at each time step $n$, with $\epsilon^0_n = -1$ with probability $0.05$, $\epsilon^0_n = 1$ with probability $0.05$ and $\epsilon^0_n = 0$ with probability $0.9$. In both cases, the initial distributions for training are picked randomly as follows. First we pick $x_{min}$ and $x_{max}$ uniformly at random between $1$ and $N_{states}$. This determines a sub-set (with periodicity) $\{x_{min},\dots,x_{max}\}$. For each point in the sub-interval, a value is picked independently and uniformly at random in $[0,1]$. Then the discrete distribution is normalized to have total mass $1$. For numerical reasons, we stop after a finite number of time steps. Here we took $100$ time steps. 

\begin{figure}[h]
	\begin{subfigure}{.9\columnwidth}
		\centering
		\includegraphics[width=\columnwidth]{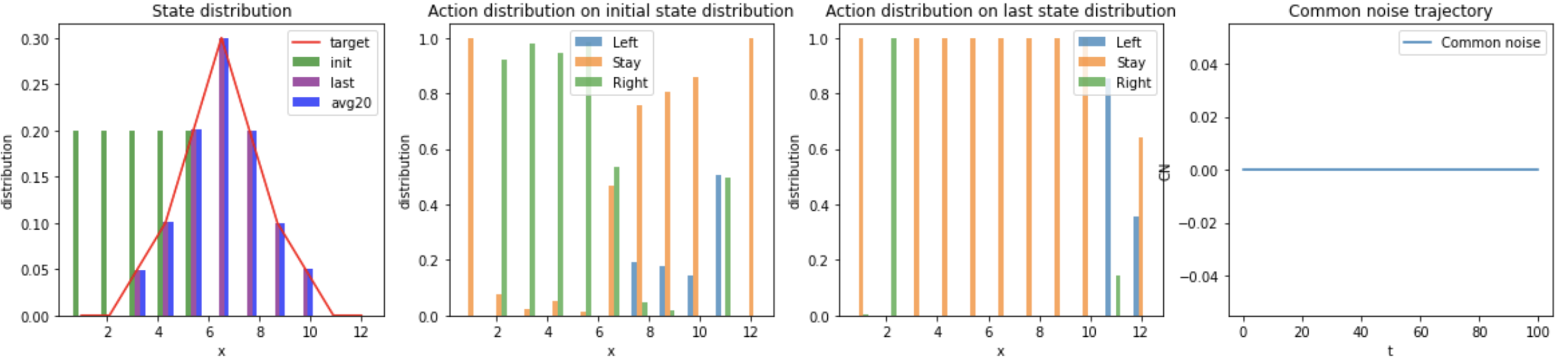}
	\end{subfigure}%
	\\
	\begin{subfigure}{.9\columnwidth}
		\centering
		\includegraphics[width=\columnwidth]{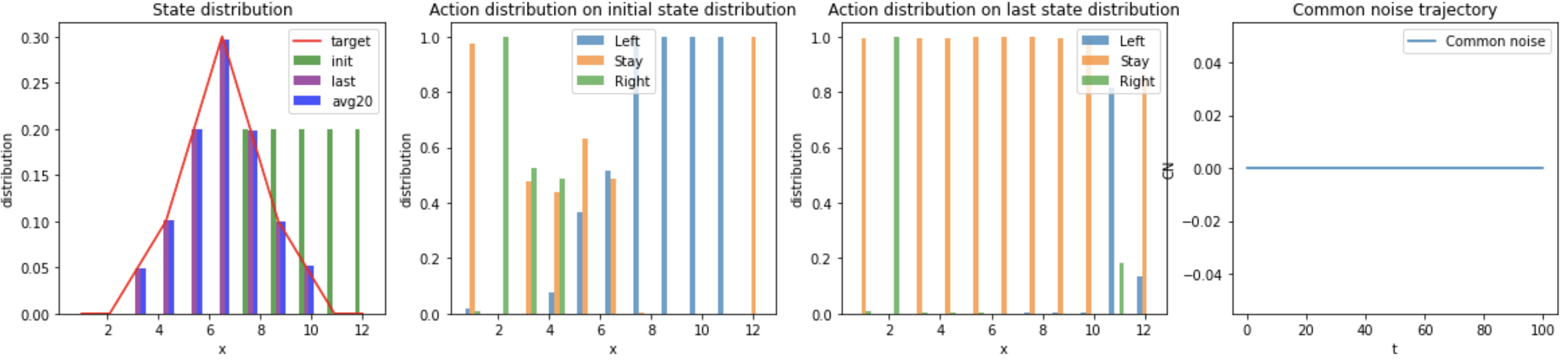}
	\end{subfigure}%
	\\
	\begin{subfigure}{.9\columnwidth}
		\centering
		\includegraphics[width=\columnwidth]{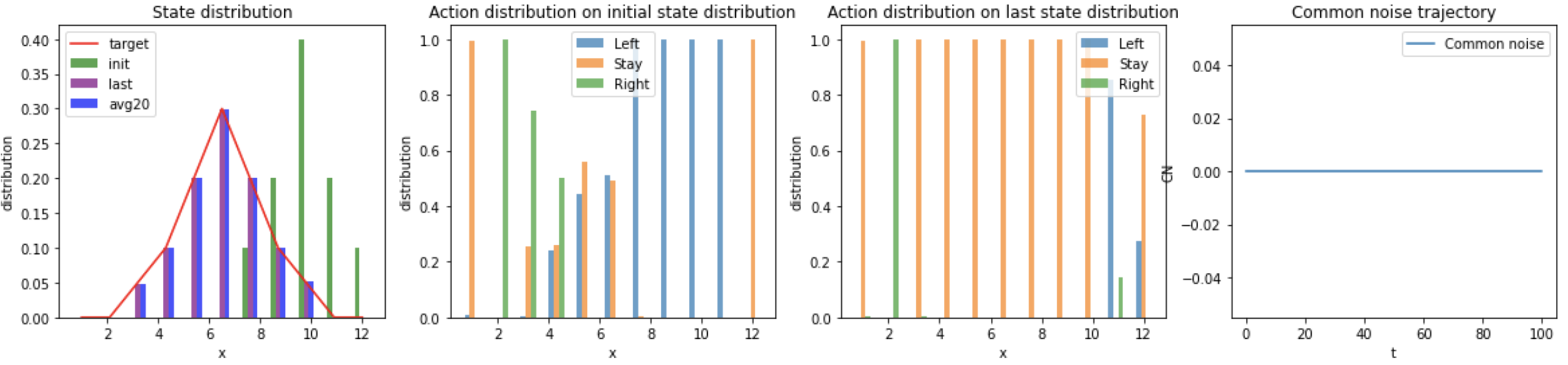}
	\end{subfigure}%
	\\
	\begin{subfigure}{.9\columnwidth}
		\centering
		\includegraphics[width=\columnwidth]{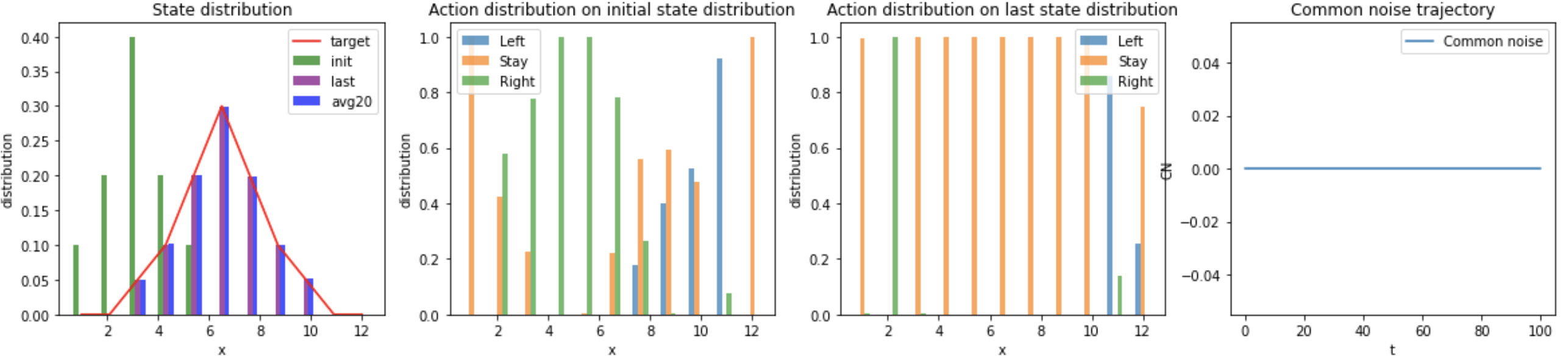}
	\end{subfigure}%
	\\
	\begin{subfigure}{.9\columnwidth}
		\centering
		\includegraphics[width=\columnwidth]{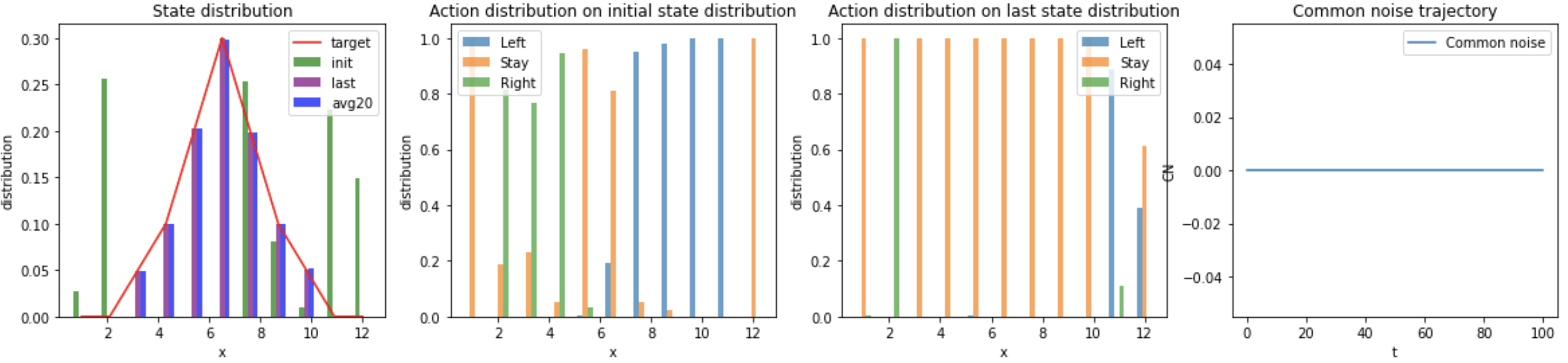}
	\end{subfigure}%
	\\
	\begin{subfigure}{.3\columnwidth}
		\centering
		\includegraphics[width=\columnwidth]{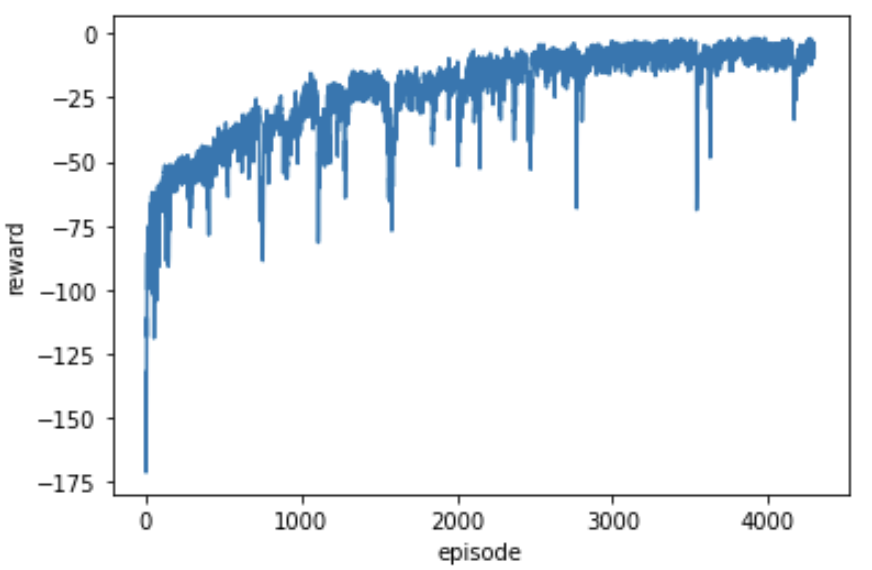}
	\end{subfigure}%
	\caption{Example 2: Discrete distribution planning. Case without common noise; $5$ testing distributions. Column 1: state distribution; columns 2 and 3: action distribution at initial and terminal time; column 4: common noise trajectory (identically $0$ here). Bottom: evolution of the reward during training.}   
 	\label{fig:finiteMFCQ-discreteplanning-DDPG-trainrand-nocn}
\end{figure}

\begin{figure}[h]
	\begin{subfigure}{.9\columnwidth}
		\centering
		\includegraphics[width=\columnwidth]{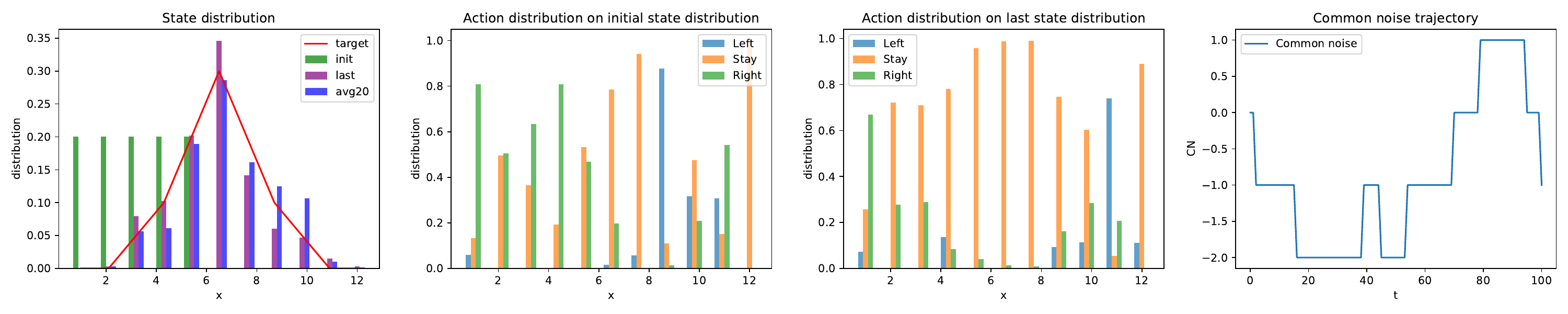}
	\end{subfigure}%
	\\
	\begin{subfigure}{.9\columnwidth}
		\centering
		\includegraphics[width=\columnwidth]{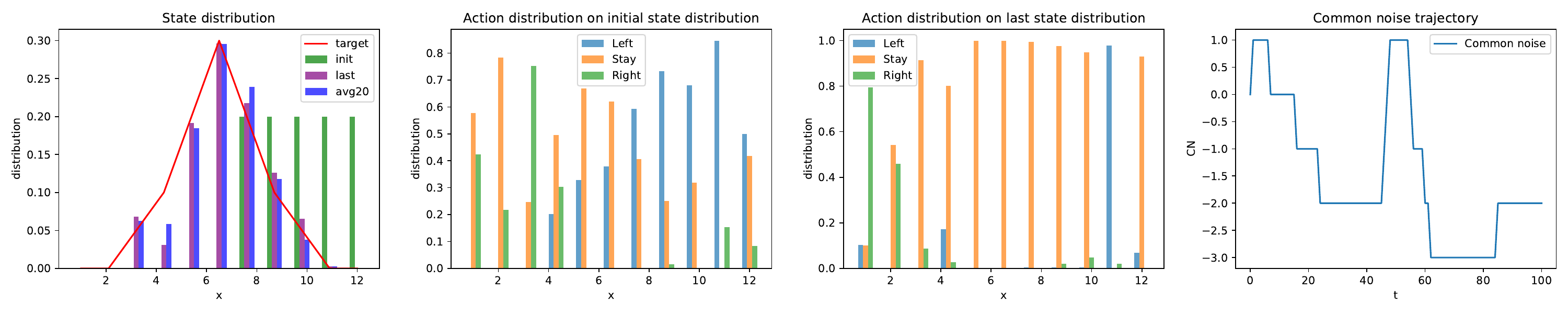}
	\end{subfigure}%
	\\
	\begin{subfigure}{.9\columnwidth}
		\centering
		\includegraphics[width=\columnwidth]{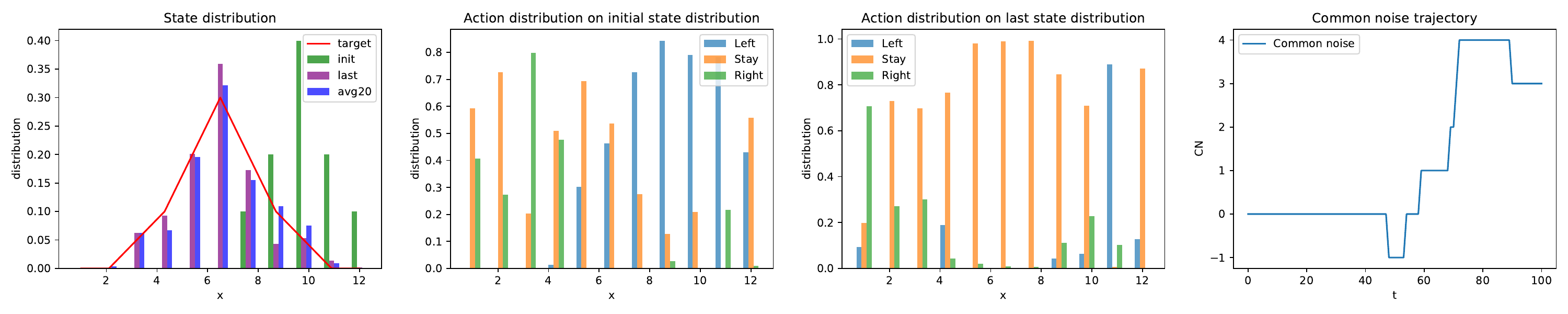}
	\end{subfigure}%
	\\
	\begin{subfigure}{.9\columnwidth}
		\centering
		\includegraphics[width=\columnwidth]{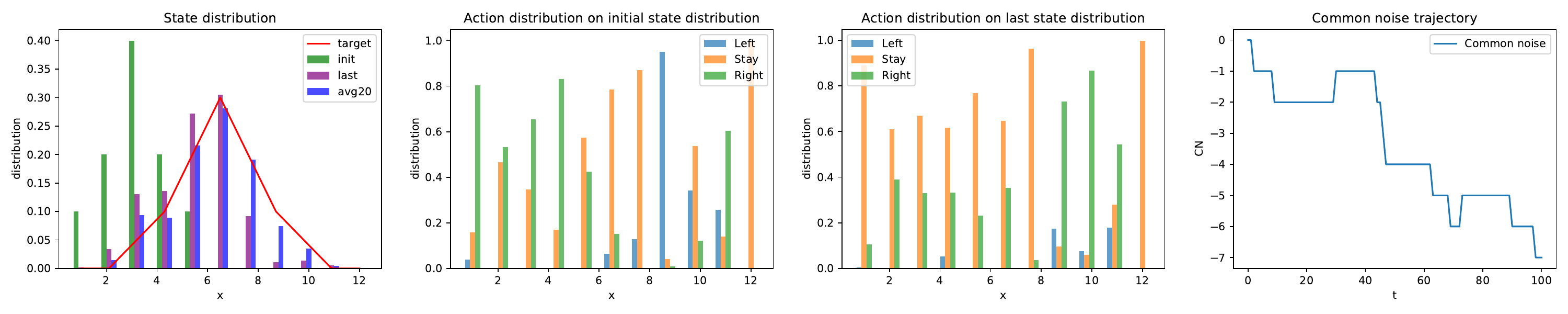}
	\end{subfigure}%
	\\
	\begin{subfigure}{.9\columnwidth}
		\centering
		\includegraphics[width=\columnwidth]{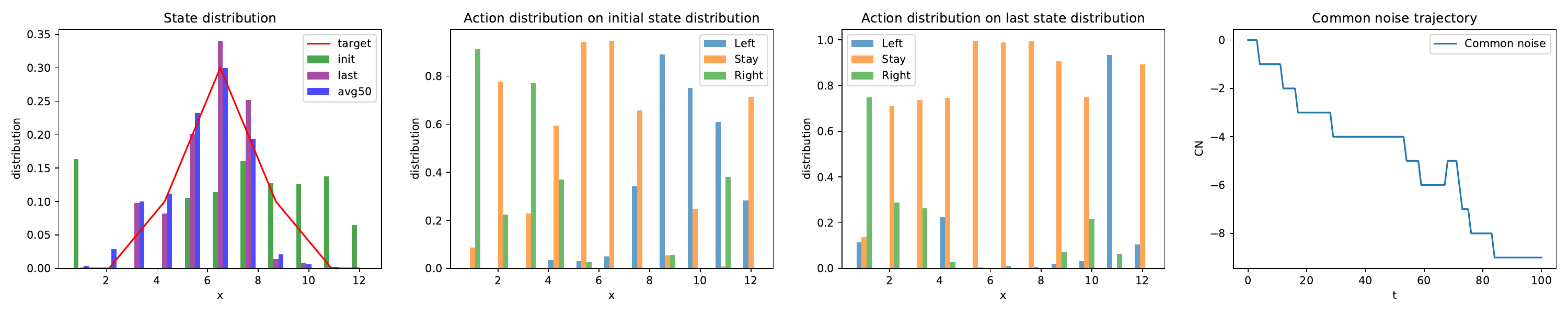}
	\end{subfigure}%
	\\
	\begin{subfigure}{.3\columnwidth}
		\centering
		\includegraphics[width=\columnwidth]{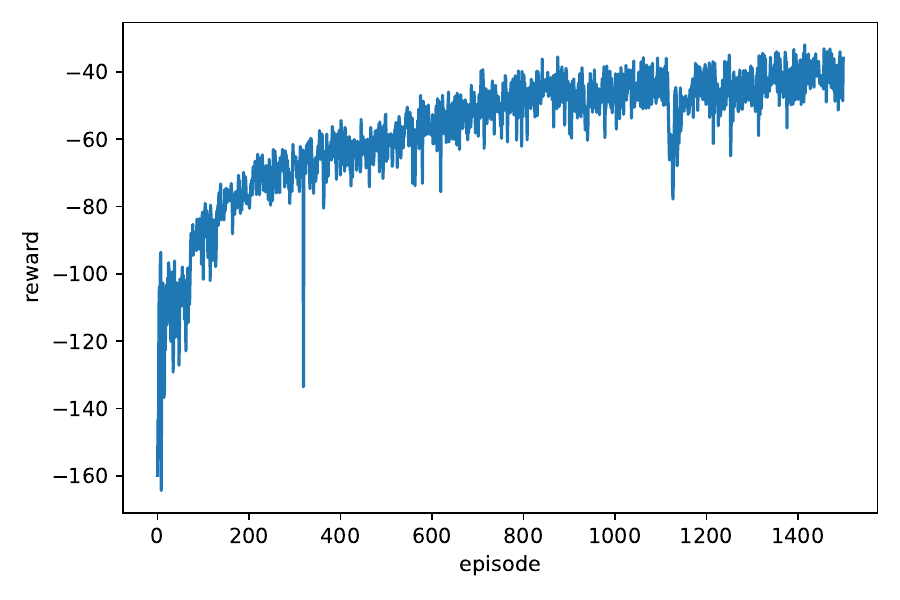}
	\end{subfigure}%
	\caption{Example 2: Discrete distribution planning. Case with common noise; $5$ testing distributions. Column 1: state distribution; columns 2 and 3: action distribution at initial and terminal time; column 4: common noise trajectory (identically $0$ here). Bottom: evolution of the reward during training.}   
 	\label{fig:finiteMFCQ-discreteplanning-DDPG-trainrand-withcn}
\end{figure}

Figures~\ref{fig:finiteMFCQ-discreteplanning-DDPG-trainrand-nocn} and~\ref{fig:finiteMFCQ-discreteplanning-DDPG-trainrand-withcn} present the results obtained without and with common noise, respectively. In each case, the results are for five different testing distributions: four fixed test distributions, as well as one random distribution. The left column displays the state distribution: initial distribution in green, target distribution in red, last distribution in purple, and the average over the last few steps before terminal time in blue. We see that the last distribution is very close to the target one so the learnt control is successful. The two columns in the middle display the control distribution at time $0$ and at terminal time. For each state, the probability of picking each action is represented by a vertical bar, with one color per possible action. We see that, at initial time, the most likely choice is to move to the right (resp. left) for states on the left (resp. right) of the domain. At terminal time, the most likely choice is to stay at the current location, because the target distribution has been reached. The right column displays the trajectory of the common noise that has affected the distribution in the present run (constant equal to $0$ in Figure~\ref{fig:finiteMFCQ-discreteplanning-DDPG-trainrand-nocn}). On Figure~\ref{fig:finiteMFCQ-discreteplanning-DDPG-trainrand-withcn}, we see that even when the common noise is rather strong, the learnt control manages to move the initial distribution close to the target distribution. At the bottom of each figure is displayed the evolution of the training reward (the negative of the MFC cost) along the training iterations (also called episodes) of DDPG.

\subsection{Example 3: Swarm motion} \,

We then turn our attention to a model in continuous state and action spaces. More precisely, we consider a model of swarm motion with aversion to crowded regions introduced in~\cite{MR3698446} (in the context of mean field games). Although many variants are possible, we define the model in the following way in order to have an analytical solution that can be used to assess the convergence of our proposed method. We take the interval $[0,1]$ with periodic boundary condition, i.e. the unit torus $\TT$, as the state space $S$. The action space is $A = \RR$.  
The dynamics of a typical agent is driven by~\eqref{eq:system_dynamics_level_0} with $F(x, a, \mu, e, e^0) = a + e + e^0$. In other words, the central planner chooses the velocity of each agent. The instantaneous reward of a typical agent at location $x$ and using action $a$ while the population's state is $\mu$, is defined as:
$
	f(x,a,\mu) = -\frac{1}{2} |a|^2 + \varphi(x) - \ln(\mu(x)).
$ 
Here, the first term penalizes a large velocity (it can be interpreted as a kind of cost proportional to the kinetic energy of the agent),  $\varphi$ encodes spatial preferences (by giving a lower cost for certain positions in space), and the last term models crowd aversion (it penalizes the fact of being at a location where the density of agents is high). We choose:
$$
	\varphi(x) = - 2 \pi^2 \left[ - \sin(2 \pi x) + |\cos(2 \pi x)|^2 \right] + 2\sin(2 \pi x),
$$
We consider that there is no common noise ($\varepsilon^0_n = 0$ for all $n$), and the idiosyncratic noises $\varepsilon_n$ have a Gaussian distribution. 
We obtain a model which, in continuous time, admits an explicit ergodic solution that we can use as a benchmark. Indeed, in this case the optimal ergodic control is given by $\tilde a(x) = 2 \pi \cos(2 \pi x)$ and the ergodic distribution of the corresponding MKV dynamics has density $\mu(x) = e^{2\sin(2 \pi x)} / \int e^{2\sin(2 \pi x')}dx'$.

The action space being continuous, here again the use of DDPG is justified. To implement this approach, we however need a finite dimensional representation of the distribution in order to pass it to the policy network and the value function network. We replace $\cP(\TT)$ by a finite dimensional simplex  $\cP(\{0,1/N_p,\dots,1-1/N_p, 1\})$ corresponding to a uniform discretization of $\TT$ with $N_p+1$ points. The environment (whose inner working is not known to the learning agent) needs to compute the evolution of the distribution. This evolution can be directly simulated with a deterministic method based, for example, on a finite difference scheme as in~\cite{MR2679575}. However, in practice, it is likely that the environment would not work in this way but would rather correspond to moving forward a large population of agents (e.g., robots). This induces extra approximations. To illustrate that our method can be applied in such situations, here we chose to implement the environment using a probabilistic approach based on Monte Carlo simulations for a large number of particles on $S$. Then, to prepare the input for the Q-function, we project their positions on $\{0,1/N_p,\dots,1-1/N_p, 1\}$ and approximate the mean field distribution by a histogram. We recall that the DDPG method uses this environment as a black-box and, for a given action $\tilde a \in \RR^{N_p}$, can only access the resulting new distribution and the associated reward. 
The actor and critic networks have been implemented using a feedforward fully connected architecture with $2$ hidden layers of width at most $300$ neurons. We used random initial states at each episode, and the noise used on the action is a Gaussian noise with mean $0$ and variance $0.02$. We used Adam optimizer with initial learning rate $0.0001$ and minibatches of size $16$. 

\begin{figure}
	\begin{subfigure}{\textwidth}
		\centering\captionsetup{width=.8\linewidth}
		\includegraphics[width=\textwidth]{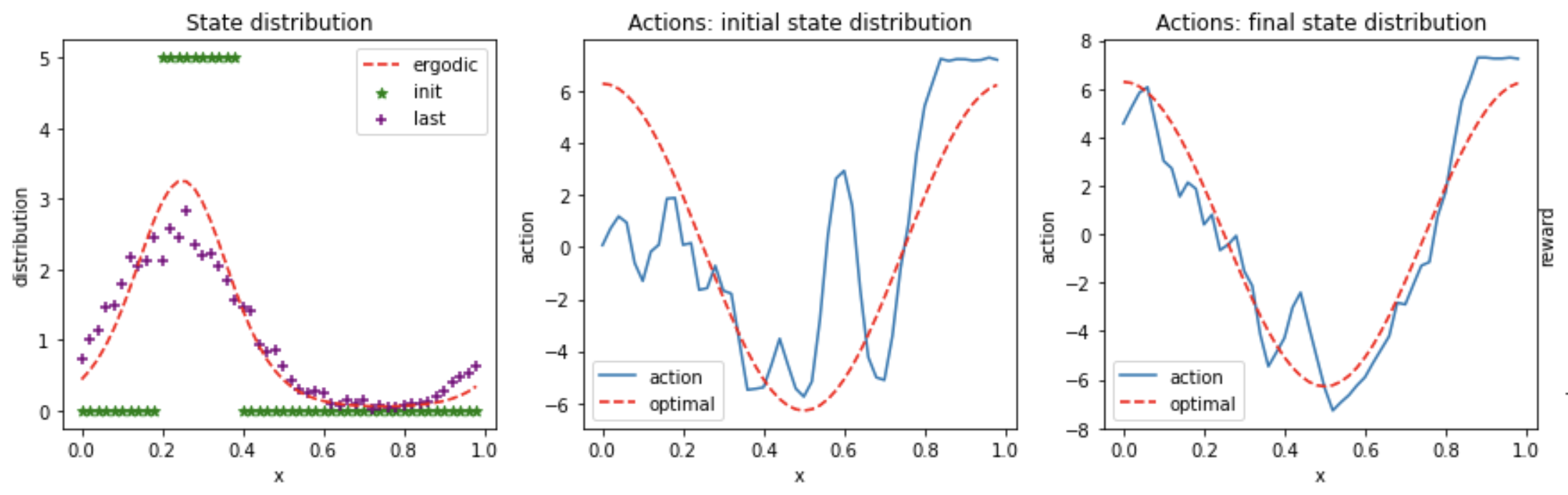}
	\end{subfigure}%
	\\
	\begin{subfigure}{\textwidth}
		\centering\captionsetup{width=.8\linewidth}
		\includegraphics[width=\textwidth]{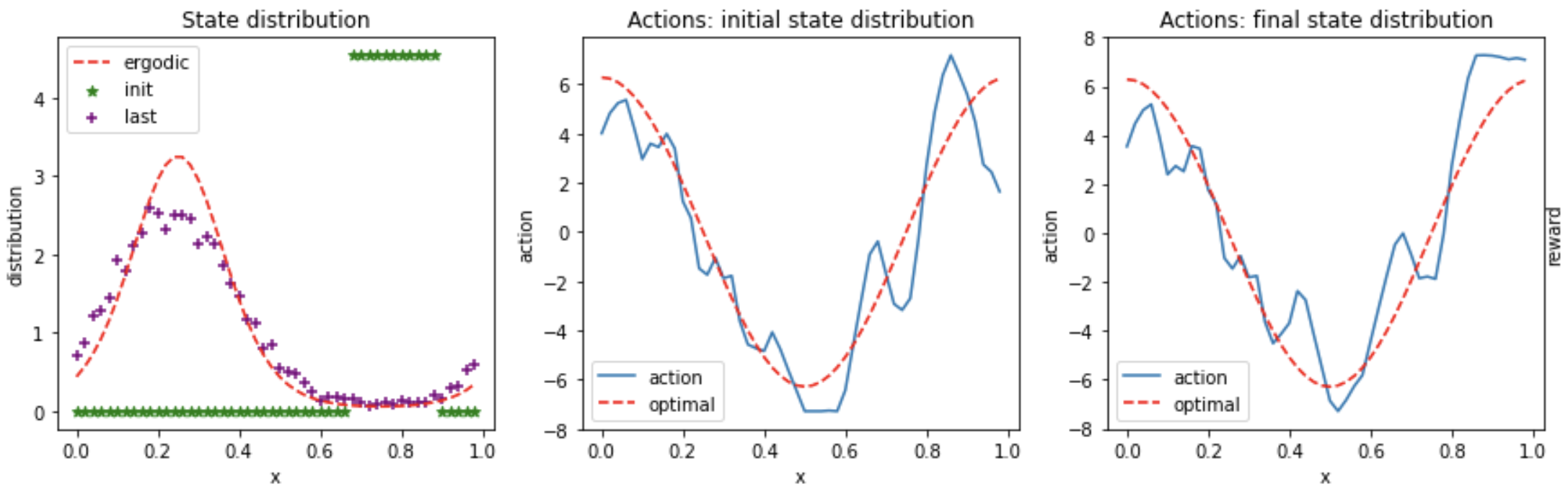}
	\end{subfigure}%
	\\
	\begin{subfigure}{\textwidth}
		\centering\captionsetup{width=.8\linewidth}
		\includegraphics[width=\textwidth]{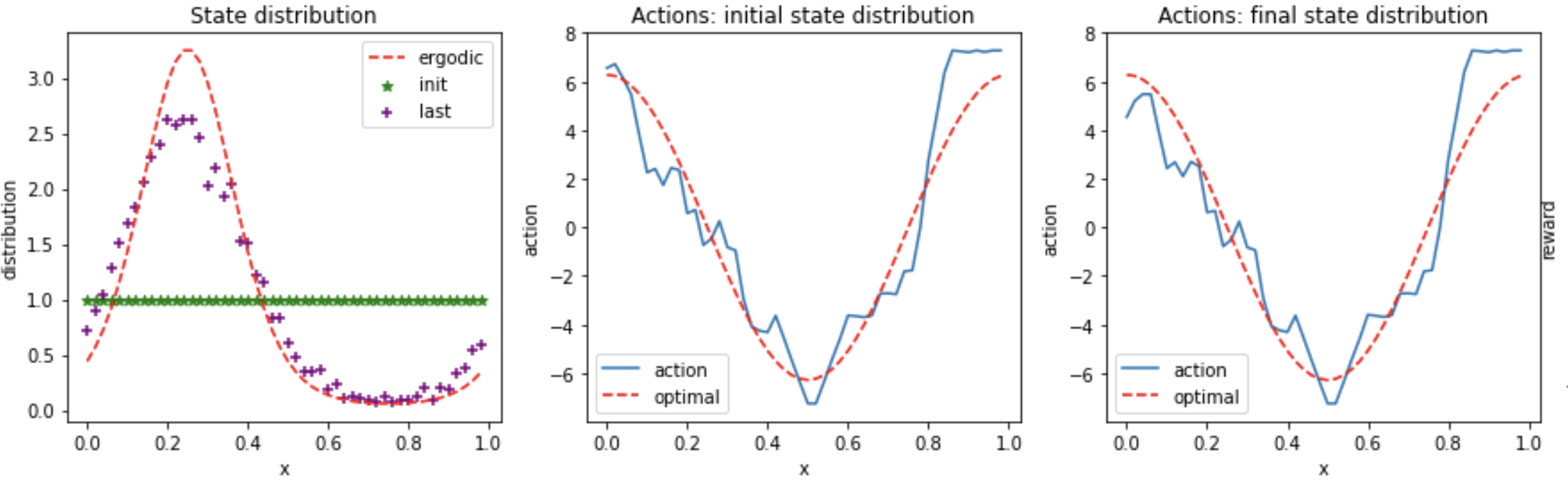}
	\end{subfigure}%
	\\
	\begin{subfigure}{.3\textwidth}
		\centering\captionsetup{width=.8\linewidth}
		\includegraphics[width=\textwidth]{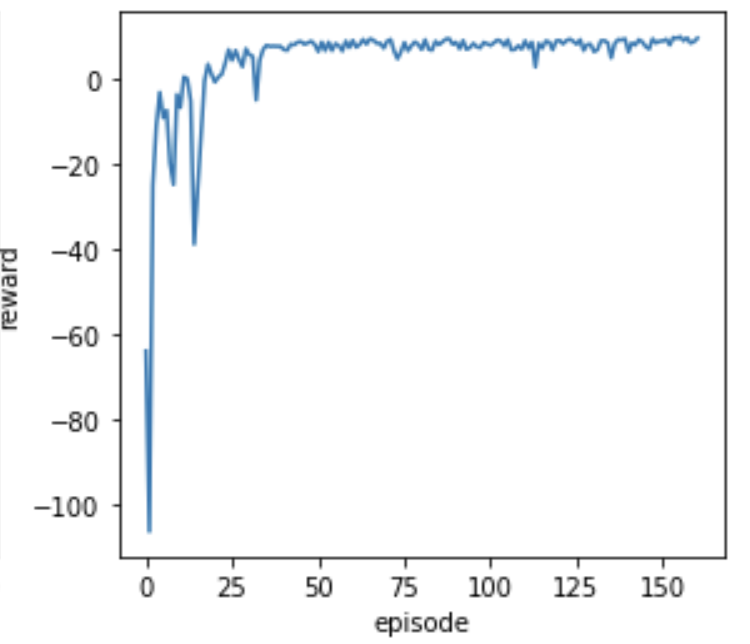}
	\end{subfigure}%
	\caption{Swarm motion: Left: Distribution induced, middle: control learnt for distribution at time $0$, right: control learnt for distribution at terminal time. For three different initial distributions. Red: ergodic analytical solution; green: initial distribution; purple: distribution at terminal time. Bottom: rewards.}
	\label{fig:swarm-motion}
\end{figure}

Figure~\ref{fig:swarm-motion} presents results obtained using this method after $160$ episodes. The system has been trained on initial distributions which are Gaussian with random mean and random variance. As illustrated in the figures (left column), the system has learnt how to drive this type of initial distributions towards the analytical stationary distribution and then how to use an approximation of the stationary optimal control (right column) in order to keep the system in the stationary regime. The middle column displays the learnt control for the initial distribution. It is not expected to match the optimal ergodic control, which should be applied when the distribution has reached the ergodic regime and here it is provided only for the sake of comparison.

\bibliographystyle{apalike}
\small
\bibliography{mfqrl-bib}

\clearpage

\appendix

\section{Auxiliary results for Section~\ref{sec:model_description}}

\begin{lemma}
	\label{lemma:identitly_level_0_action_realization}
	Given a level-0 control process $\fa$, there exists a level-0 action process $\balpha$ which is a realization of $\fa$. Moreover, if another level-0 action process $\balpha'$ is also a realization of $\fa$, then for every $n \geq 0$ and for every bounded Borel measurable function $h: A \to \RR$, we have
	\begin{equation}
		\label{eq:identity_level_0_action_realization}
		\EE \left[ h( \alpha_n' ) \, | \, \cG_n^c \right] = \int_A h(\alpha) \fa_n(d \alpha) = 	\EE \left[ h( \alpha_n ) \, | \, \cG_n^c \right], \qquad \PP-a.s..
	\end{equation}
\end{lemma}

\begin{proof}
By definition, for each $n \geq 0$, $\fa_n$ is of the form:
	$$
		\fa_n(d \alpha)= \kappa_n^{\fa}(\sU, \underline{\vartheta}_{n-1}, \underline{\vartheta}_n^0, \underline{\varepsilon}_{n}, \underline{\varepsilon}_{n}^0 )(d \alpha), \qquad\PP-a.s. \, ,
	$$
for some measurable function $\kappa_n^{\fa}$ on $\Upsilon \times \Theta^{n-1} \times (\Theta^{0})^{ n} \times E^n \times (E^{0})^{n}$ with values in $\cP(A)$. Let us denote by $\xi_n$ the random element $(\sU, \underline{\vartheta}_{n-1}, \underline{\vartheta}_n^0, \underline{\varepsilon}_{n}, \underline{\varepsilon}_{n}^0 )$, by $\rho_A$ the Blackwell-Dubins function (see Lemma~\ref{le:BlackwellDubins}) of the space $A$. Let us set $U_n=h(\vartheta_n)$ and:
	$$
		\alpha_n(\omega) := \rho_A( \kappa_n^\fa(\xi_n(\omega)), U_n(\omega)).
	$$
Then, $\alpha_n$ is $\cG_n^a$-measurable, and because the $\sigma-$field $\cG_n^c = \sigma\{\xi_n\}$ is independent of $U_n$, we have $\cL( \alpha_n | \cG_n^c) = \kappa_n^\fa(\xi_n) = \fa_n$, $\PP$-almost surely. Equality~\eqref{eq:identity_level_0_action_realization} directly follows the definition of a conditional distribution.
\end{proof}

\begin{lemma}
\label{le:policy_value}
For each $n\ge 0$, the law of the random measure $\PP^0_{(X_n,\alpha_n)}$ depends only upon the open-loop policy $\bpi$ as long as $\balpha$ is a realization of the control process generated by $\bpi$. 
\end{lemma}

\begin{proof}
Let $\balpha$ be an action process which is a realization of the control process generated by $\bpi$. For each integer $n\ge 0$, we compute $\EE\bigl[\Phi\bigl(\PP^0_{(X_n,\alpha_n)}\bigr)\bigr]$ for a family of bounded measurable functions $\Phi$ on $\cP(S\times A)$ which generate the Borel $\sigma$-field of $\cP(S\times A)$. For the sake of definiteness we work with functions $\Phi$ of the form:
$$
\Phi(\mu)=\prod_{j=1}^m\int_{S\times A}\varphi_j(x,\alpha)\,\mu(dx, d\alpha)
$$
for a finite set $\varphi_1,\cdots,\varphi_m$ of bounded continuous functions on $S\times A$. We have:
\begin{equation}
\label{fo:big_phi}
\begin{split}
\EE\bigl[\Phi\bigl(\PP^0_{(X_n,\alpha_n)}\bigr)\bigr] 
&=\EE\Bigl[\prod_{j=1}^m\int_{S\times A}\varphi_j(x,\alpha)\,\PP^0_{(X_n,\alpha_n)}(dx,d\alpha)\Bigr]\\
&=\EE\Bigl[\prod_{j=1}^m\EE\bigl[\varphi_j(X_n,\alpha_n)\,| \sigma\{\underline{\vartheta}^0_n,\underline{\varepsilon}^0_n\}\Bigr].
\end{split}
\end{equation}
Now for each $j\in\{1,\cdots,m\}$ we have 
\begin{equation}
\label{fo:phi_j}
\begin{split}
&\EE\bigl[\varphi_j(X_n,\alpha_n)\,| \sigma\{\underline{\vartheta}^0_n,\underline{\varepsilon}^0_n\}\Bigr]\\
&\hskip 35pt
=\int\cdots\int \varphi_j\bigl(X_n(u,\underline{\theta}_{n-1},\underline{e}_n,\underline{\theta}^0_{n},\underline{e}^0_n),
\alpha_n(u,\underline{\theta}_{n-1},\underline{e}_n,\underline{\theta}^0_{n},\underline{e}^0_n,\theta_n)\bigr)\\
&\hskip 125pt
\PP_{\cU}(du)\PP_{\underline{\vartheta}_{n-1}}(d\underline{\theta}_{n-1})\nu^n(d\underline{e}_n)\PP_{\vartheta_{n}}(d\theta_{n})\Bigr|_{\underline{\theta}^0_n=\underline{\vartheta}^0_n,\underline{e}^0_n=\underline{\varepsilon}^0_n}\\
&\hskip 35pt
=\int\cdots\int \Bigl(\int_A \varphi_j\bigl(X_n(u,\underline{\theta}_{n-1},\underline{e}_n,\underline{\theta}^0_{n},\underline{e}^0_n),
\alpha\bigr)\pi_n(d\alpha\,| u,\underline{\theta}_{n-1},\underline{e}_n,\underline{\theta}^0_{n},\underline{e}^0_n)\Bigr)\\
&\hskip 125pt
\PP_{\cU}(du)\PP_{\underline{\vartheta}_{n-1}}(d\underline{\theta}_{n-1})\nu^n(d\underline{e}_n)\Bigr|_{\underline{\theta}^0_n=\underline{\vartheta}^0_n,\underline{e}^0_n=\underline{\varepsilon}^0_n}\\\end{split}
\end{equation}
where we made explicit the dependence of $X_n$ on $\xi_n=(\cU,\underline{\vartheta}_{n-1},\underline{\varepsilon}_n,\underline{\vartheta}^0_{n},\underline{\varepsilon}^0_n)$ and $\alpha_n$ on $(\xi_n,\vartheta_n)$. This shows that the left hand side of \eqref{fo:phi_j}, and hence the left hand side of \eqref{fo:big_phi} only depend upon the action process $\balpha=(\alpha_n)_{n\ge 0}$ through the conditional distribution $\pi_n(d\alpha| u,\underline{\theta}_{n-1},\underline{e}_n,\underline{\theta}^0_{n},\underline{e}^0_n)$. From this we conclude that if two action processes are realizations of control processes generated by the same the open-loop policy $\bpi$, the corresponding random measures $\PP^0_{(X_n,\alpha_n)}$ have the same distribution.
\end{proof}

\begin{remark}
\label{rem:barF-Borel}
	For every $(\mu, \bar a , e^0) \in \bar \Gamma \times E^0$, $\bar F(\mu, \bar a , e^0)$ is defined as a probability measure on $S$ such that for every bounded and Borel measurable function $\phi : S \to \RR$, 
	\begin{equation}
		\label{eq:pushforward_of_F}
		\int_S \bar F(\mu, \bar a, e^0) (dx' ) \phi(x') = \int_{S \times A \times E} \bar a(dx, d\alpha) \nu(de) \phi \Big( F( x, \alpha, \bar a, e, e^0) \Big).
	\end{equation}
	It is straightforward to check that $\bar F$ is Borel measurable. See  for example \cite[Proposition~7.29]{BertsekasShreve} for a proof.
\end{remark}

\begin{lemma}
\label{le:under_A_regu-barF-barf}
Assume~\ref{assumption:basic_assumption_MFC}. It holds:
\begin{itemize}
    \item $\bar F$ is Borel measurable and for every $e^0 \in E^0$, $\bar F(\cdot, \cdot, e^0)$ is continuous in its remaining variables.
    \item $\bar f$ is bounded and lower semi-continuous.
\end{itemize}
\end{lemma}

\begin{proof}
The measurability of $\bar F$ was argued earlier (see after Definition~\ref{def:MFMDP_problem}), so we only argue the continuity for $e^0\in E^0$ fixed. Let $\bigl( (\mu_n,\bar a_n)\bigr)_{n\ge 0}$ be a sequence in $\bar \Gamma$ which converges weakly toward $(\mu,\bar a)$. Since $\mu_n=\text{pr}_1(\bar a_n)$ for each $n\ge 0$, we necessarily have $\mu=\text{pr}_1(\bar a)$ so that $(\mu,\bar a)\in\bar\Gamma$. We pick a continuous bounded function $\phi: S\mapsto \RR$ and we show that:
\begin{equation}
\label{fo:F_bar_n}
\lim_{n\to\infty}\int_S \phi(x')\bar F(\mu_n,\bar a_n,e^0)(dx')= \int_S \phi(x')\bar F(\mu,\bar a,e^0)(dx').
\end{equation}
Using Skorohod's characterization of weak convergence of probability measures, we have the existence of random variables $(Y_n,\beta_n)$ converging $\PP$-almost surely toward some $(Y,\beta)$ and such that $\PP_{(Y_n,\beta_n)}=\bar a_n$ for each $n\ge 0$ and $\PP_{(Y,\beta)}=\bar a$. Consequently, the integral in the left hand side of \eqref{fo:F_bar_n} can be rewritten as:
\begin{equation*}
\begin{split}
\int_S \phi(x')\bar F(\mu_n,\bar a_n,e^0)(dx')& = \int_{S\times A\times E} \phi\bigl(F(x,\alpha,\bar a_n,e,e^0)\bigr)\bar a_n(dx,d\alpha)\nu(de)\\
&=\int_E \Bigl(\EE\bigl[\phi\bigl(F(X_n,\alpha_n,\bar a_n,e,e^0)\bigr)\bigr]
\bigr)\nu(de),
\end{split}
\end{equation*}
which converges toward $\int_E \Bigl(\EE\bigl[\phi\bigl(F(X_n,\alpha_n,\bar a_n,e,e^0)\bigr)\bigr]\bigr)\nu(de)$
by Lebesgue's dominated convergence theorem because $F(\cdot,\cdot,\cdot,e,e^0)$ is continuous and $\phi$ is bounded continuous.

The lower semi-continuity of the one stage cost function $\bar f$  follows from  \cite[Proposition~7.31 (a)]{BertsekasShreve} which only requires the lower semi-continuity of the original one-stage cost function $f$ instead of the full continuity assumption posited in Assumption~\ref{assumption:basic_assumption_MFC}.
\end{proof}

\begin{lemma}%
	\label{lemma:MFMDP_well_definedness_value_function}
	
	Let  $\bar \bpi\in\bar\bPi$. For every $\mu \in \bar S$, let $(\bmu, \bm{\bar a})$ and $(\bmu', \bm{\bar a'})$ be two pairs of state and action processes generated by $(\bar \bpi, \mu)$. Then 
	$
		\EE \left[ \sum_{n \geq 0} \gamma^n \bar f (\mu_n, \bar a_n) \right] = \EE \left[ \sum_{n \geq 0} \gamma^n \bar f (\mu_n', \bar a_n') \right].
	$
\end{lemma}

\begin{proof}
    For any fixed initial $\mu \in \bar S$, by definition of the pair of state and action processes generated by $(\bar \bpi, \mu)$, we show by induction that, for all $n \ge 0$, $\cL(\mu_n )=\cL(\mu'_n)$ and $\cL\bigl((\mu_n, \bar a_n)\bigr)=\cL\bigl( (\mu'_n, \bar a'_n)\bigr)$.
    For $n = 0$, $\mu_0 = \mu_0' = \mu \in \bar S$.  Assume that for some $n \geq 0$, we have $\cL(\mu_n )=\cL(\mu'_n)$, then for every bounded and Borel measurable function $\phi: \bar S \times \bar A \to \RR $
	\begin{equation}
		\begin{split}
			\EE \left[ \phi ( \mu_n, \bar a_n ) \right] 
			& = \EE \left[ \EE\left[ \phi ( \mu_n, \bar a_n ) \, | \, \mu_n \right] \right]
			\\
			& = \EE \left[ \int_{\bar A} \phi( \mu_n, \bar a).  \cL(\bar a_n \, | \, \mu_n)(d \bar a) \right] 
			 \\
			 &= \EE \left[ \int_{\bar A} \phi( \mu_n, \bar a). \bar \pi_n(\mu_n) (d \bar a) \right] 
			 \\
			 &= \EE \left[ \int_{\bar A} \phi( \mu_n', \bar a).  \bar \pi_n(\mu_n') (d \bar a) \right] 
			 = \EE \left[ \EE\left[ \phi ( \mu_n', \bar a_n' ) \, | \, \mu_n' \right] \right]
			 =  \EE \left[ \phi \big( \mu_n', \bar a_n' \big) \right].
		\end{split}
	\end{equation}
	So $\cL\bigl((\mu_n, \bar a_n)\bigr)=\cL\bigl( (\mu'_n, \bar a'_n)\bigr)$. 
	Since $\varepsilon_{n+1}^0$ is independent of $(\mu_n, \bar a_n)$ and $(\mu_n', \bar a_n')$,  $\cL\bigl((\mu_n, \bar a_n, \varepsilon_{n+1}^0)\bigr) = \cL\bigl( (\mu_n', \bar a_n', \varepsilon_{n+1}^0)\bigr)$, which implies that the law of $
	    \mu_{n+1} = \bar F (\mu_n, \bar a_n, \varepsilon_{n+1}^0)$   is equal to the law of  $\bar F (\mu_n', \bar a_n', \varepsilon_{n+1}^0) = \mu_{n+1}'$. 
	Hence the conclusion.
\end{proof}

\section{Proofs for Section~\ref{sec:relation-models-closedloop-barpi}}

\begin{proof}[Proof of Lemma~\ref{lemma:dynamics_of_conditional_distribution}]    
 Let us denote by $ \zeta_{n+1}$ the right hand side of \eqref{eq:formula_dynamics_in_lemma_MFC_MFMDP}, and let  $\phi: S \to \RR$, $h_n: ( \Theta^0\times E^0 )^n \to \RR$ and $\psi_{n+1}: E^0 \to \RR$ be arbitrary bounded Borel measurable functions. We have:
    \begin{equation*}
    \begin{split}
    &\EE \left[  \psi_{n+1}(\varepsilon_{n+1}^0)  h_n(\underline \vartheta_n^0, \underline \varepsilon_n^0  )  \int_{S} \phi(x') \zeta_{n+1}(d x') \right]
    \\
     & \hskip 35pt
     = \EE \left[  \psi_{n+1}(\varepsilon_{n+1}^0)  h_n( \underline \vartheta_n^0, \underline \varepsilon_n^0  ) \int_{S} \phi(x') \Big(  \bar F( \PP^0_{X_n}, \PP^0_{(X_n, \alpha_n)}, \varepsilon_{n+1}^0 ) \Big) (d x') \right]
    \\
     & \hskip 35pt
 = \EE \left[  \psi_{n+1}(\varepsilon_{n+1}^0) h_n( \underline \vartheta_n^0, \underline \varepsilon_n^0  ) \int_{S \times A \times E} \PP^0_{(X_n, \alpha_n)}(d x, d\alpha) \nu(d e) \phi \Big( F( x, \alpha, \PP^0_{(X_n, \alpha_n)}, e, \varepsilon_{n+1}^0) \Big) \right]
     \\
     & \hskip 35pt
 = \int_{E \times E^0} \nu( d e)  \nu^0(d e^0) \psi_{n+1}(e^0) \EE \left[ h_n( \underline \vartheta_n^0, \underline \varepsilon_n^0  )  \int_{S \times A} \PP^0_{(X_n, \alpha_n)}(d x, d\alpha) \phi \Big( F( x, \alpha, \PP^0_{(X_n, \alpha_n)}, e, e^0) \Big) \right]
     \\
     & \hskip 35pt
 = \EE \left[ \psi_{n+1}( \varepsilon_{n+1}^0)  \EE \left[ h_n(\underline \vartheta_n^0, \underline \varepsilon_n^0 )
 \phi \Big( F( X_n, \alpha_n, \PP^0_{(X_n, \alpha_n)}, \varepsilon_{n+1}, \varepsilon_{n+1}^0 \big) \Big) \, \Big| \, \cF_n^0, \varepsilon_{n+1}, \varepsilon_{n+1}^0 \right]  \right]
     \\
     & \hskip 35pt
 = \EE \left[  \psi_{n+1}( \varepsilon_{n+1}^0)  h_n( \underline \vartheta_n^0, \underline \varepsilon_n^0 ) \phi (X_{n+1}) \right],
    \end{split}
    \end{equation*}
    where the first equality is by definition of $\zeta_{n+1}$, the second equality is by the definition of $\bar F$ in terms of the system function $F$ of the original MFC, the third equality is by the fact that $\varepsilon^0_{n+1}$ is independent of all the other random quantities, the fourth equality is by definition of the conditional probability $\PP^0_{(X_n, \alpha_n)}$, and  the last equality is by the tower property of conditional expectation, the fact that $X_{n+1} = F( X_n, \alpha_n, \PP^0_{(X_n, \alpha_n)}, \varepsilon_{n+1}, \varepsilon_{n+1}^0 \big)$ and the fact that $(\varepsilon_{n+1}, \varepsilon_{n+1}^0) $ is independent of $(X_n, \alpha_n) $ and $\PP^0_{(X_n, \alpha_n)}$ is measurable with respect to  $\cF_n^0 = \sigma\{ \underline \vartheta_n^0, \underline \varepsilon_n^0 \}$. This shows that  $\zeta_{n+1}=\PP^0_{X_{n+1}}$.
\end{proof}

\begin{proof}[Proof of Lemma~\ref{lemma:idenfiy_conditional_joint_dist_with_kernels}] 
 We prove equation~\eqref{eq:identify_conditional_joint_distribution} by induction. Equality \eqref{eq:identify_conditional_joint_distribution} holds for $n=0$ by the assumptions on $(\bzeta, \bm{\bar \eta} )$. For the sake of an argument by induction, let us assume that \eqref{eq:identify_conditional_joint_distribution} holds for some $n\ge 0$. We first show that $\cL(\zeta_{n+1}) = \cL(\PP^0_{X_{n+1}})$. Since $\varepsilon_{n+1}^0$ is independent of $\cF_n^0$, for every bounded Borel measurable function $\psi: \bar S \to \RR$, it holds:
    \begin{equation*}
        \begin{split}
            \EE[ \psi( \zeta_{n+1} ) ] &= \EE \left[ \EE \left[ \psi \Big( \bar F( \zeta_n, \bar \eta_n, \varepsilon_{n+1}^0 ) \Big) \, \Big| \, \zeta_n, \bar \eta_n \right] \right]
            \\
            & = \EE \left[ \int_{\bar S} \psi( \mu) P( \zeta_n , \bar \eta_n)(d \mu)   \right]
            \\
            & = \EE \left[ \int_{\bar S} \psi( \mu) P\Big( \PP^0_{X_n} , \PP^0_{(X_n, \alpha_n)} \Big) (d \mu)  \right]
            \\
            & = \EE \left[ \psi \big( \PP^0_{X_{n+1}} \big) \right], 
        \end{split}
    \end{equation*}
where the second equality is by definition of $P$ in~\eqref{def:MFMDP_transition_kernel_from_system_func}, the third equality is by the induction hypothesis, and the last equality is due to Lemma~\ref{lemma:dynamics_of_conditional_distribution}. So $\cL(\zeta_{n+1}) = \cL(\PP^0_{X_{n+1}})$. We then consider the joint law of $(\zeta_{n+1}, \bar \eta_{n+1})$.  By the assumption that $\cL( \bar \eta_n  | \zeta_n ) = \kappa_n( \zeta_n )$, we have that $(\zeta_{n+1}, \bar \eta_{n+1})$ and  $\big(\PP^0_{X_{n+1}}, \PP^0_{(X_{n+1}, \alpha_{n+1})} \big)$ share  the same regular version $\kappa_n$ of the conditional probability. We conclude that \eqref{eq:identify_conditional_joint_distribution} holds for $n+1$ instead of $n$.
 \end{proof}

\clearpage 

\section{Proofs for Section~\ref{sec:back-to-MFC-2}}
\label{app:proofs-back-MFC-2}
\begin{proof}[Proof of Lemma~\ref{lemma:construct_open_loop}]
We prove this statement by showing that there exist $\bpi \in \bPi^{\tinyol}$ and   an  open-loop action process $\balpha$ generated by $\bpi$ such that:
    $
        J^{ \bpi}= J^{ \balpha } = J^{\tilde{\bpi}}.
    $
	Let $\mu_0 \in \cP(S)$ and let $(\bX, \balpha)$ be a pair of state and action processes generated by $(\tilde{\bpi}, \mu_0)$. Let $\bfa$ be the $\cP(A)$-valued process given by:
	$$
	    \fa_n= \tilde\pi_n( X_n, \PP^0_{X_n}, \vartheta_n^0) , \qquad n \geq 0.
	$$
We recall that $\Xi_n, \xi_n$ and $\sU$ are defined in \S~\ref{subsec:open-closed-policies}. Since $\fa$ is adapted to $\GG^c$, it is an admissible level-0 control process, and for every $n \geq 0$, there exists a Borel measurable function $\pi_n: \Xi_n \to A$ satisfying
$$
     \fa_n
     = \pi_n(  \sU, (\vartheta_k, \vartheta_k^0, \varepsilon_{k+1}, \vartheta_{k+1}^0)_{k = 0, \ldots, n-1}, \vartheta_n^0 )
     = \pi_n(\xi_n), \qquad \PP-a.s. \, .
$$
Moreover, since $\balpha$ is generated by $\tilde\bpi$, it is adapted to $\GG^a$ and satisfies
$$  
    \cL( \alpha_n \, | \, \cG_n^c) = \tilde\pi_n( X_n, \PP^0_{X_n}, \vartheta_n^0) = \pi_n( \xi_n), \qquad  \PP-a.s. \qquad n \geq 0.
$$
So $\balpha$ can be viewed as an open-loop action process generated by $\bpi$. 
Meanwhile, the state process $\bX$ constructed by equation~\eqref{eq:system_dynamics_level_0} 
is also a state process associated with $(\balpha, \mu_0)$ (see Definition~\ref{def:state_process_from_control_process}). Therefore, by definition of the value function associated to an open-loop policy $\bpi$, we have:
$$
   J^{\tilde\bpi}(\mu_0)
   = J^{\balpha}(\mu_0) 
   = \EE \left[ \sum_{n\geq 0} \gamma^n f \big(X_n, \alpha_n, \PP^0_{(X_n, \alpha_n)} \big) \right].
$$
\end{proof}

\section{Disintegration of kernels}
\label{sec:disintegration}

Given a probability measure $P$ on a measurable space $(C, \cC)$ and a kernel $K$ from $(C,\cC)$ to $(D,\cD)$,  the composition of measure $P$ and kernel $K$, denoted by $P \measprod K$, is defined as a measure on the product space $(C\times D, \cC\otimes\cD)$ such that for every non-negative measurable function $f : C \times D \to \RR_+$,
$$
	(P \measprod K) f = \int_C P(dx) \int_D f(x, y) K(x, dy).
$$
Similarly, for a probability kernel $\mu: G \to \cP(C)$ and a probability kernel $K: G \times C \to \cP(D)$, the composition of kernels $\mu$ and $K$, denoted by $\mu \measprod K$, is defined as a kernel from $(G, \cG)$ to $(C \times D, \cC \otimes \cD)$ such that for every $s \in G$ and for every non-negative measurable function $f : C \times D \to \RR_+$,
$$
	( \mu \measprod K )(s) f = \int_{C} \mu(s, dx) \int_{D} f(x,y) K(s, x, dy).
$$

\section{Details on Q-learning, Section~\ref{se:Q_learning}}
\label{app:details-Qlearning}

\subsection{Proof of Theorem~\ref{thm:main-cv-tabular-pure}}
\begin{proof}[Proof of Theorem~\ref{thm:main-cv-tabular-pure}]

Recall that we denote by $\check J^*$ and $\check Q^*$ respectively the state value function and the state-action value function of the projected MFC problem defined by~\eqref{eq:generic-MFC-fctmeasure-reward-proj}--\eqref{eq:generic-MFC-fctmeasure-dyn-proj}.

	 We first note that, for every $(\mu, \tilde a) \in \frak{S} \times \cA$,
	\begin{align*}
	    \left| \check Q_{N_{\mathrm{epi}}} \big( \proj_{\check{\frak{S}}}(\mu), \tilde a \big)  - \tilde Q^* \big( \mu, \tilde a \big) \right| 
	    \leq  & \left| \check Q_{N_{\mathrm{epi}}} \big( \proj_{\check{\frak{S}}}(\mu), \tilde a \big) - \check Q^* \big( \proj_{\check{\frak{S}}}(\mu), \tilde a \big) \right| 
	    \\
	    &\quad + \left|  \check Q^* \big( \proj_{\check{\frak{S}}}(\mu), \tilde a \big) - \tilde Q^* \big(  \proj_{\check{\frak{S}}}(\mu), \tilde a \big) \right| 
	    \\
	   &  \quad + \left| \tilde Q^* \big(  \proj_{\check{\frak{S}}}(\mu), \tilde a \big) - \tilde Q^* \big( \mu, \tilde a \big)  \right|.
	\end{align*}
	We then split the proof into three steps, which consist in bounding from above each term in the right hand side.

\vskip 6pt
	\textbf{Step 1.} We first analyze the difference between $\check Q_{N_{\mathrm{epi}}}$ and $\check Q^*$. 
			This comes from standard convergence results on Q-learning for finite state-action spaces.  
			More precisely, under Assumptions~\ref{hyp:bdd-smooth-data} and \ref{hyp:covering-time}, with our choice of learning rates, and given that $N_{\mathrm{epi}}$ is of order~\eqref{eq:LB-Nepi-tabular}, we can apply Theorem~4 and Corollary~34 in~\cite{MR2247972} for asynchronous Q-learning and polynomial learning rates, and we obtain that, with probability at least $1-\delta$, 
			$$
				\|\check Q_{N_{\mathrm{epi}}} - \check Q^* \|_{\infty}= \sup_{ (\check \mu, \tilde a) \in \check{\frak{S}} \times \cA} \Big| \check Q_{N_{\mathrm{epi}}}(\check \mu, \tilde a) - \check Q^*(\check \mu, \tilde a) \Big| \le \varepsilon.
			$$

\textbf{Step 2. } We then turn our attention to the difference between $\check Q^*$ and $\tilde Q^*$. 
			The analysis amounts to say that the projection on $\check{\frak{S}}$ realized at each step does not perturb too much the value function.  
			Recall that for some given common noise $\varepsilon^0$, the operator $\check \Phi^{\varepsilon^0}: \check{\frak{S}} \times \cA \to \check{\frak{S}}$ is given by $\check \Phi^{\varepsilon^0}( \check \mu, \tilde a) = \proj_{\check{\frak{S}}} \circ \bar F ( \check \mu, \check \mu \measprod \tilde a, \varepsilon^0)$.
			Likewise, we denote the transition dynamic with $\bar F$ by a function $\Phi^{\varepsilon^0}: \frak{S} \times \cA \to \frak{S}$ such that:
			$$
				\Phi^{\varepsilon^0} (\mu, \tilde a) = \bar F( \mu, \mu \measprod \tilde a, \varepsilon^0), \qquad \forall (\mu, \tilde a) \in \frak{S} \times \cA.
			$$ 
			Let us start by noting that, for every $(\check \mu, \tilde a) \in \check{\frak{S}} \times \cA$,
			\begin{align*}
				& \left| \check Q^*(\check \mu, \tilde a) -  \tilde Q^*(\check \mu, \tilde a) \right| 
				\\
				& \le  \gamma    \EE \Bigg[ \Big| \check J^*(\check \Phi^{\varepsilon^0}(\check \mu, \tilde a)) - \bar J^*( \Phi^{\varepsilon^0}(\check \mu, \tilde a) ) \Big| \Bigg]
				\\
				& \le   \gamma   \EE \Bigg[  \Big|  \check J^*( \check \Phi^{\varepsilon^0}(\check \mu, \tilde a) ) - \bar J^*( \check \Phi^{\varepsilon^0}(\check \mu, \tilde a)) \Big|  + \Big| \bar J^* (\check \Phi^{\varepsilon^0}(\check \mu, \tilde a)) - \bar J^*( \Phi^{\varepsilon^0}(\check \mu, \tilde a)) \Big|  \Bigg]
				\\
				& \le \gamma   \EE \Bigg[  \Big|  \inf_{\tilde a' \in \cA} \check Q^* \left(\check \Phi^{\varepsilon^0}(\check \mu, \tilde a), \tilde a' \right) - \inf_{\tilde a' \in \cA } \tilde Q^* \left(\check \Phi^{\varepsilon^0}(\check \mu, \tilde a), \tilde a' \right) \Big| \Bigg]  +  \gamma L_{\bar J ^*} \EE \left[  \left\Vert \check \Phi^{\varepsilon^0}(\check \mu, \tilde a) -  \Phi^{\varepsilon^0}(\check \mu, \tilde a) \right\Vert_{d_{\frak{S}}} \right],
			\end{align*}
			where the last inequality holds by Lipschitz continuity of $\bar J^*$ on $\frak{S}$, see Assumption~\ref{hyp:smooth-V}. 
			
			The second term in the last inequality can be bounded using the simplex discretization properties and Assumption~\ref{hyp:bdd-smooth-data}: 
			\begin{align*}
			   \EE \left[  \| \check \Phi^{\varepsilon^0}(\check \mu, \tilde a) - \Phi^{\varepsilon^0}( \check \mu, \tilde a) \|_{d_{\frak{S}}}  \right]
 				& = \EE_{\varepsilon^0_1}\left[ \| \proj_{\check{\frak{S}}} \circ \bar F(\check \mu, \tilde a, \varepsilon^0) - \bar F(\check \mu, \tilde a, \varepsilon^0)  \|_{d_{\frak{S}}} \right] \leq \varepsilon_{\frak{S}}.
			\end{align*}
			 
    		For the first term, let $\check \mu' = \check \Phi^{\varepsilon^0}(\check \mu, \tilde a) \in \check{\frak{S}}$ to alleviate the notation, and let us consider $\tilde a^*_1 \in \cA$ and $\tilde a^*_2 \in \cA$ satisfying:
    		$$
    		    \check Q^*(\check \mu', \tilde a^*_1) = \inf_{\tilde a' \in \cA} \check Q^* \left(\check \mu', \tilde a' \right)
    		\qquad
    		\text{ and }
    	    \qquad
    	        \tilde Q^*(\check \mu', \tilde a^*_2) = \inf_{\tilde a' \in \cA } \tilde Q^* \left(\check \mu', \tilde a' \right).
    	    $$ 
    	   The existence of $\tilde a^*_1$ and $\tilde a^*_2$ is guaranteed respectively by finiteness of  $\frak{S} \times \cA$. 
    	   
    	   We observe that
    	    \begin{align*} 
	        & \check Q^*(\check \mu', \tilde a_1^*) - \tilde Q^*( \check \mu', \tilde a_2^*) 
    	        \\
    	        & = \Big( \check Q^*(\check \mu', \tilde a_1^*) - \check Q^*( \check \mu',  \tilde a_2^* )  \Big) + \Big( \check Q^*( \check \mu',  \tilde a_2^* ) - \tilde Q^*( \check \mu', \tilde a_2^* ) \Big) 
    	        \\
    	        & \leq 0 +  \sup_{ (\check \mu, \tilde a) \in \check{\frak{S}} \times \cA} \left| (\check Q^* - \tilde Q^*)(\check \mu, \tilde a) \right| %
    	        \\
    	        & \leq  \| \check Q^* - \tilde Q^* \|_{\infty} .
    	    \end{align*}
    	   On the other hand, 
		 \begin{align*}
	       \check Q^*( \check \mu', \tilde a_1^*) -  \tilde Q^*(\check \mu', \tilde a_2^*)
    	        = - \Big( \tilde Q^*( \check \mu' ,  \tilde a_2^* ) - \tilde Q^*( \check \mu', \tilde a_1^*) \Big) - \Big( \tilde Q^*( \check \mu' , \tilde a_1^*) - \check Q^*( \check \mu', \tilde a_1^* ) \Big)
    	        \geq - \| \check Q^* - \tilde Q^* \|_{\infty}.
    	    \end{align*}
	    
Combining the above bounds yields that for every $(\check \mu, \tilde a) \in \check{\frak{S}} \times \cA$,
			$$
			    \left| \check Q^*(\check \mu, \tilde a) - \tilde Q^*(\check \mu, \tilde a) \right| \leq \gamma \left( \| \check Q^* - \tilde Q^* \|_{\infty}   \right) + \gamma L_{\bar J^*} \varepsilon_{\frak{S}}.  
			$$
			Consequently,
			$$
			\|\check Q^* - \tilde Q^* \|_\infty
			\le
			\frac{\gamma }{1-\gamma} L_{\bar J^*} \varepsilon_{\frak{S}}.
			$$

\textbf{Step 3. } Last, we look at the difference between $\tilde Q^*( \proj_{\check{\frak{S}}}(\mu), \tilde a  )$ and $\tilde Q^*(\mu, \tilde a)$. 
			For every $\mu \in \frak{S}$ and $\tilde a \in \cA$, letting $\check \mu = \proj_{\check{\frak{S}}} (\mu)$  to alleviate the notation, we have $\| \check \mu - \mu \|_{d_{\frak{S}}} \leq \varepsilon_{\frak{S}}$. We obtain
			\begin{align*}
				& \left| \tilde Q^* (\check \mu, \tilde a) -  \tilde Q^*(\mu, \tilde a) \right|
				\\
				& \le \Bigg|  \tilde f(\check \mu, \tilde a)  -  \tilde f(\mu, \tilde a) \Bigg| +  \gamma  \EE \Bigg[  \left| \inf_{\tilde a' \in \cA} \tilde Q^* ( \Phi(\check\mu, \tilde a ), \tilde a') - \inf_{\tilde a' \in \cA}  \tilde Q^*(\Phi(\mu, \tilde a), \tilde a') \right| \Bigg]
				\\
				& \le L_{\tilde f} \| \check \mu - \mu \|_{d_{\frak{S}}} 
				 +  \gamma  \EE \left[  \left| \bar J^* ( \bar F ( \check \mu, \tilde a, \varepsilon^0 ) - \bar J^*( \bar F (\check \mu, \tilde a, \varepsilon^0 ) ) \right| \right] 
				\\
				& \le L_{\tilde f} \varepsilon_{\frak{S}}  + \gamma L_{\bar J^*} \EE \left[ \| \bar F(\check \mu, \tilde a, \varepsilon^0 ) - \bar F( \mu, \tilde a, \varepsilon^0) \|_{d_{\frak{S}}} \right]
				\\
				& \le  ( L_{\tilde f} + \gamma L_{\bar J^*} L_{\bar F} )\varepsilon_{\frak{S}}  ,
			\end{align*}
			where we used the Lipschitz continuity of $\tilde f, \bar J^*, \bar F$ and the assumption on $\check{\frak{S}}$, see Assumptions~\ref{hyp:bdd-smooth-data}, \ref{hyp:smooth-V} and the simplex discretization properties.

\end{proof}

\subsection{DDPG algorithm}

In Algorithm~\ref{algo:DDPG-MFC}, we describe the DDPG method for MFC with our notation.

\begin{algorithm}[h!]
\DontPrintSemicolon
\KwData{A number of episodes $N_{\mathrm{epi}}$; a length $T$ for each episode; a minibatch size $N_{\mathrm{batch}}$; a learning rate $\tau$.}
\KwResult{A strategy function for central planner represented by the target network $\pi'_{\omega' }$.}
\Begin{
	Randomly initialize parameters $\theta$ and $\omega$ for critic network $Q_{\theta}$ and actor network $\pi_{\omega}$\;
	Initialize  $\theta' \leftarrow \theta$ and $\omega' \leftarrow \omega$ for target networks $Q'_{\theta'}$ and $\pi'_{\omega'}$ \;
  \For{ $k =0, 1, \dots N_{\mathrm{epi}}-1$}{
    Initial distribution $\check\mu_0$ \;
    Initialize replay buffer $R_{\mathrm{buffer}}$ \;
    \For{$n=0, 1, \dots T-1$}{
    Select an action $\bar a_n = \pi_{\omega}(\check\mu_n) + \epsilon^a_{n+1} \in \RR^{N_p}$, where $\epsilon^a_{n+1}$ is the exploration noise \;
    Execute $\bar a_n$, observe cost $c_n = \tilde f_n(\check\mu_n, \bar a_n)$ and  $\check\mu_{n+1}$\;
    Store transition $(\check\mu_n, \bar a_n, c_n, \check\mu_{n+1})$ in $R_{\mathrm{buffer}}$ \;
    Sample a random minibatch of $N_{\mathrm{batch}}$ transitions $(\check\mu_i, \bar a_i, c_i, \check\mu_{i+1})$ from $R_{\mathrm{buffer}}$ \;
    Set $y_i = c_i + \gamma Q'_{\theta'}(\check\mu_{i+1}, \pi'_{\omega'}(\check\mu_{i+1}))$, for $i=1\dots,N_{\mathrm{batch}}$ \;
    Update the critic by minimizing the loss: $L(\theta)=\frac{1}{N_{\mathrm{batch}}}\sum \limits_i (y_i - Q_{\theta}(\check\mu_i, \bar a_i))^2$ \;
    Update the actor policy using the sampled policy gradient $\nabla_{\omega} J$: 
    $$
    	\nabla_{\omega} J(\omega) \approx
	\frac{1}{N_{\mathrm{batch}}}\sum \limits_i \nabla_{\bar a} Q_{\theta}(\check\mu_i, \pi_{\omega}(\check\mu_i)) \nabla_{\omega} \pi_{\omega} (\check\mu_i)
  $$
    Update target networks: $ \theta' \leftarrow \tau \theta + (1-\tau)\theta'$ and 
    $\omega' \leftarrow \tau \omega + (1-\tau)\omega'$\;
    }
    }
  \KwRet{$\pi'_{\omega'}$}
  }
\caption{DDPG for MFC}
\label{algo:DDPG-MFC}
\end{algorithm}

\end{document}